\documentclass[11pt, a4paper]{article}
\usepackage{authblk}
 
\usepackage{lettrine}
\usepackage{mathpazo}
\usepackage[margin=0.8in]{geometry}

\usepackage[utf8]{inputenc} 
\usepackage[T1]{fontenc}    
\usepackage{hyperref}       
\usepackage{url}            
\usepackage{booktabs}       
\usepackage{amsfonts}       
\usepackage{nicefrac}       
\usepackage{microtype}      
\usepackage{lipsum}
\usepackage{amsmath}
\usepackage{tabularx}
\usepackage{fancyhdr}
\usepackage{setspace}
\usepackage{graphicx}
\usepackage{multirow}
\usepackage{multicol}
\usepackage[tight]{subfigure}
\usepackage{graphicx}%
\usepackage{multirow}%
\usepackage{amsmath,amssymb,amsfonts}%
\usepackage{mathrsfs}
\usepackage{amsthm}%
\usepackage{mathrsfs}%
\usepackage[title]{appendix}%
\usepackage{xcolor}%
\usepackage{textcomp}%
\usepackage{manyfoot}%
\usepackage{booktabs}%
\usepackage{algorithm}%
\usepackage{algorithmicx}%
\usepackage{algpseudocode}%
\usepackage{listings}%
\usepackage{comment}

\newcommand{\bea}{\begin{eqnarray}}
\newcommand{\eea}{\end{eqnarray}}
\newcommand{\beq}{\begin{equation}}
\newcommand{\eeq}{\end{equation}}
\newcommand{\bdm}{\begin{displaymath}}
\newcommand{\edm}{\end{displaymath}}

\usepackage{url}            
\usepackage{booktabs}       
\usepackage{amsfonts}       
\usepackage{nicefrac}       
\usepackage{microtype}      
\usepackage{lipsum}
\usepackage{amsmath}
\usepackage{tabularx}
\usepackage{fancyhdr}
\usepackage{setspace}
\usepackage{graphicx}
\usepackage{multirow}
\usepackage{multicol}
\usepackage[tight]{subfigure}
\usepackage{lscape} 
\usepackage[nolist]{acronym}
\usepackage{pifont}
\usepackage{longtable}
\usepackage{makecell}
\usepackage{soul}

\usepackage{xcolor}

\newcommand{\bfx}{\mathbf{x}}

\definecolor{myOrange}{RGB}{253,174,97}
\definecolor{myYellow}{RGB}{255,255,191}
\definecolor{myGreen}{RGB}{161,218,180}
\definecolor{myLtBlue}{RGB}{44,127,184}
\definecolor{myDkBlue}{RGB}{37,52,148}
\newcommand{\gmk}{\textcolor{myLtBlue}{\ding{110}}}
\newcommand{\omk}{\textcolor{myYellow}{\ding{110}}}
\newcommand{\rmk}{\textcolor{myOrange}{\ding{110}}}
\newcommand{\bmk}{\textcolor{myDkBlue}{\ding{110}}}
\newcommand{\ymk}{\textcolor{myGreen}{\ding{110}}}
\newcommand{\nmk}{$\cdot$}

\title{Analytical benchmark problems and methodological framework for the assessment and comparison of\\
multifidelity optimization methods}

\author[1,2,3*]{\bf L. Mainini}
\author[4]{\bf A. Serani}
\author[1,5,6]{\bf H. Pehlivan}
\author[1,3]{\bf F. Di Fiore}
\author[7]{\bf M. P. Rumpfkeil}
\author[8]{\bf E. Minisci}
\author[9]{\bf D. Quagliarella}
\author[10]{\bf S. Yildiz}
\author[4,11]{\bf S. Ficini}
\author[4]{\bf R. Pellegrini}
\author[12]{\bf A. Thelen}
\author[13]{\bf D. Bryson}
\author[10]{\bf M. Nikbay}
\author[4]{\bf M. Diez}
\author[13]{\bf P. Beran}

\affil[1]{Department of Aeronautics, Imperial College London, London, UK}
\affil[2]{{Department of Aeronautics and Astronautics}, {Massachusetts Institute of Technology}, {{Cambridge, MA}, {USA}}}
\affil[3]{{Department Mechanical and Aerospace Engineering}, {Politecnico di Torino}, {{Turin}, {Italy}}}
\affil[4]{{National Research Council-Institute of Marine Engineering}, {{Rome}, {Italy}}}
\affil[5]{{Department of Naval Architecture and Ocean Engineering}, {Istanbul Technical University}, {{Istanbul}, {Turkey}}}
\affil[6]{{LHEEA}, {\'Ecole Centrale de Nantes}, {{Nantes}, {France}}}
\affil[7]{{Department of Mechanical and Aerospace Engineering}, {University of Dayton}, {{Ohio}, {USA}}}
\affil[8]{{Mechanical and Aerospace Engineering}, {University of Strathclyde}, {{Glasgow}, {UK}}}
\affil[9]{{Italian Aerospace Research Centre}, {{Capua}, {Italy}}}
\affil[10]{{Department of Aeronautics and Astronautics}, {Istanbul Technical University}, {{Istanbul}, {Turkey}}}
\affil[11]{{Roma Tre University}, {{Roma}, {Italy}}}
\affil[12]{{Applied Mechanics Division}, {University of Dayton Research Institute}, {{Ohio}, {USA}}}
\affil[13]{{Wright-Patterson Air Force Base}, {Air Force Research Laboratory}, {{Ohio}, {USA}}}

\affil[*]{\textit{Corresponding author}: l.mainini@imperial.ac.uk; andrea.serani@cnr.it; h.pehlivansolak@imperial.ac.uk; f.difiore@imperial.ac.uk; markus.rumpfkeil@udayton.edu; edmondo.minisci@strath.ac.uk; d.quagliarella@cira.it; yildizsih@itu.edu.tr; simone.ficini@uniroma3.it; riccardo.pellegrini@cnr.it; andrew.s.thelen@nasa.gov; dean.bryson@us.af.mil; nikbay@itu.edu.tr; matteo.diez@cnr.it; philip.beran@us.af.mil}

\date{}

\begin{document}

\maketitle

\begin{abstract}
As engineering systems increase in complexity and performance demands intensify, Multidisciplinary Design Optimization (MDO) methodologies are becoming essential for integrating models from multiple disciplines to optimize complex multi-physics systems. Within this context, major challenges remain in selecting appropriate disciplinary fidelity levels, and how to couple them effectively. Multifidelity methods offer a promising path forward by strategically combining information sources of varying fidelity - whether computational or experimental - to enable efficient and scalable design exploration and optimization. Despite the development of numerous multifidelity methods, their comparative performance remains difficult to assess due to the absence of standardized benchmark frameworks capable of evaluating performance across diverse optimization tasks. To address this gap, this paper introduces a comprehensive benchmarking framework that includes: (i) a suite of analytical benchmark optimization problems designed to stress-test and validate multifidelity methods; (ii) a set of assessment metrics for quantifying and comparing performance over measurable objectives; and (iii) the classification, evaluation, and comparison  of several families of multifidelity optimization methods and frameworks using the proposed benchmarks to identify their respective strengths and weaknesses in real-world scenarios. The proposed benchmark problems are analytically defined functions carefully selected to capture mathematical challenges commonly encountered in real-world applications, including high dimensionality, multimodality, discontinuities, and noise. Their closed-form nature ensures computational efficiency, high reproducibility, and a clear separation of algorithmic behavior from numerical artifacts. The accompanying performance metrics support the systematic evaluation of multifidelity methods, measuring both optimization effectiveness and global approximation accuracy. By providing a rigorous, reproducible, and accessible benchmarking framework, this work aims to enable the broader community to understand, compare, and advance multifidelity optimization methods for complex problems in science and engineering.
\end{abstract}


\section{Introduction}\label{s:Introduction}
Innovative product requirements are evolving rapidly, reflecting the technological advances in many engineering disciplines. The accelerating nature of this change is driven by the growing demand of product performance which, in turn, leads to higher systems complexity and cost. To meet emerging requirements and fast development paces, design processes are required to (i) thoroughly and accurately explore design spaces of large size (high dimensional, broad intervals spanned by the design variables \cite{hwang2019large, 2018-JMeST-Li_etal, 2016-JCS-Liu_etal}, (ii) leverage potentially complex physical interactions to unlock unprecedented performance and avoid deleterious physical couplings that may greatly increase product cost when leading to late discovery of defects \cite{sarkar2019multifidelity, perdikaris2017nonlinear, ng2014multifidelity}. 

The international community widely acknowledges that substantial performance and cost improvements can be achieved by capturing, modelling, and coupling multiple disciplinary aspects earlier in the design process \cite{li2022machine, khatouri2022metamodeling, 2022-CMAME-Guo_etal, 2019-IISE-Song_etal}. This philosophy underpins the framework of multidisciplinary design optimization (MDO), which seeks to systematically integrate knowledge from various engineering domains to optimize the design of complex multi-physics systems. However, there is currently no universally established mathematical scheme to uniquely determine which disciplines to include, the appropriate level of coupling, or the required model fidelity to accurately capture the critical physical phenomena for a given system \cite{de2007state, vanaret2017consequences}. These modeling decisions still heavily rely on expert judgment and domain-specific experience. While this approach may be effective for systems with strong historical precedent, it often breaks down when applied to novel or radically innovative vehicles and technologies, where the relevant interactions and constraints may be poorly understood or entirely unprecedented. 

In this regard, one of the long-term challenges of multidisciplinary design optimization is to enable the use of models at increased level of fidelity earlier in the design exercise, to capture physical aspects that are critical and that could constrain or unlock particular novel solutions \cite{foumani2023multi, 2021-CMAME-Zhang_etal, 2019-IJCFD-Serani_etal}. Indeed, it has been observed and demonstrated that exclusively relying on low fidelity models for the exploration and search of large design spaces may lead to the identification of solutions that are infeasible of significantly sub-optimal because the physics is not sufficiently modeled or resolved \cite{he2021multi, lucia2004reduced}. Nevertheless, simply replacing these models with higher-fidelity models during the optimization is often not a practical strategy, because of the higher computational cost associated with these more informative representations.

Multifidelity methods offer the conceptual approach to efficiently optimize products by judiciously using a limited number of high-fidelity analyses while leveraging the information provided by low-fidelity methods \cite{kennedy2001bayesian, park2017remarks, 2021-KBS-Liu_etal, 2021-JMLR-Moss_etal, 2020-ASOC-Li_etal, 2020-SMO-Yi_etal}. This paper considers multifidelity approaches to be part of a larger class of methods that manipulate a set of information sources to accelerate the overall computational tasks of design exploration and search. These information sources could, in principle, represent the response of a system by using either computational approaches (i.e., mathematical descriptions and the concomitant numerical analyses) and/or non-computational approaches (e.g., physical experiments, analytical solutions, and expert analysis).

Despite the development of quite a large number of multifidelity methods, their capabilities are still under discussion and their potential is still under-explored \cite{peherstorfer2018survey,fernandez2019-AIAA}. This motivates the interest in benchmark problems that could support the comparative and rigorous, reproducible and comparative assessment of these methods. Beran et al. \cite{beran2020comparison} propose to classify use cases and test problems into three classes: L1 problems, computationally cheap analytical functions with exact solutions; L2 problems, simplified engineering applications problems that can be executed with a reduced computational expense; and L3 problems, more complex engineering use cases, usually including multiphysics couplings. 

Following this classification, this paper captures the growing demand from the community and provides an overview of L1 benchmarks conceived as analytical problems with no explicit resemblance to actual engineering problems but that present mathematical features of relevance to cross-domain investigations. A large number of L1 benchmark problems have been proposed in the literature, mostly in conjunction with the presentation of a novel multifidelity method \cite{2018-IEEE-Wang_etal,2017-SMO-Durantin_etal,2017-AIAA-Cai_etal,rumpfkeil2020-AIAA,grassi2021-AIAA,Bryson2018a,bryson2016-AIAA}. However, a comprehensive framework of computationally efficient benchmarks is not yet available. 

This paper addresses the current lack of a cohesive and well-documented suite of L1 benchmarks proposing and discussing a suite of analytical benchmark problems and assessment metrics specifically formulated and selected to stress-test and assess the capabilities of a broad spectrum of multifidelity methods. The framework is intended to provide a set of standard problems, recommended experimental setups and performance assessment metrics to support the rigorous test and comparison of different computational methods. Specifically, the benchmark suite is composed of the Forrester function (continuous and discontinuous), the Rosenbrock function, the Rastrigin function (shifted and rotated), the Heterogeneous function, a coupled spring-mass system, and the Pacioreck function affected by noise; the set has been implemented and validated for Matlab, Fortran and Python environment with code made available to the public [https://gitlab.com/qudo046/avt-331-l1-benchmarks].

The L1 benchmark problems proposed in this work are carefully selected analytical functions that, while not derived from physical systems, reflect key mathematical characteristics and behaviors commonly encountered in simulation-based real-world optimization. These include challenges such as the curse of dimensionality and the resulting scalability issues for multifidelity methods \cite{bellman1957dynamic, koziel2012scaling}, as well as the need to navigate objective functions that may exhibit localized optima, multimodal landscapes, discontinuities, or irregular gradients \cite{lee2024comprehensive}. Some benchmarks also account for the presence of noise in the objective function evaluations, mimicking uncertainties that frequently arise in practical settings \cite{pellegrini2022multi}. Importantly, these problems are computationally inexpensive, relying on closed-form expressions that require no numerical solvers, which makes them ideal for rapid prototyping, extensive parametric studies, and integration into diverse multi-fidelity optimization frameworks. Their analytical nature also ensures a high degree of reproducibility across implementations and experimental setups. Moreover, the absence of simulation-related issues, such as solver instabilities, convergence failures, or discretization errors, allows researchers to evaluate the performance of multifidelity methods in a controlled environment where algorithmic behavior can be studied in isolation. Finally, because the global optima of these benchmarks are known by construction, the performance of optimization methods can be assessed with precise, objective metrics that quantify both global approximation fidelity and the accuracy and efficiency in locating the true optimal solution, capabilities that are otherwise difficult to evaluate in real-world applications where the ground truth is unknown or inaccessible.

The remainder of the paper is organized as follows. Section~\ref{sec:problems} illustrates the individual benchmark problems including their formulations and their distinguishing mathematical features. Section~\ref{sec:setup} presents the recommendations about the setup of the benchmark experiments for a fair and meaningful comparison of the methods. Section~\ref{sec:methods} summarizes a brief introduction to the assessed multifidelity approximations along with the utilized acquisition functions. Section~\ref{sec:metrics} discusses the different metrics and criteria to assess and compare the performance of multifidelity modelling and optimization strategies. Finally, Section~\ref{sec:results} presents the numerical results and concluding remarks are discussed in Section~\ref{sec:conclusion}.

\section{Analytical Benchmarks for Multifidelity Optimization}\label{sec:problems}
The analytical benchmarks proposed would exemplify potential objective functions to be optimized, somehow related to system-level design for complex industrial/military applications, thus solvable by using goal-driven multifidelity methods. Specifically, we will consider a box-constrained optimization problem in the form:
\beq \label{e:OptimizationProblem}
 \min_{{\bf x} \in \mathcal{A}} f({\bf x}), \qquad \mathrm{with} \qquad \bf{l}\leq\mathbf{x}\leq\bf{u},
\eeq
where $\bfx \in {\mathbb R}^D$ is a design point in the feasible domain $\mathcal{A}$ bounded by $\bf l$ (lower bound) and $\bf u$ (upper bound), $\{x_k\}$ are the elements of $\bfx$ ($k$ is an integer satisfying $1\leq k\leq D$ or $k\in[1,D]$), $D$ is the dimensionality of the parameter space, $f(\bfx) \in {\mathbb R}$ is the objective function, and $\bfx^\star$ is the optimum design point satisfying:
\beq \label{e:DesignPoint}
 \bfx^\star={\underset{\bfx \in \mathcal{A}}{\rm argmin}}f(\bfx),
\eeq
where $f^{\star}\equiv f(\bfx^\star)$.

Within the multifidelity settings of the benchmarks, the $f({\bf x})$ to minimize would be the highest-fidelity function $f_1({\bf x})$, while all the other possible $L$ representations would be considered as cheaper-to-evaluate approximations of the objective function, thus providing a fidelity spectrum from $f_1(\bfx)$ up to $f_L(\bfx)$, where the latter is the lowest-fidelity level available. 

\begin{table}[!t]
\small
\caption{Analytical benchmarks main features}
\centering
\begin{tabular}{cp{2.5cm}p{2.2cm}p{2.cm}p{2.1cm}p{1.cm}}
\toprule
\textbf{ID} & {\bf Name} & {\bf Behaviors} & {\bf Scalability} &  {\bf Discrepancy} & {\bf Noise}\\  
\midrule 
MF1 & Forrester & {Local / \newline (Dis)continous} & -          & (non)linear & no\\
MF2 & Rosenbrock & Local & Parametric & nonlinear   & no\\
\\
MF3 & {Shifted-rotated \newline Rastrigin} & Multi-modal & Parametric / \newline Fidelity  & nonlinear  & no  \\
MF4 & Heterogeneous & Local / \newline Multi-modal & Parametric  & nonlinear  & no    \\
MF5 & Spring-Mass  \newline system & Multi-modal  & Parametric / \newline Fidelity & nonlinear  & no  \\ 
MF6 & Paciorek & Multi-modal & Fidelity & nonlinear & yes\\
\bottomrule
\end{tabular}
\label{t:ExKind}
\end{table}

The multifidelity benchmark problems are selected to capture fundamental mathematical characteristics and properties that also mimic real-world engineering problems. The distinguishing mathematical features of the selected benchmark problems can be listed as behaviors, scalability, discrepancy type, and noise. Function behaviors include multi-modality, discontinuities, and atypical local behaviors, that cannot be neglected a priori, especially for real-world problems where it could be important to represent the whole variable domain. In particular, multi-modality and discontinuities are challenging from both a modeling and optimization viewpoint. Scalability takes into consideration both the function parameterization and the fidelity spectrum. The former is one of the important criteria for the comprehensive assessment of the performance of the multifidelity methods which enables to represent the same parametric function with different input dimensions, whereas the latter is useful to demonstrate how the modeling process can be improved depending on the fidelity level available. This last point is particularly relevant, because of its relation with the discrepancy type, which describes the relation among fidelities. In general, a linear discrepancy is simpler to model than a nonlinear one. For this reason, the discrepancy type allows for a deeper assessment of the multifidelity methods, since the number of fidelities available has to be correlated to the associated discrepancy type. Finally, in real-world engineering problems, noise may exist, it is undesired but inescapable in the overall response of a system which may exhibit abrupt changes within the solution domain. It is difficult or even impossible to distinguish the individual impact of the noise in the overall response and eliminate it. Thus, the user has to deal with a function with some embedded noise, and is important to assess the ability of multifidelity methods to model a noisy function.  
Considering all these mathematical characteristics over the benchmark problems will allow to evaluate the performance metrics to assess the strengths and weaknesses of the multifidelity methods employed. The mathematical characteristics that are considered in this work are briefly summarized in Tab. \ref{t:ExKind}.

The following subsections present each benchmark formulation in turn. 
\begin{figure}[!b]
\centering
\includegraphics[width=0.45\textwidth]{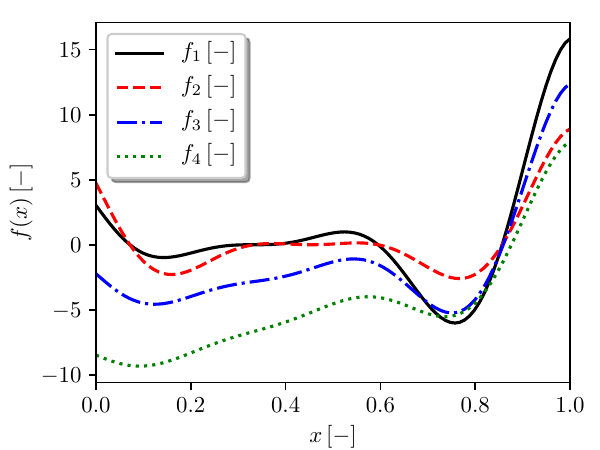} 
\includegraphics[width=0.45\textwidth]{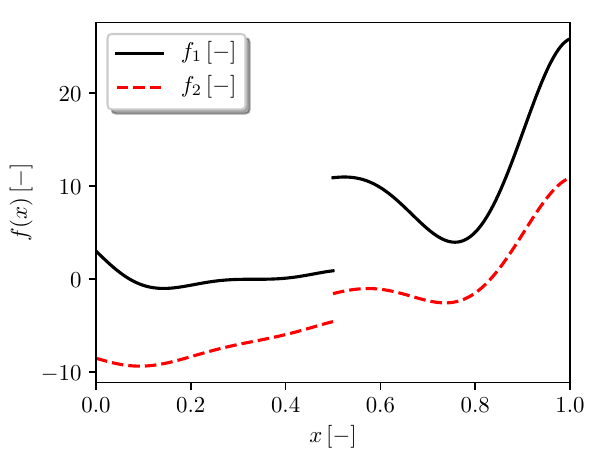} 
\caption{Forrester function (left) and discontinuous variant (right)}
\label{fig:forrester}
\end{figure}
\subsection{Forrester function}
The proposed Forrester \cite{forrester2007-PRSA} multifidelity benchmark (MF1.1) is a well-known one-dimensional benchmark for multifidelity methods, described by the following equations (from 1 highest fidelity to 4 lowest fidelity level) and shown in Fig.~\ref{fig:forrester}.

\begin{align}
f_1(x) &= \left(6x-2\right)^2\sin(12x-4)    \\
f_2(x) &= \left(5.5x-2.5\right)^2\sin(12x-4)   \\
f_3(x) &= 0.75f_1(x)+5(x-0.5)-2  \\
f_4(x) &= 0.5f_1(x)+10(x-0.5)-5   
\end{align}

The function is defined in the domain $0\leq x\leq 1$ and the minimum is located at $x^\star=0.75724876$ and given by $f(x^\star)= -6.020740$.  

To observe the performance of the multifidelity methods in problems with discontinuous behavior, the discontinuous Forrester function \cite{DiscForr2018} is also selected as one of the benchmarks. The benchmark function (MF1.2) is derived from the revision of Forrester's function and is also called the Forrester function with jump \cite{JumpForr}. The discontinuous Forrester function is described by the following equation 
\begin{equation}
f_{1}(x)=\left\{\begin{array}{lc}
(6x-2)^2 \sin(12x-4), & \ 0 \leq x \leq 0.5 \\
(6x-2)^2 \sin(12x-4)+10, & \ 0.5<x \leq 1
\end{array}\right.
\end{equation}
\begin{equation}
f_{2}(x)=\left\{\begin{array}{lc}
0.5 f_1(x) + 10 (x-0.5) -5, & \ 0 \leq x \leq 0.5 \\
0.5 f_1(x) + 10 (x-0.5) -2, & \ 0.5<x \leq 1
\end{array}\right.
\end{equation}

The function is defined in the domain $0 \leq x \leq 1$ and the minimum is located $x^* = 0.1426$ and given by $f(x^*)=-0.9863$.

\subsection{Rosenbrock function}
The Rosenbrock function, also referred to as the Valley or Banana function, is a well-known $D$-dimensional optimization benchmark problem described by the following equation
\begin{equation}
    f_1(\mathbf{x})=\sum_{i=1}^{D-1} 100\left(x_{i+1}-x_i^2\right)^2 + \left(1-x_i\right)^2
\end{equation}
The global minimum is inside a long, narrow, parabolic-shaped flat valley. 

The function is unimodal, and the global minimum is located at $\mathbf{x}^\star=\{1,\dots,1\}^{\sf T}$ and equal to $f(\mathbf{x}^\star)=0$. However, even though this valley is easy to find, convergence to the minimum is difficult. Note that there is also a local minimum at $\{-1,1,\dots,1\}^{\sf T}$ for $4 \le D \le 7$. For the present problem, the variable domain is defined as $-2\leq x_i \leq 2$ for $i=1,\dots D$.

The extension to multifidelity purposes is described by the following equations, where $f_2$ can be considered as a medium-fidelity level \cite{ficini2021-AIAA} and $f_3$ as the lowest fidelity \cite{bryson2016-AIAA}. 

\begin{equation}
    f_2(\mathbf{x})=\sum_{i=1}^{D-1}  50\left(x_{i+1}-x_i^2\right)^2 + \left(-2-x_i\right)^2 - \sum_{i=1}^D 0.5x_i
\end{equation}

\begin{equation}
    f_3(\mathbf{x})= \dfrac{f_1(\mathbf{x})-4-\sum_{i=1}^D 0.5x_i}{10+\sum_{i=1}^D 0.25x_1} 
\end{equation}

\begin{figure}[!t]
\centering
\includegraphics[width=0.325\textwidth]{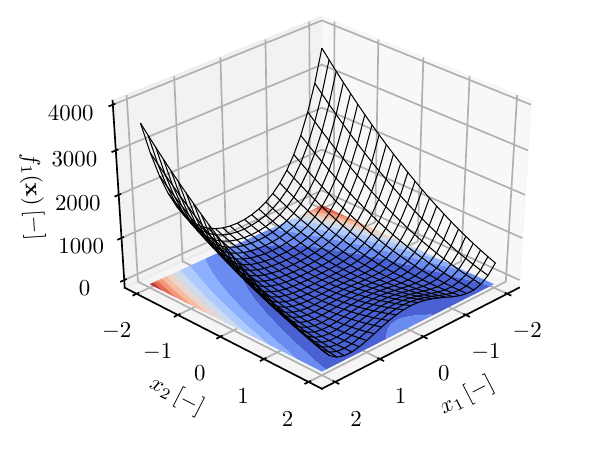} 
\includegraphics[width=0.325\textwidth]{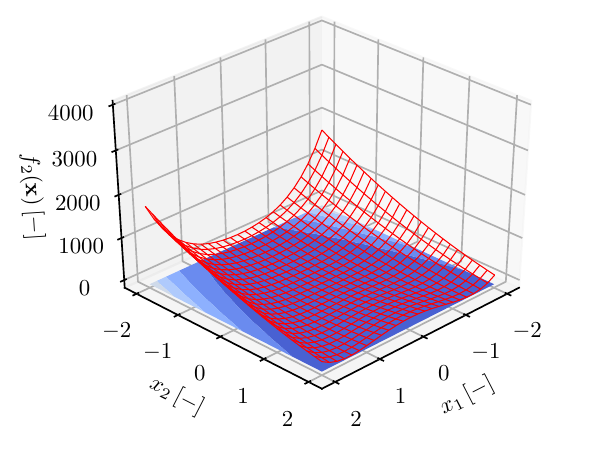} 
\includegraphics[width=0.325\textwidth]{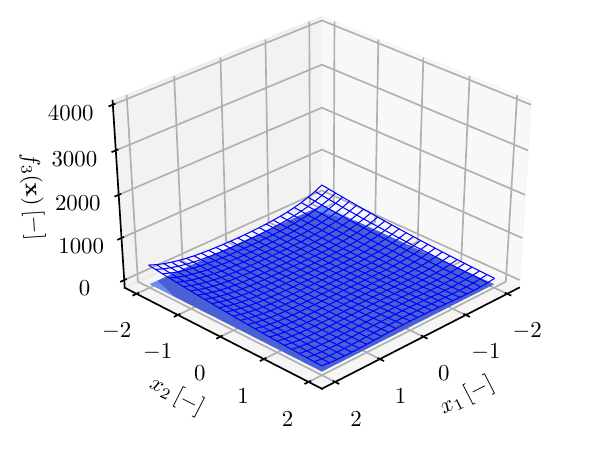} 
\caption{Rosenbrock Function: from left to right, $f_1$ (highest-fidelity), $f_2$ , and $f_3$ (lowest-fidelity)}
\label{fig:rosenbrock}
\end{figure}

The three-fidelity levels are shown for two dimensions in Figure \ref{fig:rosenbrock}.

\begin{figure}[!b]
\centering
\includegraphics[width=0.325\textwidth]{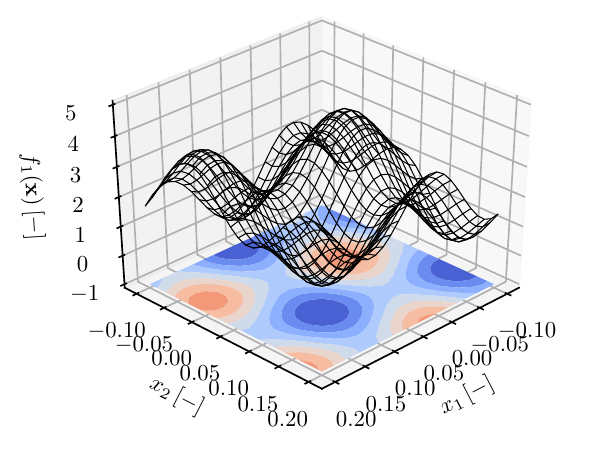} 
\includegraphics[width=0.325\textwidth]{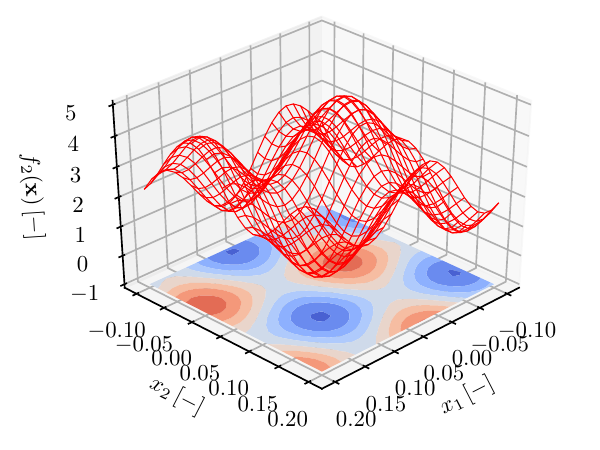} 
\includegraphics[width=0.325\textwidth]{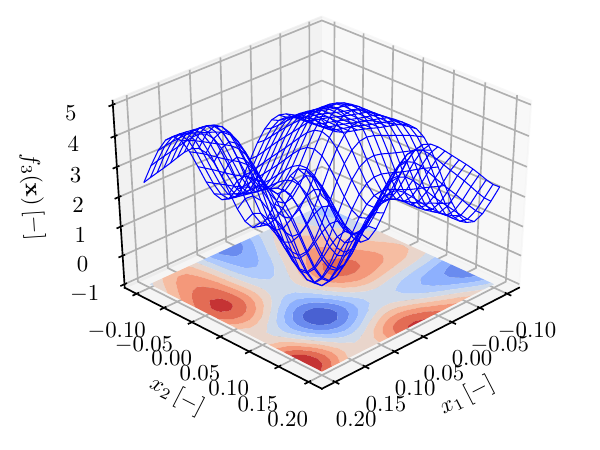} 
\caption{Shifted-rotated Rastrigin Function: from left to right, $f_1$ (highest-fidelity), $f_2$, and $f_3$ (lowest-fidelity)}
\label{fig:srrastrigin}
\end{figure}
\subsection{Shifted-rotated Rastrigin function}
To address real-world optimization problems, where the landscape of the objective function is usually multi-modal, the Rastrigin function is selected as a benchmark. The function is shifted to change the position of the minimum and rotated to change the properties of the function itself within the variable space. The equation of the shifted-rotated Rastrigin function reads as follows
\begin{equation}
    f_1(\mathbf{z})=\sum_{i=1}^D \left(z_i^2 +1 - \cos(10\pi z_i)\right)
\end{equation}
with
\begin{equation}
    \mathbf{z} = R(\theta)(\mathbf{x}-\mathbf{x}^\star) \,\,\,\,\, \mathrm{with} \,\,\,\,\, R(\theta)= \begin{bmatrix}
  \cos\theta & -\sin\theta \\
  \sin\theta &  \cos\theta \\
\end{bmatrix}
\end{equation}
where $R$ is the rotation matrix in 2$D$ and can be extended to arbitrary dimension by using the Aguilera-Perez algorithm \cite{aguilera2004general}.

The variable ranges is defined such as $-0.1 \leq x_i \leq 0.2$ for $i=1,\dots D$, rotation angle $\theta=0.2$, and optimum equal to $f(\mathbf{x}^\star)=0$ at $\mathbf{x}^\star=\{0.1,\dots,0.1\}^{\sf T}$.

The fidelity levels can be defined following the work of Wang et al. \cite{2018-IEEE-Wang_etal}, where a resolution error is defined as follows
\begin{equation}
    e_r(\mathbf{z},\phi)=\sum_{i=1}^D a(\phi)\cos^2(w(\phi)z_i+b(\phi)+\pi)
\end{equation}
with $a(\phi)=\Theta(\phi)$, $w(\phi)=10\pi\Theta(\phi)$, $b(\phi)=0.5\pi\Theta(\phi)$, and $\Theta(\phi)= 1-0.0001\phi$.
The fidelity levels are thus described as follows and depicted in Fig. \ref{fig:srrastrigin}.
\begin{equation}
    f_{i}(\mathbf{z},\phi_i)=f_1(\mathbf{z})+e_r(\mathbf{z},\phi_i) ~~~ \mathrm{for} ~~~ \, i=1,2,3
\end{equation}
with $\phi_1=10000$ (high-fidelity), $\phi_2=5000$ (medium-fidelity), and $\phi_3=2500$ (low-fidelity).

\subsection{Heterogeneous function}
This problem employs heterogeneous non-polynomial analytic functions defined on unit hypercubes ($0 \leq x_i \leq 1$ for $i=1,\dots D$) in one, two (taken from \cite{clark2016-AIAA}), and three dimensions, with low-fidelity functions obtained by using linear additive and multiplicative bridge functions. The benchmark reads for $D=1$ 
\beq\begin{cases}
f_1(x) & =\sin[30(x-0.9)^4] \cos[2(x-0.9)]+(x-0.9)/2\\
f_2(x) & =(f_1(x)-1.0+x)/(1.0+0.25x)
\end{cases}
\eeq
and for $D\geq2$

\beq\begin{cases}
f_1(\bfx)&=\sin[21(x_1-0.9)^4] \cos[2(x_1-0.9)]+(x_1-0.7)/2 + \\ 
& + \sum_{i=2}^D ix_i^i \sin \left(\prod_{j=1}^i x_j\right) \\
f_2(\bfx)&=(f_1(\bfx)-2.0+ \sum_{i=1}^D x_i)/(5.0 + \sum_{i=1}^2 0.25ix_i -\sum_{i=3 \atop D > 2}^D 0.25ix_i)
\end{cases}
\eeq

\begin{figure}[h]
\centering
\includegraphics[width=0.45\textwidth]{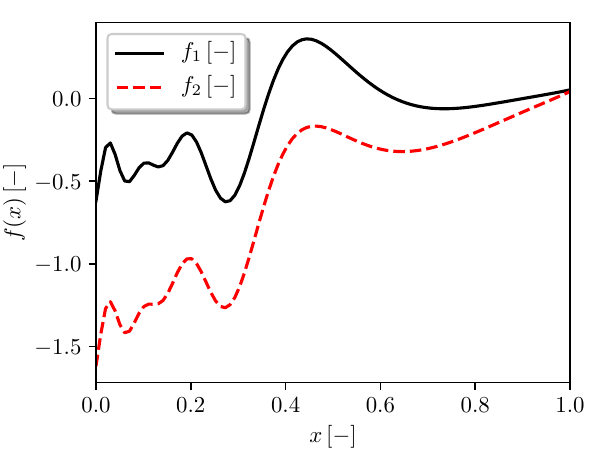} 
\caption{Heterogeneous function ($D=1$)}
\label{fig:clarkbae1d}
\end{figure}

\begin{figure}[h]
\centering
\includegraphics[width=0.325\textwidth]{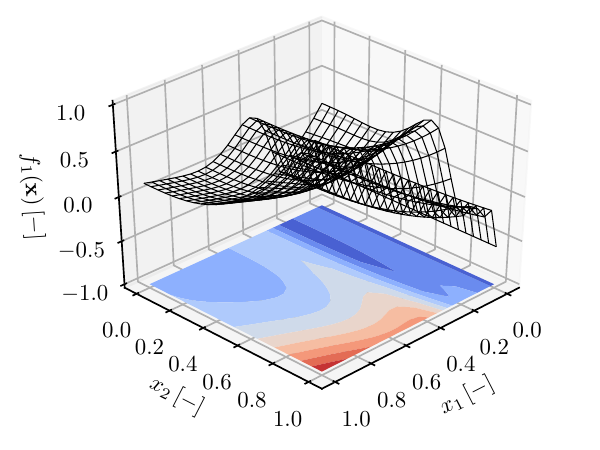} 
\includegraphics[width=0.325\textwidth]{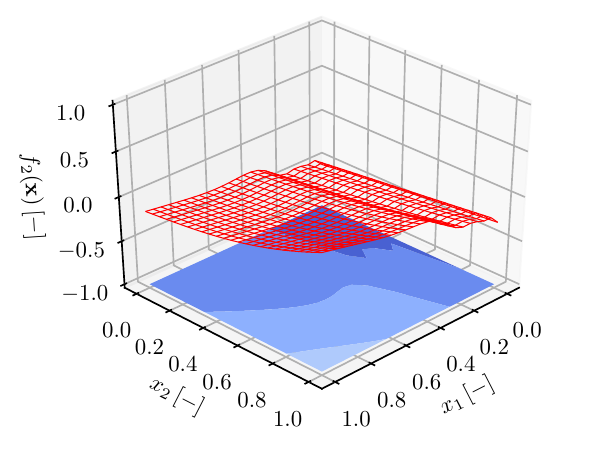} 
\caption{Heterogeneous function ($D=2$)}
\label{fig:clarkbae2d}
\end{figure}

The one- and two-dimensional functions are displayed together with their low-fidelity approximations in Figures \ref{fig:clarkbae1d} and \ref{fig:clarkbae2d}. The optimum for $D=1$ is equal to $f(\mathbf{x}^\star)=-0.625$ at $\mathbf{x}^\star=0.27550$, while for $D=2,3$ it is equal to $f(\mathbf{x}^\star)=-0.5627123$ at $\mathbf{x}^\star=\{0,\dots,0\}^{\sf T}$.


\subsection{Coupled Spring-Mass system}
\begin{figure}[h]
\centering
\includegraphics[width=0.5\textwidth]{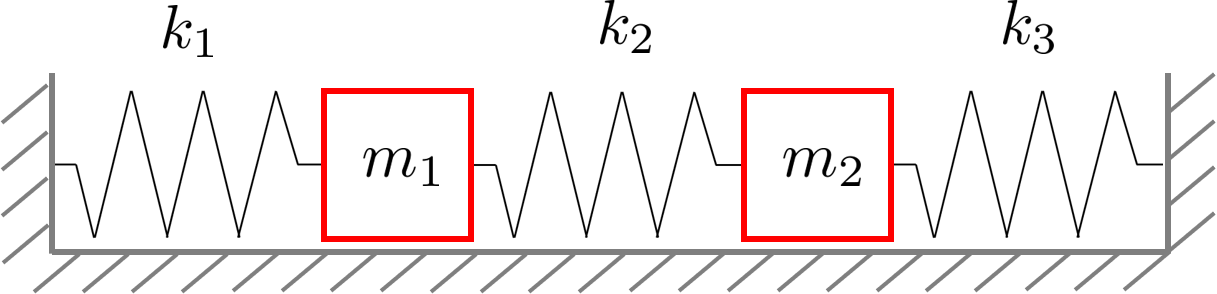} 
\caption{Coupled spring-masses system}
\label{3masses}
\end{figure}

The proposed problem represents a general coupled spring-mass system. Consider two point masses $m_1$ and $m_2$ concentrated at their center of gravity and attached by three springs, that operate according to Hooke’s law with constants $k_1$, $k_2$, and $k_3$. The mass of each spring is negligible and they are restored after compression and extension. The masses can slide along a frictionless horizontal surface, while the first and last springs ($k_1$ and $k_3$) are attached to fixed walls, as shown in Figure \ref{3masses}.
\begin{figure}[h]
\centering
\includegraphics[width=0.325\textwidth]{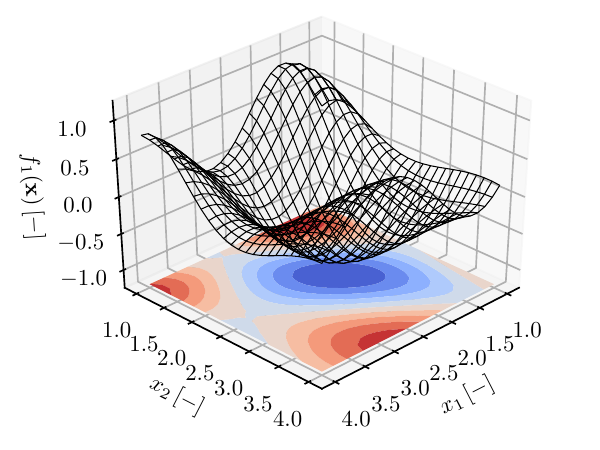} 
\includegraphics[width=0.325\textwidth]{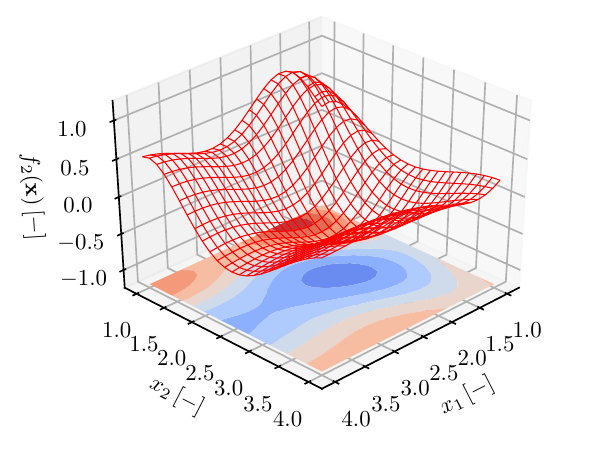} 
\caption{Spring-mass system ($D=2$, springs only)}
\label{fig:spring}
\end{figure}

Assume that $x_1(t)$ and $x_2(t)$ denote the mass positions along the horizontal surface, measured from their equilibrium positions, positive right and negative left. The equations of motion are given by the following system
\beq\begin{cases}\label{eqofmot}
m_1 \ddot{x}_1(t) &= -k_1 x_1 (t) + k_2 [x_2 (t) - x_1 (t)] \\
m_2 \ddot{x}_2(t) &= -k_2 [x_2 (t) - x_1 (t)] - k_3 x_2 (t)
\end{cases}
\eeq
Eq.~(\ref{eqofmot}) can be written as a second-order vector-matrix system
\beq
{\bf M} \ddot{{\bf x}} (t) = {\bf Kx}(t)
\label{secorderODE}
\eeq
where the displacement ${\bf x}$, mass matrix ${\bf M}$ and stiffness matrix $\bf K$ are defined by the following formulae
\bdm
{\bf x} = \left\{ \begin{array}{c} x_1 \\ x_2 \end{array} \right \} \quad {\bf M} = \left[ \begin{array}{cc} m_1 & 0 \\ 0 & m_2 \end{array} \right] \quad 
{\bf K} = \left[ \begin{array}{cc} -k_1-k_2 & k_2 \\ k_2 & -k_2-k_3 \end{array} \right]
\edm

This is a constant-coefficient homogeneous system of second-order ODEs the solution of which is given by 
\beq
\label{solutionODE}
{\bf x} (t)= \sum_{i=1}^2 [a_i \cos(\omega_i t) + b_i \sin(\omega_i t)] {\bf z}_i
\eeq
where $\omega_i = \sqrt{-\lambda_i}$ and $\lambda_i$ are the eigenvalues of the matrix ${\bf M}^{-1} {\bf K}$ and ${\bf z}_i$ are the corresponding eigenvectors. The constants $a_i$ and $b_i$ are determined by the initial conditions ${\bf x} (t=0) = {\bf x}_0$ and $\dot{\bf x} (t=0) = \dot{\bf x}_0$.

Converting Eq. \ref{secorderODE} into a system of first-order ODEs and using the fourth-order accurate Runge-Kutta time-marching method yields a multifidelity analysis problem by varying the time-step size $\Delta t$. The proposed benchmark uses the initial conditions ${\bf x}_0 = \{1 \; 0\}^{\sf T}$ and $\dot{\bf x}_0 = \{0 \; 0\}^{\sf T}$ with two fidelity levels, defined by the time-step size, specifically equal to $\Delta t = 0.01$ and 0.6. Two tests are proposed, considering the position of the first mass $x_1$ at time $t=6$ as the objective function: (MF5.1) springs $k_1$ and $k_2$ are the independent input variables with $1 \le (k_1,k_2) \le 4$ and $k_3=k_1$, while the masses are constant ($m_1=m_2=1$), the optimum is equal to $f(\mathbf{x}^\star)=-1$ at $\mathbf{x}^\star=\{2.467401, 2.193245\}^{\sf T}$; (MF5.2) springs $k_1$ and $k_2$ and masses $m_1$ and $m_2$ are the independent input variables with $1 \le (k_1,k_2,m_1,m_2) \le 4$ and $k_3=k_1$, the optimum is equal to $f(\mathbf{x}^\star)=-1$ at $\mathbf{x}^\star=\{1.000000, 3.946018, 4.000000, 3.286277\}^{\sf T}$. The two-dimensional (springs only) design space is shown in Fig. \ref{fig:spring}. 

\subsection{Paciorek Function}
The Paciorek equation \cite{Toal2014}, which has localized and multi-modality properties, is considered and a normally distributed random noise parameter is added to the high- and low-fidelity equation to model the noise. The Paciorek function with noise term is defined as following

\beq
f_1(\bfx)=\sin\left(\prod_{i=1}^D x_i\right)^{-1} +  rand.normal(0,\alpha_1)
\eeq

\beq
f_2(\bfx)=f_1(\bfx)-9A^2\cos\left(\prod_{i=1}^D x_i\right)^{-1} + rand.normal(0, \alpha_2)
\eeq

\begin{figure}[h]
\centering
\includegraphics[width=0.325\textwidth]{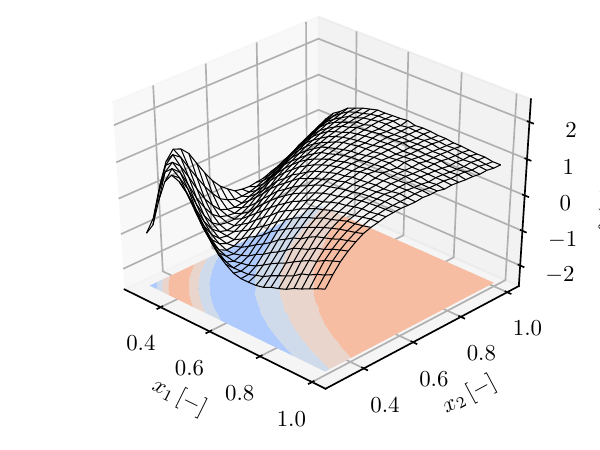}
\includegraphics[width=0.325\textwidth]{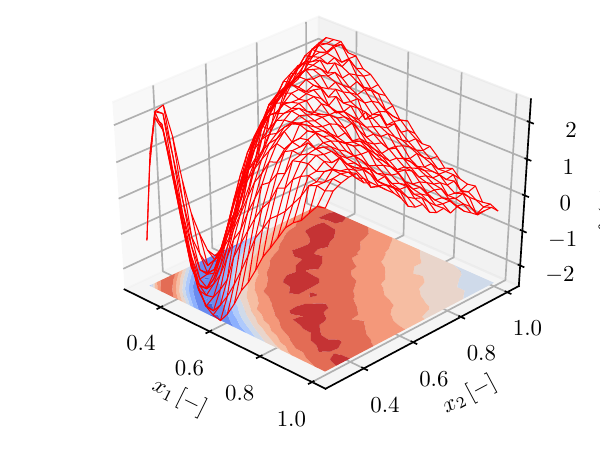}
\caption{Paciorek function with noise: $f_1$ (left) and $f_2$ (right)}
\label{fig:PaciorekPlot}
\end{figure}

$A$ is a parameter that models the error among the fidelities and can vary between 0 and 1, when $A=0$ the low- and high-fidelity equations are the same, while the change between the low-fidelity and high-fidelity models increases as the value of $A$ increases. In this study, it is set as $A=0.5$, while $\alpha_i$ is the coefficient that defines the noise level. The benchmark is defined considering the variable ranges as $0.3 \leq x_i \leq 1.0$ for $\forall i$. For the high-fidelity case, a noise level is added to the equation which corresponds to approximately 5\% of the response interval, as $\alpha_1=0.0125$. Accordingly, for the low-fidelity case much higher level of noise is added to the equation corresponding to 10\% of the response interval, $\alpha_2=0.075$. The resulting response surfaces are shown in Fig.~\ref{fig:PaciorekPlot}.

\section{Setup of the Numerical Experiments}\label{sec:setup}
Table \ref{t:ExSetup} summarizes the recommended setup for the assessment and comparison of the performance of multifidelity methods. The following subsections provide an overview of the rationale and the criteria motivating the recommendations along with additional complementary suggestions. The overview encompasses criteria for the assignment of the evaluation costs to the different fidelity levels as a fraction of the computational expense of the highest fidelity function, criteria for the initialization of the search, and criteria to terminate the search. 

\begin{table}[h]
\caption{Experiments setup summary}
\centering
\begin{tabular}{lccccccc}
\toprule
\multirow{3}{*}{\bf Function} & \multirow{3}{*}{\bf MF ID} & \multirow{3}{*}{$D$}  & \multirow{3}{*}{\bf Budget } & \multicolumn{4}{c}{\bf Fidelity cost} \\ 
\cmidrule{5-8}
&      &    &     & $f_1$ & $f_2$  & $f_3$      & $f_4$ \\ 
\midrule 
Forrester      & MF1.1 & 1  & 100 & 1.00E-0 & 5.00E-1 & 1.00E-1 & 5.00E-2  \\
Jump Forrester & MF1.2 & 1  & 100 & 1.00E-0 & 2.00E-1 & -          & -     \\
\midrule
               & MF2.1 & 2  & 200 & 1.00E-0 & 5.00E-1 & 1.00E-1 & - \\
Rosenbrock     & MF2.2 & 5  & 500 & 1.00E-0 & 5.00E-1 & 1.00E-1 & - \\
               & MF2.3 & 10 & 1000& 1.00E-0 & 5.00E-1 & 1.00E-1 & - \\
\midrule
\multirow{3}{*}{Shifted-rotated Rastrigin} & MF3.1 & 2  & 200 & 1.00E-0 & 6.25E-2 & 3.91E-3 & - \\
               & MF3.2 & 5  & 500 & 1.00E-0 & 6.25E-2 & 3.91E-3 & - \\
               & MF3.3 & 10 & 1000& 1.00E-0 & 6.25E-2 & 3.91E-3 & - \\ 
\midrule
               & MF4.1 & 1  & 100 & 1.00E-0 & 2.00E-1 &  -         & - \\
Heterogeneous  & MF4.2 & 2  & 200 & 1.00E-0 & 2.00E-1 &  -         & - \\
               & MF4.3 & 3  & 300 & 1.00E-0 & 2.00E-1 &  -         & -\\
\midrule
Springs        & MF5.1 & 2  & 200 & 1.00E-0 & 1.67E-2 &  -         & -   \\
Springs-masses & MF5.2 & 4  & 400 & 1.00E-0 & 1.67E-2 &  -         & -    \\
\midrule
Pacioreck      & MF6   & 2  & 200 & 1.00E-0 & 2.00E-1 &  -         & - \\ 
\bottomrule
\end{tabular}
\label{t:ExSetup}
\end{table}

\subsection{Fidelity cost assignment criteria}\label{s:FidelityCost}
The $f_1({\bf x})$ functions are given the unitary cost $\lambda_1=1$, while their lower fidelity representations are assigned cost values as fractions of the $f_1({\bf x})$ cost. Table \ref{t:ExSetup} proposes a set of cost assignments for the different fidelity representations of each of the benchmark problems. The values indicated for the shifted-rotated Rastrigin are determined according to the non-linear function proposed by Wang et al. \cite{2018-IEEE-Wang_etal} for the allocation of cost values to an arbitrary number of fidelities:
\beq \label{e:CostF1}
\lambda_l = \left(1/2^{l-1}\right)^4 \qquad \forall \ l\geq 1
\eeq
where higher values of the integer $l$ indicate representations of the objective functions at progressively lower levels or fidelity. Differently, the values indicated for the Forrester, the Rosenbrock, the Heterogeneous, and the Spring-Mass problems are determined on the experiences shared by the international community of methods developers and practitioners.

\subsection{Initialization criteria} \label{s:Init}
To assure a fair comparison across different families of multifidelity methods, statistics over a set of different starting samples/points is highly recommended. 
The cardinality and composition of the initial set are not constrained but determined for or by the specific multifidelity method. This approach is preferred over the definition of specific initial samples to be used. The performance assessment would then consider statistics over the set of all the convergence histories. 

\subsection{Termination criteria} \label{s:Term}
Real-world design problems are constrained by the limited time and computing resources available to conduct analysis, search, and optimization of the alternative candidate solutions. 
This motivates the choice to recommend conducting the experiments at given computational budgets assigned to the overall modeling and optimization task, rather than prescribing the maximum number of allowed iterations. The termination condition is reached when no computational budget is left. Table \ref{t:ExSetup} indicates the computational budget assigned to each experiment in terms of an equivalent number of evaluations of the high-fidelity representation.

\section{Multifidelity Approximations}\label{sec:methods}
From a wider perspective, multifidelity strategies fold into two categories, as surrogate lean and surrogate rich methods. In the context of the study, various different multifidelity approaches and original strategies have been developed, which may or may not rely on the use of a surrogate model combining information from different sources. This section provides a comprehensive overview of the considered surrogate lean and rich methods compared in this study followed by the acquisition functions that govern sample selection.

\subsection{Surrogate Lean}\label{sec:sLean}
This group of multifidelity methods is designed to minimize reliance on surrogate models while emphasizing direct utilization of information from models with varying fidelity levels. These methods adopt a combination of low-fidelity and high-fidelity model evaluations, often informed by adaptive sampling techniques. In this context, it is particularly critical to minimize the use of the costly high-fidelity models, which low-fidelity models guide the search process to keep the cost of optimization feasible. In this category of methods, the Broyden–Fletcher–Goldfarb–Shanno (BFGS) and Asymmetrical Multi-Objective Genetic Algorithm (AMOGA) are implemented and assessed over benchmark suit.

\paragraph{Broyden–Fletcher–Goldfarb–Shanno (BFGS) algorithm}
The multifidelity BFGS method extends the single-fidelity quasi-Newton method to incorporate various fidelity levels via the unconstrained optimization approach by using a penalty function. This method levereges the the concept of an expected point to construct a localised, two-point correction along the search direction resulting from bound-constrained quadratic sub-problem. 

High-fidelity gradient information is used to construct Hessian to guide the search towards objective-reducing direction. To reduce the cost of optimization line search is conducted based on high-fidelity corrected low-fidelity data. Once the optimization fails the process continues with the high-fidelity data without any loss of data. The BFGS algorithm minimizes a quadratic approximate at the $k$-th iterate, $m^{(k)}(\mathbf{x})$,
\begin{equation}
    m^{(k)}(\mathbf{x})=f_h(\mathbf{x}^{(k)})+\nabla f_h^{(k)T}(\mathbf{x}-\mathbf{x}^{(k)})+\frac{1}{2}(\mathbf{x}-\mathbf{x}^{(k)})^T
\mathbf{H}^{(k)}(\mathbf{x}-\mathbf{x}^{(k)})\end{equation}
where $f_h$ and $\nabla f_h$ are high-fidelity objective and associated gradient, and $\mathbf{H}$ is the approximate Hessian. Multifidelity polynomial chaos is utilized to implement low-fidelity based on first-order polynomials to limit the computational cost \cite{Bryson2018a}.

\paragraph{Asymmetrical Multi-Objective Genetic Algorithm (AMOGA)} 
This method is a multifidelity optimization based on utilizing the trade-off between accurate but expensive high-fidelity $f$ with cheaper but less accurate low-fidelity evaluations $g$. The obective function space divided into zones based on guiding the search towards minimising $g$. By generating a small initial population objective is evaluated and crossover and mutation applied. 

Once the acquisition treshold is reached the population element is used as starting point $x_i^0$ for the single-objective evolutionary procedure with $g$ as the fitness function. The best element of the single-objective $x_i^f$ used in the external level and evaluate the population. In order to resolve failure in converging towards global extremum a chained evaluation is used which is the utilisation of the best element of single-objective run to re-evaluate f in the next generation once the better result achieved than the starting individual.

\subsection{Surrogate Rich}\label{sec:sRich}
Modern engineering optimization tasks often involve expensive, black-box objective functions—commonly defined through simulations—that limit the number of function evaluations feasible under computational constraints. In such settings, surrogate-based optimization has emerged as a powerful paradigm, wherein a computationally inexpensive model approximates the true objective function and guides the selection of future evaluations. 

In this group of methods dedicated to surrogate modeling techniques, Polynomial Chaos Expansions (PCE) and their sparse variants, Kriging and co-Kriging models, Gaussian Processes (GPs) and their multifidelity extensions using autoregressive schemes, and Radial Basis Functions (RBF), including the stochastic RBF formulation are assessed over the benchmark functions.

\paragraph{Polynomial Chaos Expansion and Sparse Polynomial Chaos Expansion}
Polynomial Chaos Expansion (PCE) is a spectral surrogate modeling method that represents a system response as a weighted sum of orthogonal polynomial basis functions. These basis functions are selected according to the probability distributions of the input random variables, making PCE particularly well-suited for uncertainty quantification tasks.

Let $\boldsymbol{\xi} = (\xi_1, \xi_2, \dots, \xi_d)$ denote a vector of $d$ input random variables. The model response $Y = f(\boldsymbol{\xi})$ is approximated by the truncated expansion:
\begin{equation}
    Y \approx \hat{Y} = \sum_{i=0}^{P} c_i \Psi_i(\boldsymbol{\xi}),
\end{equation}
where $\{\Psi_i(\boldsymbol{\xi})\}_{i=0}^P$ are multivariate orthogonal polynomials and $c_i$ are the deterministic coefficients to be computed. The number of terms $P$ grows combinatorially with $d$ and the polynomial degree $p$, as $P = \binom{d + p}{p}$ in the case of a full tensor product basis.

Sparse Polynomial Chaos Expansion (SPCE) mitigates this growth by identifying and retaining only a subset of significant terms in the expansion. This is typically achieved using sparsity-promoting algorithms such as Least Angle Regression (LARS), LASSO, or compressed sensing, which solve:
\begin{equation}
    \min_{\mathbf{c}} \| \mathbf{Y} - \mathbf{\Psi} \mathbf{c} \|_2^2 + \lambda \|\mathbf{c}\|_1,
\end{equation}
where $\mathbf{\Psi}$ is the design matrix and $\lambda$ is a regularization parameter. SPCE is particularly advantageous in high-dimensional settings where many polynomial terms are negligible.

\paragraph{Kriging and Co-Kriging}
Kriging is an interpolation-based surrogate modeling method that treats the model output as a realization of a Gaussian process (GP). It is especially effective for deterministic functions and provides both function estimates and uncertainty quantification.

Given a deterministic model $f(\mathbf{x})$, the Kriging model assumes:
\begin{equation}
    f(\mathbf{x}) = \mu(\mathbf{x}) + Z(\mathbf{x}),
\end{equation}
where $\mu(\mathbf{x})$ is a deterministic trend function (often assumed constant or linear) and $Z(\mathbf{x})$ is a zero-mean GP with covariance function $k(\mathbf{x}, \mathbf{x}')$. The GP prior implies that the predictions at new points are normally distributed with analytically computable mean and variance.

Given training data $\mathcal{D} = \{ \mathbf{X}, \mathbf{y} \}$, the Kriging predictor at a new point $\mathbf{x}^*$ is given by:
\begin{align}
    \hat{f}(\mathbf{x}^*) &= \mu + \mathbf{r}^\top \mathbf{R}^{-1}(\mathbf{y} - \mu \mathbf{1}), \\
    \sigma^2(\mathbf{x}^*) &= k(\mathbf{x}^*, \mathbf{x}^*) - \mathbf{r}^\top \mathbf{R}^{-1} \mathbf{r},
\end{align}
where $\mathbf{r}$ is the vector of covariances between $\mathbf{x}^*$ and training points, and $\mathbf{R}$ is the covariance matrix of the training data.

Co-Kriging generalizes Kriging to handle multiple fidelity levels. A common formulation for two fidelity levels is the autoregressive model:
\begin{equation}
    f_H(\mathbf{x}) = \rho f_L(\mathbf{x}) + \delta(\mathbf{x}),
\end{equation}
where $f_H$ and $f_L$ are the high- and low-fidelity functions respectively, $\rho$ is a scaling factor, and $\delta(\mathbf{x})$ is an independent GP modeling the discrepancy between fidelities.

This structure allows co-Kriging to combine inexpensive low-fidelity data with sparse high-fidelity observations to improve prediction accuracy. The co-Kriging model requires estimating multiple covariance structures and introduces additional hyperparameters, but it provides a powerful tool for multifidelity optimization and modeling under limited data budgets.

\paragraph{Gaussian Process Regression and Multifidelity Gaussian Process Models}
Gaussian Process Regression (GPR) is a flexible and powerful non-parametric method for surrogate modeling. A Gaussian Process defines a distribution over functions:
\begin{equation}
    f(\mathbf{x}) \sim \mathcal{GP}(m(\mathbf{x}), k(\mathbf{x}, \mathbf{x}')),
\end{equation}
where $m(\mathbf{x})$ is the mean function and $k(\mathbf{x}, \mathbf{x}')$ is the covariance (kernel) function. The choice of kernel controls properties such as smoothness and correlation length, with common choices including the squared exponential and Matérn kernels.

Given training data $\{ \mathbf{X}, \mathbf{y} \}$, the GP posterior at a new point $\mathbf{x}^*$ provides a predictive distribution:
\begin{align}
    \hat{f}(\mathbf{x}^*) &= \mathbf{k}_*^\top (\mathbf{K} + \sigma_n^2 \mathbf{I})^{-1} \mathbf{y}, \\
    \sigma^2(\mathbf{x}^*) &= k(\mathbf{x}^*, \mathbf{x}^*) - \mathbf{k}_*^\top (\mathbf{K} + \sigma_n^2 \mathbf{I})^{-1} \mathbf{k}_*,
\end{align}
where $\mathbf{k}_*$ is the covariance vector between $\mathbf{x}^*$ and training points, $\mathbf{K}$ is the covariance matrix of the training data, and $\sigma_n^2$ is the observation noise variance.

In multifidelity modeling, Gaussian processes can be extended using autoregressive schemes such as the Kennedy-O'Hagan (KOH) model:
\begin{equation}
    f_H(\mathbf{x}) = \rho f_L(\mathbf{x}) + \delta(\mathbf{x}),
\end{equation}
where $f_L$ is a low-fidelity GP model, $\rho$ is a correlation parameter, and $\delta(\mathbf{x}) \sim \mathcal{GP}(0, k_\delta(\mathbf{x}, \mathbf{x}'))$ is a GP modeling the discrepancy. This framework can be recursively extended for more than two fidelity levels.

Hyperparameters of the kernels and scaling coefficients are typically learned via maximum likelihood estimation. Multifidelity Gaussian Process models retain the Bayesian nature of GPR, offering uncertainty estimates and principled mechanisms for active learning or adaptive sampling.

\paragraph{Radial Basis Function Models and Stochastic Radial Basis Functions}
Radial Basis Function (RBF) surrogate models are interpolation techniques based on radial symmetry around sampled data points. These models are particularly advantageous in high-dimensional design spaces due to their simplicity, smoothness, and ability to approximate complex nonlinear behaviors.

Given a set of $N$ input-output pairs $\{ \mathbf{x}_i, y_i \}_{i=1}^N$, an RBF model approximates the target function as:
\begin{equation}
    \hat{f}(\mathbf{x}) = \sum_{i=1}^N \lambda_i \phi(\| \mathbf{x} - \mathbf{x}_i \|),
\end{equation}
where $\lambda_i$ are weights determined by solving an interpolation system, and $\phi(\cdot)$ is a radial basis function depending only on the Euclidean distance between input vectors. Common choices for $\phi(r)$ include:
\begin{itemize}
    \item Gaussian: $\phi(r) = \exp(-\epsilon^2 r^2)$,
    \item Multiquadric: $\phi(r) = \sqrt{r^2 + \epsilon^2}$,
    \item Thin Plate Spline: $\phi(r) = r^2 \log(r)$,
\end{itemize}
where $\epsilon$ is a shape parameter that controls the spread of the basis function.

To improve numerical stability and allow for global trends, a polynomial term $p(\mathbf{x})$ is often added:
\begin{equation}
    \hat{f}(\mathbf{x}) = \sum_{i=1}^N \lambda_i \phi(\| \mathbf{x} - \mathbf{x}_i \|) + p(\mathbf{x}),
\end{equation}
and the interpolation conditions are enforced using a system of equations that includes the RBF kernel matrix and a polynomial constraint matrix.

Stochastic Radial Basis Function (SRBF) methods build on the RBF foundation by incorporating randomness and uncertainty modeling, particularly in surrogate-assisted optimization frameworks. SRBFs are often used in global optimization strategies, where a new point is selected not only based on predicted performance but also on the potential to improve exploration.

A common approach integrates RBF surrogates with probabilistic sampling strategies such as:
\begin{equation}
    \mathbf{x}_{\text{new}} = \arg\min_{\mathbf{x} \in \mathcal{X}} \hat{f}(\mathbf{x}) + \beta \cdot \sigma(\mathbf{x}),
\end{equation}
where $\hat{f}(\mathbf{x})$ is the RBF prediction, $\sigma(\mathbf{x})$ is an uncertainty measure derived from the surrogate, and $\beta$ controls the trade-off between exploitation and exploration. Since classical RBFs do not provide a natural uncertainty estimate, $\sigma(\mathbf{x})$ is typically approximated using distance-based heuristics or ensemble models.

\subsubsection{Acquisition Functions}
Building upon these surrogate lean and rich models, the subsequent sections introduce a suite of acquisition functions designed to strategically balance exploration and exploitation. These include Expected Improvement (EI) and Probability of Improvement (PI), information-theoretic criteria such as Entropy Search (ES) and Max-Value Entropy Search (MES), Lower Confidence Bound (LCB) strategies and evolution control (EC) mechanisms, and Global-Local Regression Discrepancy (GLRD) along with their multifidelity variants. 

\paragraph{Expected Improvement and Multifidelity Extension}
Expected Improvement (EI) quantifies the expected amount by which a new evaluation at a point $\mathbf{x} \in \mathbb{R}^d$ is anticipated to improve upon the best observed objective value $f_{\min}$. Assuming that the objective function $f(\mathbf{x})$ is modeled using a Gaussian process (GP), with predictive mean $\mu(\mathbf{x})$ and standard deviation $\sigma(\mathbf{x})$, the EI can be computed in closed form as:

\begin{equation}
\text{EI}(\mathbf{x}) = \mathbb{E}\left[\max(f_{\min} - f(\mathbf{x}), 0)\right] = (f_{\min} - \mu(\mathbf{x})) \Phi(z) + \sigma(\mathbf{x}) \phi(z),
\end{equation}

\noindent where $z = \frac{f_{\min} - \mu(\mathbf{x})}{\sigma(\mathbf{x})}$, and $\Phi(\cdot)$ and $\phi(\cdot)$ denote the cumulative distribution function (CDF) and probability density function (PDF) of the standard normal distribution, respectively. The first term in the EI expression promotes exploitation by favoring areas with low predicted mean, while the second term encourages exploration by promoting regions of high uncertainty.

Multifidelity Expected Improvement (MFEI) generalizes EI to such settings by integrating information from multiple fidelity levels indexed by $l \in \{1, 2, \dots, L\}$, where $L$ denotes the highest-fidelity level. A general form of the MFEI criterion is given by:

\begin{equation} \label{eq:MFEI_general}
\text{MFEI}(\mathbf{x}, l) = \text{EI}(\mathbf{x}, L) \cdot \alpha_1(\mathbf{x}, l) \cdot \alpha_2(\mathbf{x}, l) \cdot \alpha_3(l),
\end{equation}

\noindent where $\text{EI}(\mathbf{x}, L)$ is the expected improvement computed using the high-fidelity GP surrogate, and the three multiplicative factors $\alpha_1$, $\alpha_2$, and $\alpha_3$ represent fidelity-specific utility corrections.

The utility term $\alpha_1(\mathbf{x}, l)$ encodes the statistical correlation between the objective function evaluated at fidelity level $l$ and at the highest fidelity $L$:

\begin{equation} \label{eq:MFEI_alpha1}
\alpha_1(\mathbf{x}, l) = \text{corr} \left[ f^{(l)}(\mathbf{x}), f^{(L)}(\mathbf{x}) \right].
\end{equation}

\noindent This term captures how informative an evaluation at fidelity $l$ is with respect to the true objective. High correlation implies that low-fidelity evaluations are predictive of high-fidelity behavior and therefore valuable.

The second term, $\alpha_2(\mathbf{x}, l)$, reflects the contribution of the new evaluation at fidelity $l$ to the reduction of the GP posterior uncertainty. It accounts for the expected variance reduction from incorporating a new point in the multifidelity Gaussian process model:

\begin{equation} \label{eq:MFEI_alpha2}
\alpha_2(\mathbf{x}, l) = 1 - \frac{\sigma_\epsilon}{\sqrt{\sigma^{2(l)}(\mathbf{x}) + \sigma^2_\epsilon}},
\end{equation}

\noindent where $\sigma^{2(l)}(\mathbf{x})$ is the predictive variance at fidelity $l$, and $\sigma^2_\epsilon$ denotes the model discrepancy or observational noise at that fidelity. This term tends toward one in regions where the variance is large compared to the noise, and it vanishes when the noise dominates the signal.

The third term, $\alpha_3(l)$, encodes the relative cost-effectiveness of fidelity level $l$ compared to the highest fidelity level $L$:

\begin{equation} \label{eq:MFEI_alpha3}
\alpha_3(l) = \frac{\lambda^{(L)}}{\lambda^{(l)}},
\end{equation}

\noindent where $\lambda^{(l)}$ is the computational cost of performing a function evaluation at fidelity $l$. This term favors the selection of lower-fidelity queries when their cost is substantially smaller and the other utility terms remain sufficiently high.

\paragraph{Probability of Improvement and Multifidelity Extension}
Probability of Improvement (PI) is designed to evaluate how likely a candidate point $\mathbf{x} \in \mathbb{R}^d$ is to improve upon a given target value, typically the current best observed function value $f_{\min}$. When the objective function is modeled by a Gaussian process (GP) surrogate, with predictive mean $\mu(\mathbf{x})$ and standard deviation $\sigma(\mathbf{x})$, PI is defined as:

\begin{equation}
\text{PI}(\mathbf{x}) = \mathbb{P}(f(\mathbf{x}) < f_{\min} - \xi) = \Phi\left( \frac{f_{\min} - \mu(\mathbf{x}) - \xi}{\sigma(\mathbf{x})} \right),
\end{equation}

\noindent where $\Phi(\cdot)$ denotes the cumulative distribution function of the standard normal distribution, and $\xi \geq 0$ is a hyperparameter used to control the trade-off between exploration and exploitation. Small values of $\xi$ result in more aggressive exploitation of the current minimum, while larger values encourage exploration in regions with greater predictive uncertainty.

The multifidelity variant of PI incorporates fidelity-specific information by modulating the base PI score computed at the highest fidelity level $L$ with fidelity-aware correction factors. Let $l \in \{1, \dots, L\}$ denote the index of a lower-fidelity model. The multifidelity PI acquisition function is defined as:

\begin{equation}
\text{MFPI}(\mathbf{x}, l) = \text{PI}(\mathbf{x}, L) \cdot \beta_1(\mathbf{x}, l) \cdot \beta_2(\mathbf{x}, l) \cdot \beta_3(l),
\end{equation}

\noindent where $\text{PI}(\mathbf{x}, L)$ represents the high-fidelity acquisition score, and the $\beta_i$ terms weight the contribution of fidelity $l$ to guide adaptive sampling. 

The first correction factor, $\beta_1(\mathbf{x}, l)$, captures the statistical correlation between the predictions of the $l$-th fidelity model and those of the highest fidelity:

\begin{equation}
\beta_1(\mathbf{x}, l) = \text{corr}\left[ f^{(l)}(\mathbf{x}), f^{(L)}(\mathbf{x}) \right],
\end{equation}

\noindent favoring fidelities that are more predictive of the high-fidelity response at the query location.

The second term, $\beta_2(\mathbf{x}, l)$, incorporates a penalization based on sampling density, designed to discourage redundant sampling in regions already well explored at fidelity $l$. This is modeled using a radial correlation function $R(\cdot, \cdot)$ over the $n_l$ samples available at level $l$:

\begin{equation} \label{e:MFPI3}
\beta_2(\mathbf{x}, l) = \prod_{i=1}^{n_l} \left[ 1 - R\left(\mathbf{x}, \mathbf{x}_i^{(l)} \right) \right].
\end{equation}

\noindent This term diminishes the acquisition score in high-density regions, effectively implementing a local repulsion mechanism that promotes broader exploration of the input domain.

The third term, $\beta_3(l)$, encodes the cost-efficiency of evaluating at fidelity $l$ relative to the highest fidelity:

\begin{equation}
\beta_3(l) = \frac{\lambda^{(L)}}{\lambda^{(l)}},
\end{equation}

\noindent where $\lambda^{(l)}$ and $\lambda^{(L)}$ are the computational costs of a single function evaluation at fidelity $l$ and $L$, respectively. This weighting encourages the use of lower-fidelity models when their cost is substantially lower and their informativeness (captured by $\beta_1$ and $\beta_2$) remains high.

\paragraph{Entropy-based Acquisition Functions: Entropy Search, Max-value Entropy Search, and Multifidelity Extensions}
Entropy-based acquisition functions aim to reduce uncertainty about the location of the global optimum $\mathbf{x}^*$ of an unknown objective function $f: \mathbb{R}^d \rightarrow \mathbb{R}$ modeled via a Gaussian process (GP). These methods view Bayesian optimization as a sequential information-gathering process and select new evaluation points to maximally reduce the posterior entropy of $\mathbf{x}^*$ or related quantities.

Entropy Search (ES) explicitly seeks to reduce the Shannon entropy $\mathcal{H}[p(\mathbf{x}^* \mid \mathcal{D})]$ of the probability distribution over the minimizer $\mathbf{x}^*$, given the current data $\mathcal{D}$. At each iteration, ES selects the point $\mathbf{x}$ whose evaluation is expected to result in the greatest reduction of this entropy:

\begin{equation}
\text{ES}(\mathbf{x}) = \mathbb{E}_{f(\mathbf{x})} \left[ \mathcal{H}[p(\mathbf{x}^* \mid \mathcal{D})] - \mathcal{H}[p(\mathbf{x}^* \mid \mathcal{D} \cup \{(\mathbf{x}, f(\mathbf{x}))\})] \right].
\end{equation}

Computing this acquisition function requires approximating the posterior distribution over $\mathbf{x}^*$ and performing nested expectations, making ES computationally intensive, especially in high dimensions. Nonetheless, it provides a principled framework for active learning of the optimizer.

Max-value Entropy Search (MES) offers a simplification of ES by focusing on the entropy of the minimum value $f^* = \min_{\mathbf{x} \in \mathcal{X}} f(\mathbf{x})$ rather than the minimizer's location. The MES acquisition function is defined as the mutual information between $f^*$ and the function value at a candidate location $\mathbf{x}$:

\begin{equation}
\text{MES}(\mathbf{x}) = \mathcal{H}[f^* \mid \mathcal{D}] - \mathbb{E}_{f(\mathbf{x})} \left[ \mathcal{H}[f^* \mid \mathcal{D} \cup \{(\mathbf{x}, f(\mathbf{x}))\}] \right].
\end{equation}

Under the GP prior, this mutual information can be approximated efficiently using samples from the distribution of $f^*$ and the GP predictive distribution. MES achieves a favorable balance between exploration and exploitation while offering improved scalability over ES.

In the presence of multiple fidelity models or heterogeneous information sources indexed by $l \in \{1, \dots, L\}$, entropy-based acquisition functions can be extended to select both a location $\mathbf{x}$ and a fidelity level $l$ that jointly maximize the expected information gain per unit cost.

The multifidelity variant of MES seeks to maximize the information gain on $f^*$ while accounting for the informativeness and cost of querying a lower-fidelity model. A general form of the multifidelity MES acquisition function is given by:

\begin{equation}
\text{MF-MES}(\mathbf{x}, l) = \frac{I(f^*; f^{(l)}(\mathbf{x}))}{\lambda^{(l)}},
\end{equation}

\noindent where $I(f^*; f^{(l)}(\mathbf{x}))$ denotes the mutual information between the global minimum $f^*$ and the function value at location $\mathbf{x}$ from fidelity level $l$, and $\lambda^{(l)}$ is the cost of evaluating fidelity $l$. This formulation ensures that the acquisition strategy remains cost-effective by favoring fidelity levels with higher information efficiency.

Approximating $I(f^*; f^{(l)}(\mathbf{x}))$ typically involves Monte Carlo estimation based on GP predictions and sampled realizations of $f^*$, possibly using the highest-fidelity surrogate $f^{(L)}(\mathbf{x})$ as reference. Extensions to this formulation can also include fidelity-specific correction factors to account for the statistical discrepancy between $f^{(l)}$ and $f^{(L)}$.

For multi-source settings, similar entropy-based acquisition criteria can be used by defining source-specific predictive models and computing their contributions to the reduction of uncertainty about either $f^*$ or $\mathbf{x}^*$. A generic multi-source entropy acquisition function can be written as:

\begin{equation}
\text{MS-ES}(\mathbf{x}, s) = \frac{I(T; f^{(s)}(\mathbf{x}))}{\lambda^{(s)}},
\end{equation}

\noindent where $s$ indexes the source, $T$ represents the target quantity (either $f^*$ or $\mathbf{x}^*$), and $\lambda^{(s)}$ is the corresponding cost.

\paragraph{Lower Confidence Bound (LCB) and Multifidelity Extension}
The Lower Confidence Bound (LCB) is a popular acquisition function that seeks to balance exploration and exploitation by selecting the next evaluation point with the lowest lower bound on the objective function's possible values. Given a Gaussian process (GP) model with predictive mean $\mu(\mathbf{x})$ and standard deviation $\sigma(\mathbf{x})$, the LCB acquisition function at a candidate point $\mathbf{x}$ is defined as:

\begin{equation}
\text{LCB}(\mathbf{x}) = \mu(\mathbf{x}) - \kappa \sigma(\mathbf{x}),
\end{equation}

where $\kappa > 0$ is a hyperparameter that controls the trade-off between exploration and exploitation. A small value of $\kappa$ emphasizes exploitation of the regions with the highest predicted mean, while larger values of $\kappa$ encourage exploration of regions with high uncertainty, where the model's predictions are more uncertain.

In the multifidelity optimization context, the LCB acquisition function can be extended to choose both the location $\mathbf{x}$ and the fidelity level $l$ that minimizes the lower confidence bound of the objective function at that location, taking into account both the predictive uncertainty from each fidelity model and the associated evaluation costs. The multifidelity LCB is defined as:

\begin{equation}
\text{MF-LCB}(\mathbf{x}, l) = \frac{\mu^{(l)}(\mathbf{x}) - \kappa \sigma^{(l)}(\mathbf{x})}{\lambda^{(l)}},
\end{equation}

where $\mu^{(l)}(\mathbf{x})$ and $\sigma^{(l)}(\mathbf{x})$ are the predictive mean and standard deviation of the objective function at $\mathbf{x}$ under the $l$-th fidelity model. $\kappa$ is the same hyperparameter used in the standard LCB, controlling the exploration-exploitation trade-off and $\lambda^{(l)}$ is the computational cost associated with evaluating the $l$-th fidelity model.

\paragraph{Evolution Control (EC)}
The Evolution Control acquisition function uses a strategy where the candidate point $\mathbf{x}$ for the next evaluation is selected to maximize an acquisition function that reflects both the uncertainty about the objective function and the expected gain from future evaluations. The general form of the EC acquisition function can be expressed as:

\begin{equation}
\text{EC}(\mathbf{x}) = \mu(\mathbf{x}) + \gamma \sigma(\mathbf{x}),
\end{equation}

where $\mu(\mathbf{x})$ is the predictive mean, $\sigma(\mathbf{x})$ is the predictive standard deviation, and $\gamma$ is a hyperparameter that controls the balance between exploration and exploitation. A higher value of $\gamma$ promotes exploration by favoring locations with higher uncertainty, while smaller values encourage exploitation of regions where the objective function is better understood. The evolution control strategy can be integrated into the optimization loop by adjusting the acquisition function iteratively based on the historical evaluation results. Specifically, the acquisition function at iteration $t$ is adapted to balance between exploration of new areas of the design space and exploitation of regions known to yield good results. This adjustment can be implemented as a dynamic weighting of the predictive mean and variance:

\begin{equation}
\text{EC}_t(\mathbf{x}) = \mu_t(\mathbf{x}) + \gamma_t \sigma_t(\mathbf{x}),
\end{equation}

where $\gamma_t$ is a time-dependent parameter that evolves with the optimization process, reflecting the current state of exploration versus exploitation. Early in the optimization, larger values of $\gamma_t$ are chosen to encourage exploration, and as the optimization progresses, $\gamma_t$ is gradually reduced to focus on exploitation of the most promising regions of the design space.

The multifidelity version of EC balances the contributions from each fidelity model based on their accuracy and computational cost:

\begin{equation}
\text{MF-EC}^{\text{cost}}(\mathbf{x}, l) = \frac{\mu^{(l)}(\mathbf{x}) + \gamma^{(l)} \sigma^{(l)}(\mathbf{x})}{\lambda^{(l)}},
\end{equation}

where $\mu^{(l)}(\mathbf{x})$ and $\sigma^{(l)}(\mathbf{x})$ represent the predictive mean and standard deviation under the $l$-th fidelity model, and $\lambda^{(l)}$ is the cost of evaluating the objective function at fidelity level $l$. The parameter $\gamma^{(l)}$ controls the balance between exploration and exploitation for the $l$-th fidelity level, and it may vary across different fidelity levels. A higher value of $\gamma^{(l)}$ promotes exploration within that fidelity, while a lower value encourages exploitation of regions with higher certainty about the objective function.

\paragraph{Global-Local Regression Discrepancy (GLRD)}
Global-Local Regression Discrepancy (GLRD) is an acquisition function designed to address the challenge of identifying regions in the design space where the surrogate model's predictive performance is poor. This is achieved by incorporating both the global and local information from the surrogate model’s predictions. The global regression term refers to the general behavior of the surrogate model across the entire design space, while the local regression term focuses on the deviations or discrepancies in small subregions of the design space. The GLRD acquisition function can be expressed as:

\begin{equation}
\text{GLRD}(\mathbf{x}) = \left| \mu_{\text{global}}(\mathbf{x}) - \mu_{\text{local}}(\mathbf{x}) \right| + \lambda \sigma(\mathbf{x}),
\end{equation}

where $\mu_{\text{global}}(\mathbf{x})$ represents the global predictive mean of the surrogate model, $\mu_{\text{local}}(\mathbf{x})$ is the local predictive mean (often obtained from a local regression model or local Kriging), and $\sigma(\mathbf{x})$ is the standard deviation. The term $\lambda$ is a hyperparameter that controls the weight between the local regression discrepancy and the uncertainty in the predictions, balancing exploration and exploitation. The function $\mu_{\text{global}}(\mathbf{x})$ is typically computed from the surrogate model (e.g., Gaussian Process or Kriging), which provides the overall trend of the objective function. In contrast, $\mu_{\text{local}}(\mathbf{x})$ is derived from a local regression approach, focusing on a neighborhood around $\mathbf{x}$ to capture the immediate variations in the objective function.

In the context of multifidelity optimization, the GLRD function can be extended to incorporate multiple fidelity models. This extension is valuable when different fidelities offer complementary information that can guide the optimization process more efficiently. The multifidelity version of GLRD combines information from different levels of fidelity to better capture both local and global discrepancies:

\begin{equation}
\text{MF-GLRD}(\mathbf{x}, l) = \left| \mu_{\text{global}}^{(l)}(\mathbf{x}) - \mu_{\text{local}}^{(l)}(\mathbf{x}) \right| + \lambda^{(l)} \sigma^{(l)}(\mathbf{x}),
\end{equation}

where $\mu_{\text{global}}^{(l)}(\mathbf{x})$ is the global predictive mean at the $l$-th fidelity level, $\mu_{\text{local}}^{(l)}(\mathbf{x})$ is the local predictive mean, and $\sigma^{(l)}(\mathbf{x})$ is the standard deviation at fidelity level $l$. The term $\lambda^{(l)}$ is a fidelity-specific parameter controlling the balance between local discrepancy and predictive uncertainty at each level of fidelity. 

The multifidelity approach can incorporate both low- and high-fidelity models, providing an efficient mechanism to guide the optimization towards regions of high discrepancy or uncertainty. The discrepancy between global and local predictions at each fidelity level highlights regions where the local approximation diverges from the global trend, suggesting that further exploration is needed in those regions. This is particularly beneficial when low-fidelity models are used to guide exploration and high-fidelity models are used to exploit the best solutions found so far.

Table~\ref{tab_ch06:l1_methods} summarizes the main characteristics of each method used. Specifically, two surrogate-lean versus seventeen surrogate-rich methods are compared. Among the latter group, the majority (14 out of 17) are GP-based, whereas 12 are online (including an adaptive sampling approach).

\begin{landscape}
\begin{table}[h]
\centering
\caption{Multifidelity Methods Summary}
\label{tab_ch06:l1_methods}
\tiny
\begin{tabular}{llp{1.7cm}lp{1.15cm}lp{1.7cm}p{1.65cm}lp{2.1cm}p{1cm}}
\toprule
\bf Classification& \bf Acronym & \bf Surrogate$^*$              & \bf Off/Online & \bf Adaptive sampling$^\dag$ & \bf DoE$^\ddag$ (D/S) & \bf Optimization algorithm$^\mathsection$ (D/S) & \bf Derivative$^{\ddag\ddag}$ informed/free&\bf Organization & \bf Notes & \bf Refs.\\
\midrule
Surrogate lean  & BFGS    & --                               & Online & --   & LHS (S)   & BFGS  (D)    & Informed & AFRL   & Highest- and lowest-fidelity only used & \cite{Bryson2018a}\\
                & AMOGA   & --                               & Online & --   & RBP (S)   & AMOGA (S)    & Free     & CIRA   & First (highest) and second-fidelity levels only used & \cite{quagliarella2004aerodynamic}\\
\midrule
\midrule
Surrogate rich  & GP      & Auto-regressive Gaussian process & Offline & --  & HSOS (S)    & MS\_SQP (D) & Informed        & ITU    & Highest- and lowest-fidelity only used & \cite{yildiz2022-ECCOMAS}\\
                & GP-PI   & Auto-regressive Gaussian process & Online  & PI  & LHS (S)    & GA+SQP (S)       & Informed    & PoliTo & Highest- and lowest-fidelity only used & \cite{ruan2020variable}\\
                & GP-EI   & Auto-regressive Gaussian process & Online  & EI  & LHS (S)    & GA+SQP (S)       & Informed        & PoliTo & Highest- and lowest-fidelity only used & \cite{huang2006sequential} \\
                & GP-EIM  & Auto-regressive Gaussian process & Online  & EIM & HSOS (S)    & MS\_SQP (D) & Informed        & ITU    & Highest- and lowest-fidelity only used & \cite{yildiz2022-ECCOMAS}\\
                & GP-MVES & Auto-regressive Gaussian process & Online  & MVES& LHS (S)    & GA+SQP (S)       & Free        & PoliTo & Highest- and lowest-fidelity only used & \cite{takeno2020multi}\\
                & GP-MSES & Auto-regressive Gaussian process & Online  & MSES& LHS (S)    & GA+SQP (S)       & Informed        & \it EMUKIT & Highest- and lowest-fidelity only used; tests by ITU, PoliTo & \cite{emukit2019}\\
                & GP-LCB  & Auto-regressive Gaussian process & Online  & LCB & CCF$^*$ (D)& DPSO+LS (D) & Free     & CNR-INM& Used with interpolation & \cite{ficini2021-AIAA}\\
           & NGP& Nonlinear auto-regressive Gaussian process & Offline & --  & HSOS (S)    & MS\_SQP (D) & Informed        & ITU    & Highest- and lowest-fidelity only used & \cite{perdikaris2017nonlinear}\\
                & KG      & Kriging                          & Offline & --  & HSS (D)    & GA (S)      & Free     & \it DAKOTA & Tests by CNR-INM & \cite{dakota2021}\\
                & KG-GLRD & Kriging                          & Online  & GLRD& FFD (D)    & --          & --       & AFRL   & Used for global accuracy only (NRMSE) & \cite{Rumpfkeil2016}\\ 
                & KG-EC   & Kriging                          & Online  & EC  & LHS (S)    & GA+SQP (S)  & Informed & UoS    & -- & \cite{Minisci:2013} \\
                & KG-AEI  & Kriging                          & Online  & AEI & CCF (D)    & L-BFGS-B (D)& Informed & MARIN  & Highest- and lowest-fidelity only used & \cite{huang2006sequential}\\    
                & RCK-EI  & Regression co-Kriging            & Online  & EI & LHS (S)    & GA (S)       & Free        & AFRL   & Highest- and lowest-fidelity only used & \cite{RCKthelen2020}\\
                & SRBF-LCB& Stochastic radial basis functions& Online  & LCB & CCF$^*$ (D)& DPSO+LS (D) & Free     & CNR-INM& Regressive for MF6 only & \cite{pellegrini2022multi,wackers2022efficient}\\
                & PCE     & Polynomial chaos expansion       & Offline & --  & SGP (D)    & --          & --       & AFRL   & Used for global accuracy only (NRMSE) & \cite{Bryson2017b} \\
                & SPCE    & Sparse polynomial chaos expansion& Offline & --  & SGP (D)    & --          & --       & AFRL   & Used for global accuracy only (NRMSE) & \cite{Rumpfkeil2020b} \\
&SPCE\_KG-GLRD&Sparse polynomial chaos expansion with Kriging& Online  & GLRD& FFD (D)    & --          & --       & AFRL   & Used for global accuracy only (NRMSE) & \cite{Rumpfkeil2022a} \\
\bottomrule
\multicolumn{11}{p{20cm}}{{Note: D/S = Deterministic/Stochastic behavior; metric medians and standard deviations are given for stochastic methods.}}\\
\multicolumn{11}{p{20cm}}{$^\dag$ {AEI = augmented expected improvement, EC = evolution control, EI = expected improvement, EIM = expected improvement + function minimum, GLRD = global/local regression discrepancy, LCB = lower confidence bounding, MVES = maximum value entropy search, MSES = multi-source entropy search, PI = probability of improvement.}}\\
\multicolumn{11}{p{20cm}}{$^\ddag$ {For online methods is the initial DoE: CCF = face-centered central composite design ($^*$only on the lowest fidelity, center of the domain for the others), FFD = factorial fractional design, HSS = Hammersley sequence sampling, HSOS = Halton sampling with owen-scramble LHS = Latin hypercube sampling, RBP = random baseline perturbation, SGP = sparse Gauss-Patterson grid.}}\\
\multicolumn{11}{p{20cm}}{$^\mathsection$ {Optimization algorithm used to solve Eq. \ref{e:OptimizationProblem}; for surrogate-rich methods, the algorithm is applied on surrogate prediction.}}\\
\multicolumn{11}{p{20cm}}{$^{\ddag\ddag}$ {Use of derivatives by the optimization algorithm.}}\\
\multicolumn{11}{p{20cm}}{$^*$ {None of the surrogates use derivatives to construct the function approximation.}}\\
\end{tabular}
\end{table}
\end{landscape}

\section{Performance Assessment Metrics}\label{sec:metrics} 
The proposed metrics are selected to offer a comprehensive framework and to be able to compare a broad spectrum of different goal-driven multifidelity methods regardless of their specific features. 
Considering the initialization criteria recommended in Subsection~\ref{s:Init}, the performance of the different multifidelity methods will be assessed through statistics over the set of all the convergence histories. 

Two types of metrics are defined: goal-sensitive and goal-insensitive.  Goal-sensitive metrics evaluate the accuracy of the optima $\bfx^\star$ computed with MF approximations, whereas goal-insensitive metrics address the global accuracy of MF approximations over the design space. The ability to compute these metrics is very dependent on the nature of the design space, the benchmark complexity, and the methods employed for optimization. Construction and use of these metrics thus involve compromises, which are hopefully struck in a manner that hits a sweet spot of generality, usefulness, and feasibility. 

With regard to goal-insensitive metrics, the global accuracy of MF approximations can be well interrogated for small $D$ and low benchmark complexity. As $D$ and benchmark complexity increase, the quantification of global accuracy rapidly becomes untenable through the effects of the \textit{curse of dimensionality} (CoD) \cite{bellman1957dynamic}. 

In contrast, global accuracy is not explicitly measured by goal-sensitive metrics, which simply assess computed optima. However, computing optima with global methods is also afflicted by the CoD, as is the global approximation of the functions that are optimized (whose cost can be mitigated by reduced sampling in unproductive areas). 

Central to quantifying accuracy is a scaling of the design space. By scaling variables to equivalent ranges of variation (on the unit hypercube), the influence of parameter sensitivities can be balanced and the relative accuracy of computed optima better characterized. Unless stated otherwise, scaling of each design parameter is performed linearly between lower and upper limits {\it supplied with the benchmark definition}:  
\beq
\tilde{x}_k \equiv {x_k-l_k \over u_k-l_k}\ \ \ (k=1,...,D),
\eeq
or $\tilde{\bfx} = \bf{S}\left(\bfx-{\bf l}\right)$, where $\bf{S}$ is a diagonal scaling matrix and $\bf l$ is the lowest valued corner of the un-scaled design space. For notational convenience, the tilde notation is neglected unless needed for clarity. Lower and upper limits are provided for scaling and simplification, even when a parameter might potentially be unbounded. Scaling of the design space is independent of constraint surfaces existing in the benchmark problem.  

Goal-insensitive metrics are defined first, which do not require knowledge of $\bfx^\star$. There are several competing goals for defining these metrics:
\begin{itemize}
\item {\it High consistency}, all methods should be evaluated in the same way;
\item {\it Low bias}, the evaluation strategy should not be biased;
\item {\it High utility}, the evaluation of metrics should be highly accurate and informative; 
\item {\it Affordable cost}, computing metrics should be affordable.
\end{itemize}

The goal-insensitive metric $\mathcal{E}_{\rm RMSE}$ records the  root-mean-squared error between the highest fidelity model of the objective $f$, to an approximation to $f$, $\widehat{f}$, over the design space:
\beq
   \mathcal{E}_{\rm RMSE} \equiv {1\over f_{\max}-f_{\min}} {\sqrt{\frac{1}{S}{\sum_{i=1}^{S}}(f(\bfx_i)-\widehat{f}(\bfx_i))^2}},
\label{e:RMSE}
\eeq

where $i$ is the sample index, $S$ is the number of samples, and $f_{\min}$ and $f_{\max}$ are minimum and maximum values of $f$ observed over the training data:
\beq 
\label{e:fminmax}
 f_{\min} \equiv \min_{\bfx_i} f(\bfx_i), \ \ \ 
 f_{\max} \equiv \max_{\bfx_i} f(\bfx_i).
\eeq

The sampling plan attempts to achieve two compromises: balancing computational cost with accuracy and balancing consistency with low bias. The design spaces corresponding to the (inexpensive) analytical benchmarks are exhaustively sampled to rigorously assess method accuracy and performance for different functions. As benchmark complexity is increased, the amount of sampling can be reduced.   Global accuracy cannot be computed for high-fidelity-based computational solutions, for their obvious high cost of sampling, but can be expected that MF-method strengths and weaknesses, in the context of global accuracy, can be extensively addressed with the analytical benchmarks proposed.

Goal-sensitive metrics are divided between metrics for the location of the optimum point (design accuracy) and metrics for the objective value (goal accuracy). The metrics are expressed as errors, which incorporate information about the optima assumed to be known {\it apriori}.

Three error metrics are defined for the benchmark problems, where optima are typically known analytically. These characterize normalized error in the design space, the objective function, and Euclidean distance in the normalized $\bfx$-$f$ hyperspace, respectively: 
\beq
\label{e:error_x}
\mathcal{E}_{x} \equiv \dfrac{\|\hat{\bfx}^\star-\bfx^\star\|}{\sqrt{N}},
\eeq

\beq
\label{e:error_f}
\mathcal{E}_{f} \equiv {f(\hat{\bfx}^\star)-f_{\min} \over f_{\max}-f_{\min}},
\eeq
\beq
\label{e:error_t}
\mathcal{E}_{t} \equiv \sqrt{\frac{\mathcal{E}_{x}^{2}+\mathcal{E}_{f}^{2}}{2}},
\eeq
where $\bfx$ is the {\it scaled} array of design variables on the unit hypercube, $\hat{\bfx}^\star$ is the location of the optimum of the approximation to $f$ (in scaled parameters), and $f_{\min}$ and $f_{\max}$ are the minimum and maximum of $f$ in $\mathcal{A}$. The error metrics $\mathcal{E}_{x}$ and $\mathcal{E}_{f}$ evaluate the design and goal accuracy, whereas the aggregated metric $\mathcal{E}_{t}$ introduced by Serani et. al.~\cite{serani2016-ASOC} evaluates accuracy in both senses, which can be of heightened significance when optima are in very flat or very peaky portions of the design space. Like the limits on design variables, $f_{\min}$ and $f_{\max}$ are provided with the benchmark cases to ensure consistent application by different MF methods. Reference values for the evaluation of the goal-sensitive metrics are summarized in Tab. \ref{t:referv}. It should be emphasized that Eq. (\ref{e:error_f}) utilizes an evaluation of the model of the highest-fidelity level at $\hat{\bfx}^\star$, a quantity computed with the MF approximation to $f$.  Evaluation of the full-fidelity model is favored over that of the approximation to better characterize the true objective realized through the MF approximation.  

\begin{table}[h]
\caption{Summary of multifidelity benchmark values for the evaluation of the goal-sensitive metrics}\label{t:referv}
\scriptsize
\centering
\begin{tabular}{lccccc}
\toprule
\bf Function & \bf MF ID & $D$  & $\mathbf{x}^\star$ & $f_{\min}$ & $f_{\max}$ \\ 
\midrule 
Forrester     & MF1.1 & 1  & 7.5752E-1       & -6.0207E-0 & 1.5830E+1 \\
Jump Forrester & MF1.2 & 1  & 1.4260E-1       & -9.8630E-1 & 2.5830E+1 \\
\midrule
              & MF2.1 & 2  & \{1.0000, 1.0000\}$^\mathsf{T}$ & 0.000E-0 & 3.6090E+3 \\
Rosenbrock    & MF2.2 & 5  & \{1.0000, \dots, 1.0000\}$^\mathsf{T}$ & 0.000E-0 & 1.4436E+4 \\
              & MF2.3 & 10 & \{1.0000, \dots, 1.0000\}$^\mathsf{T}$ & 0.000E-0 & 3.2481E+4 \\
\midrule
\multirow{3}{*}{Shifted-rotated Rastrigin}
              & MF3.1 & 2  & \{0.1000, 0.1000\}$^\mathsf{T}$ & 0.000E-0 & 4.0200E-0 \\
              & MF3.2 & 5  & \{0.1000, \dots, 0.1000\}$^\mathsf{T}$ & 0.000E-0 & 1.0050E+1 \\
              & MF3.3 & 10 & \{0.1000, \dots, 0.1000\}$^\mathsf{T}$ & 0.000E-0 & 2.0100E+1 \\
\midrule
              & MF4.1 & 1  & 2.7550E-1       & -6.2500E-1 & 3.6151E-1\\
Heterogeneous & MF4.2 & 2  & $\forall x_1=$ 0.000E-0 & -5.6271E-1 & 1.8350E-0\\
              & MF4.3 & 3  & $\forall x_1=$ 0.000E-0 & -5.6271E-1 & 4.3594E-0\\
\midrule
Springs       & MF5.1 & 2  & \{2.4674, 2.1932\}$^\mathsf{T}$ & -1.0000E-0 & 1.0000E-0   \\
Springs-masses & MF5.2 & 4  & \{1.0000, 3.9460, 4.0000, 3.2863\}$^\mathsf{T}$ & -1.0000E-0 & 1.0000E-0    \\
\midrule
Pacioreck     & MF6   & 2  & $\forall x_1x_2=2/[(3+j)\pi]$ with $j=0,4$ & -1.0000E-0 & 1.0000E-0    \\ 
\bottomrule
\end{tabular}
\end{table}

\section{Results and Discussion}\label{sec:results}
While most methods were used to find the global optimum in the design space and the NRMSE is just computed as a byproduct, the KG-GLRD, SPCE-KG-GLRD, PCE, and SPCE methods were only used to build globally accurate surrogates. Thus, a direct comparison in the NRMSE performance is not necessarily fair. Note that all methods are intended to be an MF approach and no comparison with their single-fidelity (high-fidelity) counterpart is given since it is beyond the scope of this work. 

In the following, the results for each benchmark problem will be shown and discussed. Note that lines are plotted between markers as a visual aid, but do not necessarily reflect performance at intermediate points. For the statistical assessment of stochastic methods, ten runs are conducted and results are provided in terms of the median (or mean if the median was not provided by a group) and standard deviation. 

The x-axis always shows the percentage (0 to 100) of the total computational budget given in Table~\ref{t:ExSetup}. The y-axis shows the percentage error in the four metrics (from left to right: variable, function, and global distances as well as NRMSE) discussed previously. The top row shows the median (or mean) and the bottom row displays the standard deviation. Obviously, no statistics are given for deterministic approaches.  

Before discussing benchmark by benchmark, some comments are required to better understand the results of the BFGS. Since this method is gradient-based only local optima can be reliably found which is not a problem if the function is not multi-modal (such as Rosenbrock in 2D). However, when the function is multi-modal (e.g., the Rastrigin function, the spring-mass problem, or the Rosenbrock function in 5D) the acquired optimum depends on the starting point. Here, the center of the domain and nine other Latin-hypercube sampled starting points are employed and the mean and standard deviations of these ten optimization runs are displayed. If, for example, nine of these runs found the global optimum correctly but one of the runs converges to a local optimum (possibly far away in both variable and function space from the global optimum) then the statistics would be dominated by this one run and the convergence history would look fairly flat. In addition, the cost of BFGS is not controlled directly but is rather dictated by the number of required iterations and line-search evaluations in each iteration. Thus, the cost of each iteration varies and whenever a run converges it is not included in the averaging procedure any longer. 

\begin{figure}[H]
\centering
\includegraphics[width=0.9\textwidth]{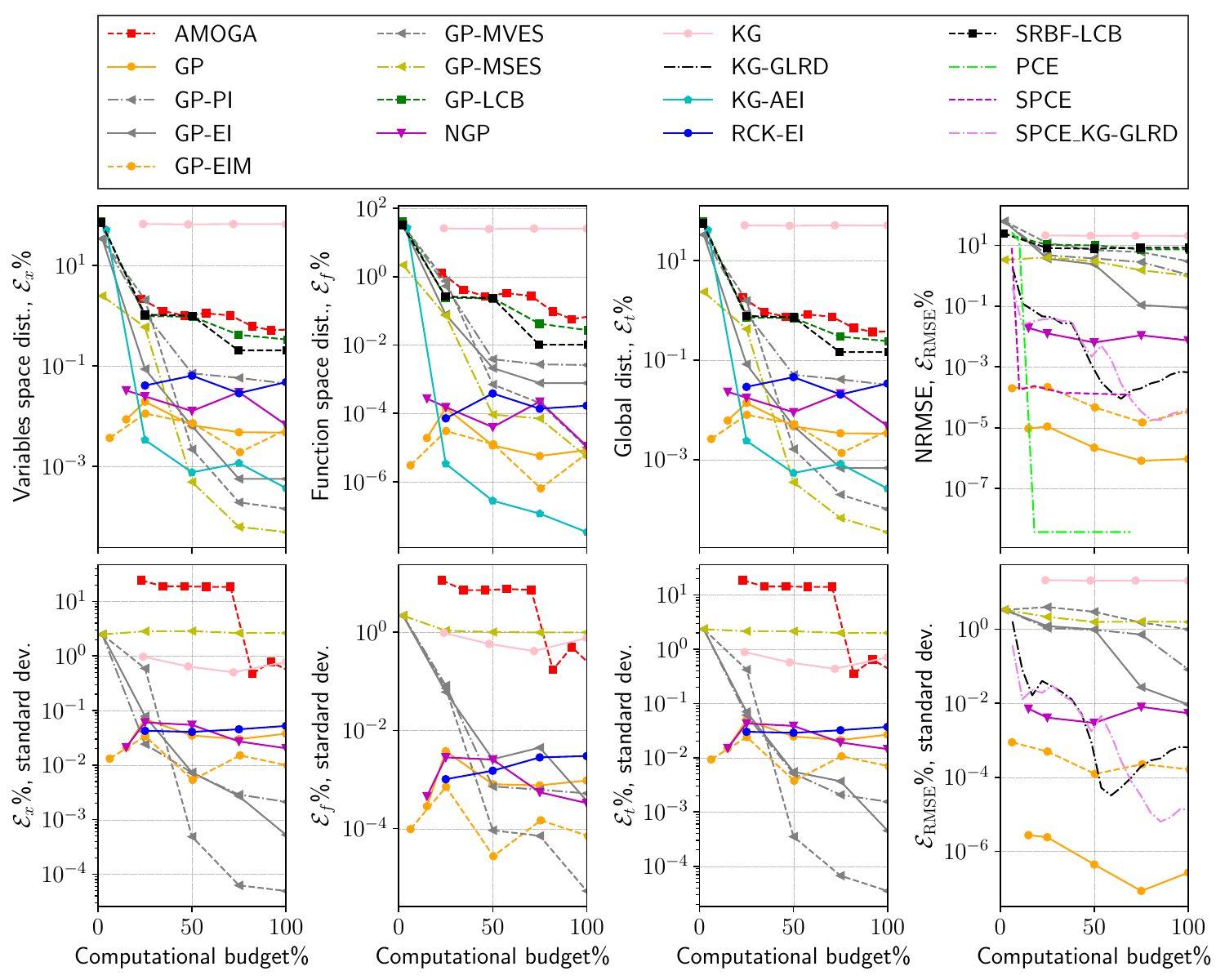}
\caption{MF1.1 metrics convergence}
\label{fig_ch06:l1-mf1.1}
\centering
\includegraphics[width=0.9\textwidth]{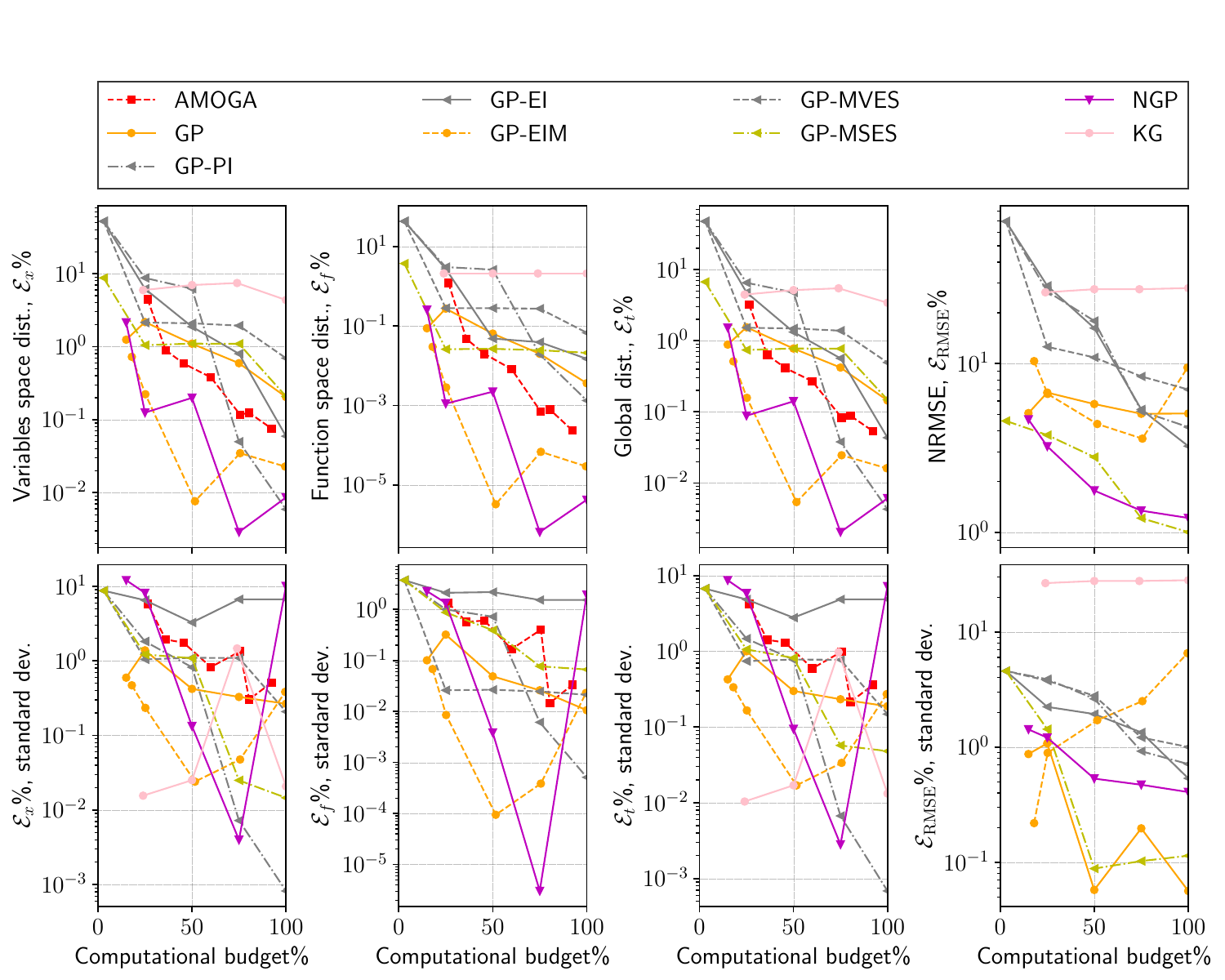}
\caption{MF1.2 metrics convergence}
\label{fig_ch06:l1-mf1.2}
\end{figure}

The Forrester and Jump Forrester results are shown in Figures~\ref{fig_ch06:l1-mf1.1} and \ref{fig_ch06:l1-mf1.2}, respectively. Looking at the optimization of the Forrester function (without jump) RCK-EI and GP-EIM are performing best whereas a globally accurate surrogate of the function is best represented by PCE, GP, and again GP-EIM. For the Jump Forrester function, GP-EIM and NGP are performing relatively well for the optimization. However, none of the methods give a particularly good globally accurate model but NGP and GP-MSES are doing better than the other methods.

\begin{figure}[H]
\centering
\includegraphics[width=0.9\textwidth]{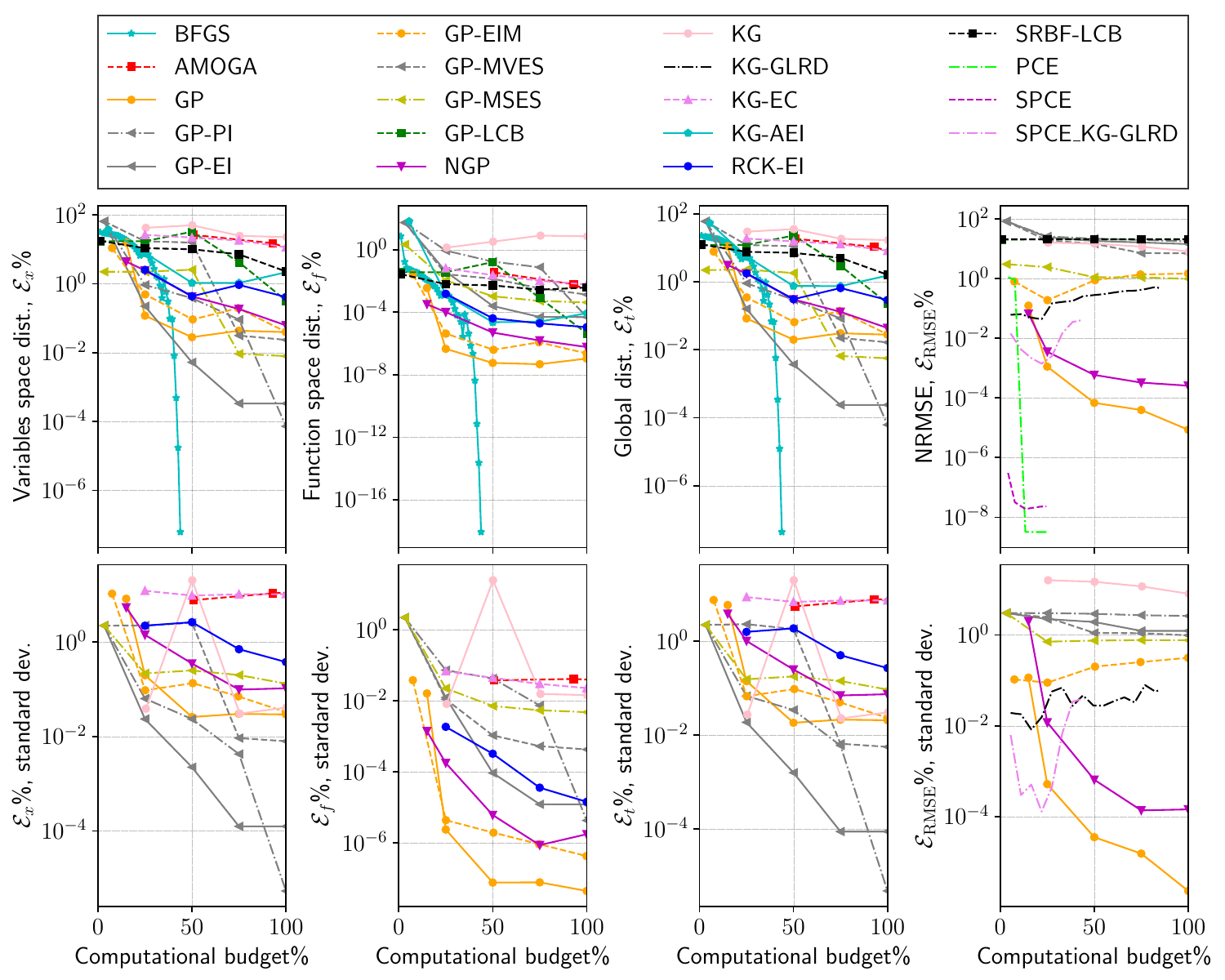}
\caption{MF2.1 metrics convergence}
\label{fig_ch06:l1-mf2.1}
\centering
\includegraphics[width=0.9\textwidth]{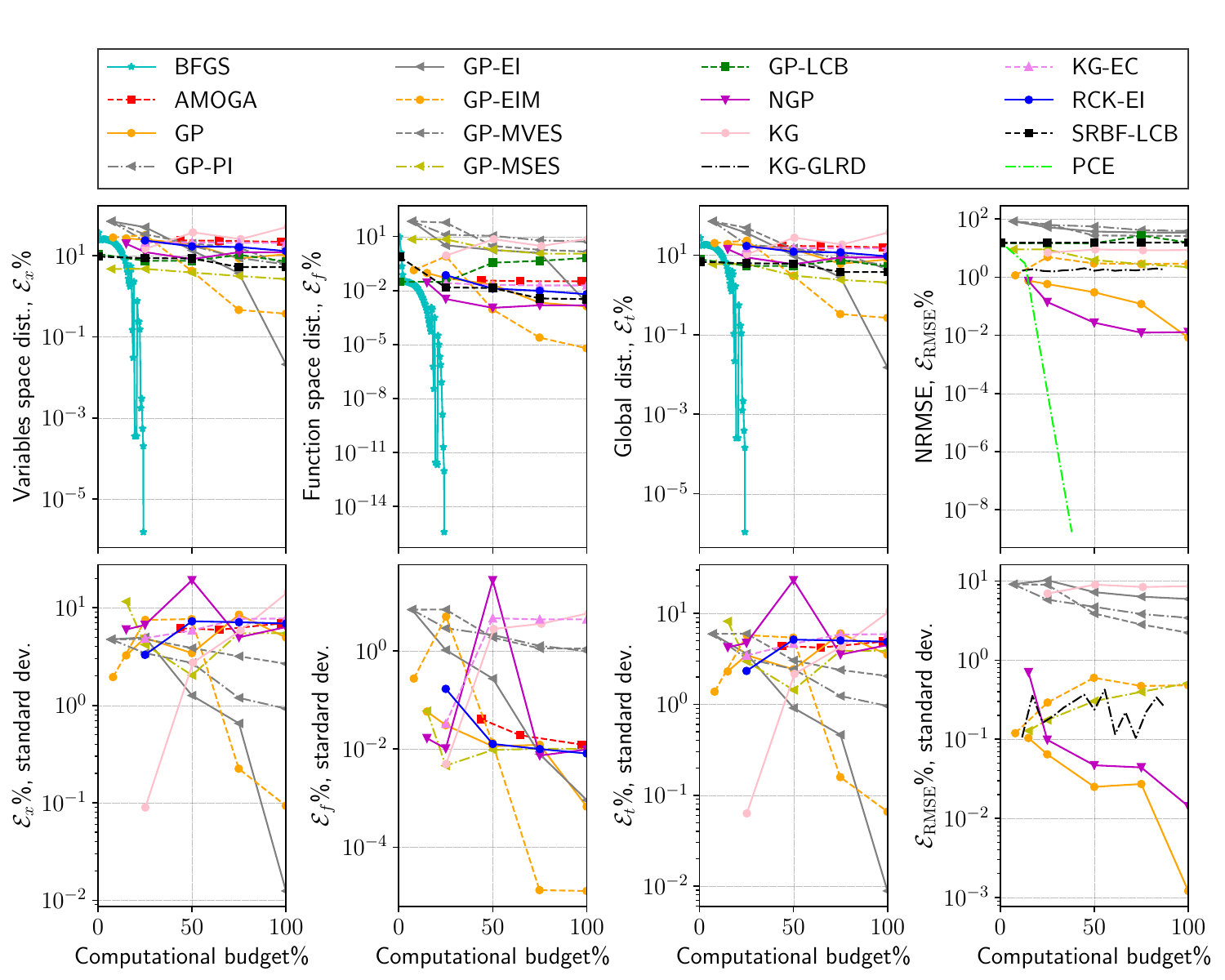}
\caption{MF2.2 metrics convergence}
\label{fig_ch06:l1-mf2.2}
\end{figure}

\begin{figure}[H]
\centering
\includegraphics[width=0.9\textwidth]{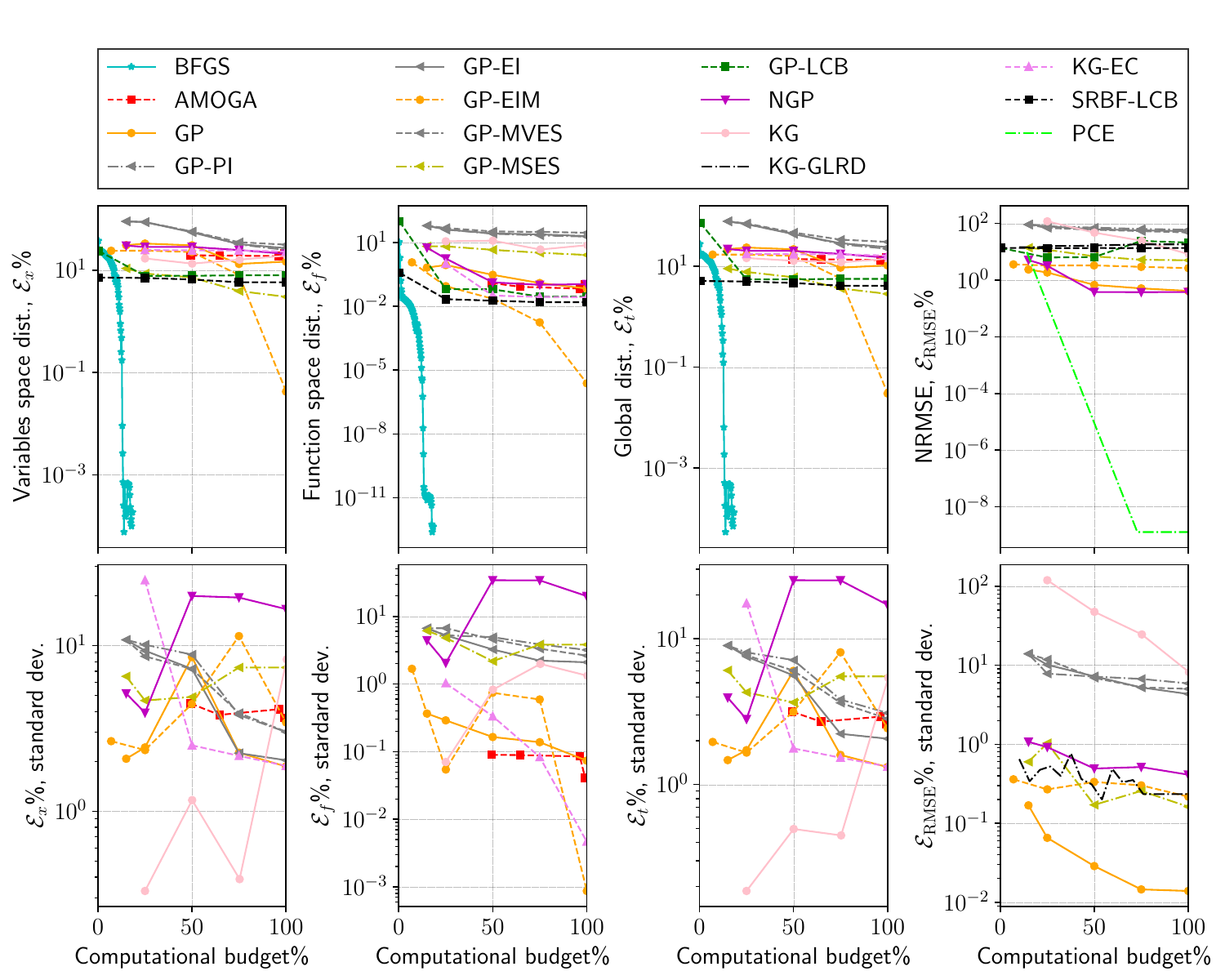}
\caption{MF2.3 metrics convergence}
\label{fig_ch06:l1-mf2.3}
\centering
\includegraphics[width=0.9\textwidth]{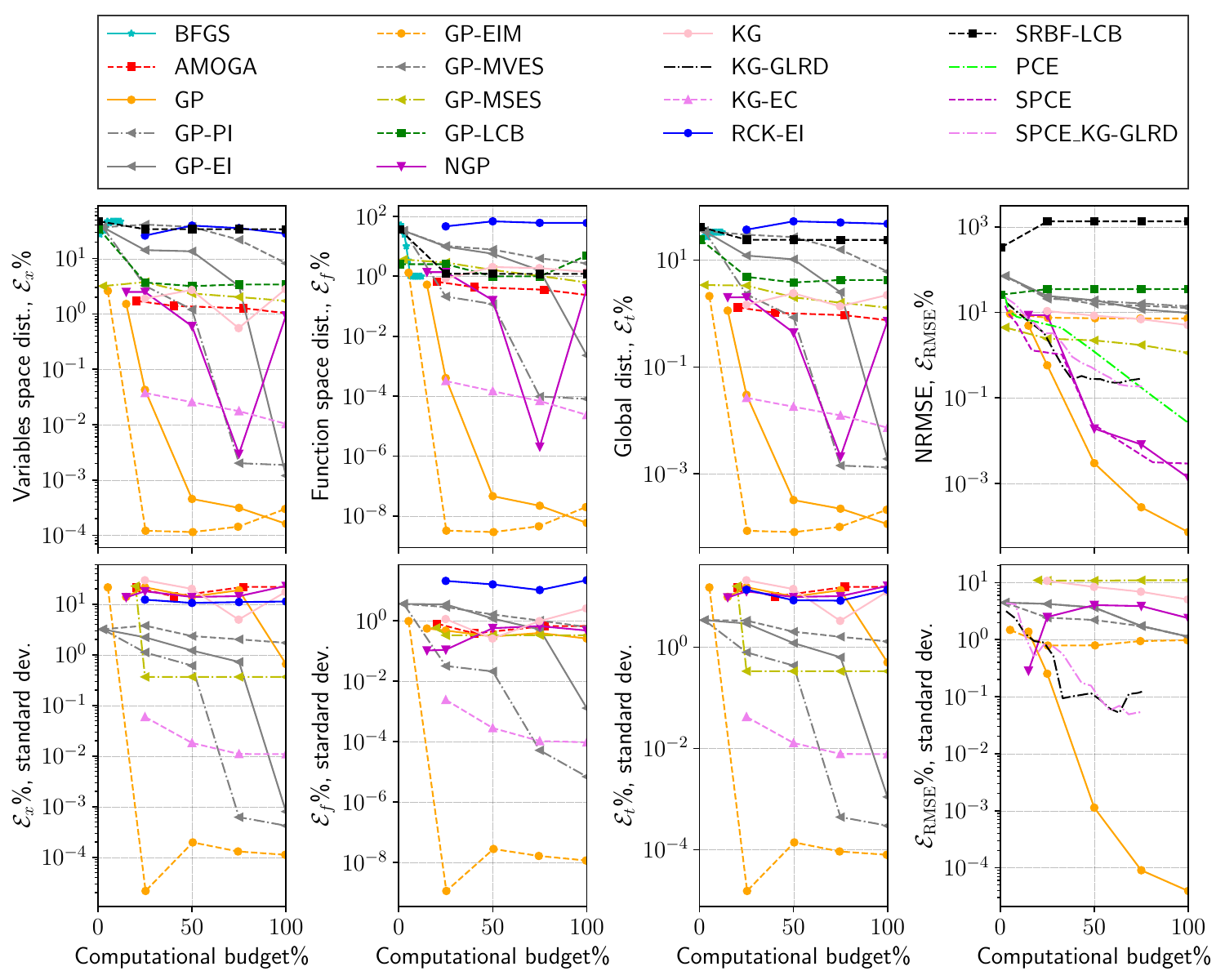}
\caption{MF3.1 metrics convergence}
\label{fig_ch06:l1-mf3.1}
\end{figure}

The Rosenbrock function results are shown in Figures~\ref{fig_ch06:l1-mf2.1} to \ref{fig_ch06:l1-mf2.3} for 2, 5, and 10 dimensions, respectively. This function has only one optimum and thus the statistics for MF BFGS converge well compared to the surrogate-rich methods, especially in higher dimensions. Among the surrogate-rich methods, GP-EIM is performing best over all three dimensions considered here.

\begin{figure}[H]
\centering
\includegraphics[width=0.9\textwidth]{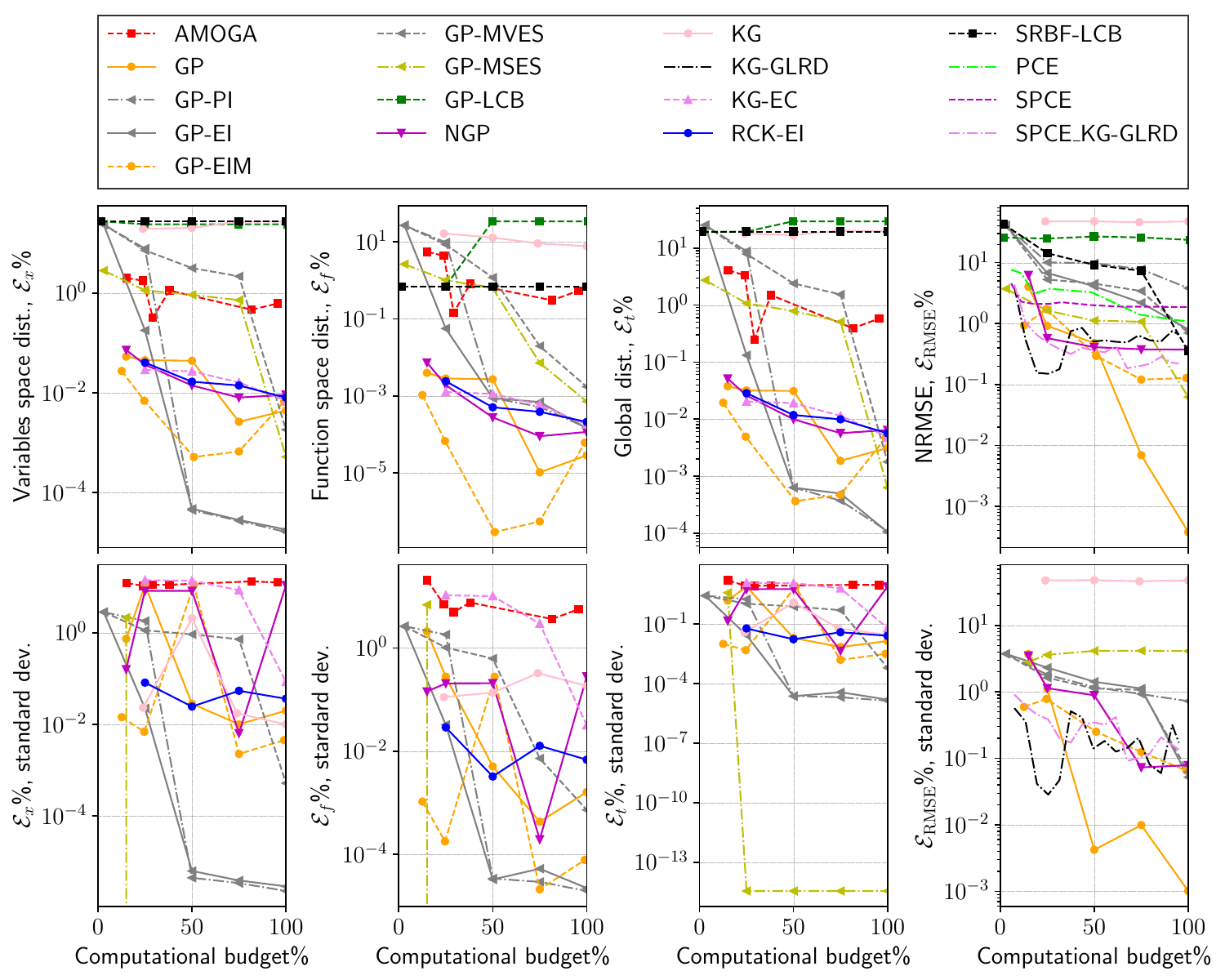}
\caption{MF4.1 metrics convergence}
\label{fig_ch06:l1-mf4.1}
%
\centering
\includegraphics[width=0.9\textwidth]{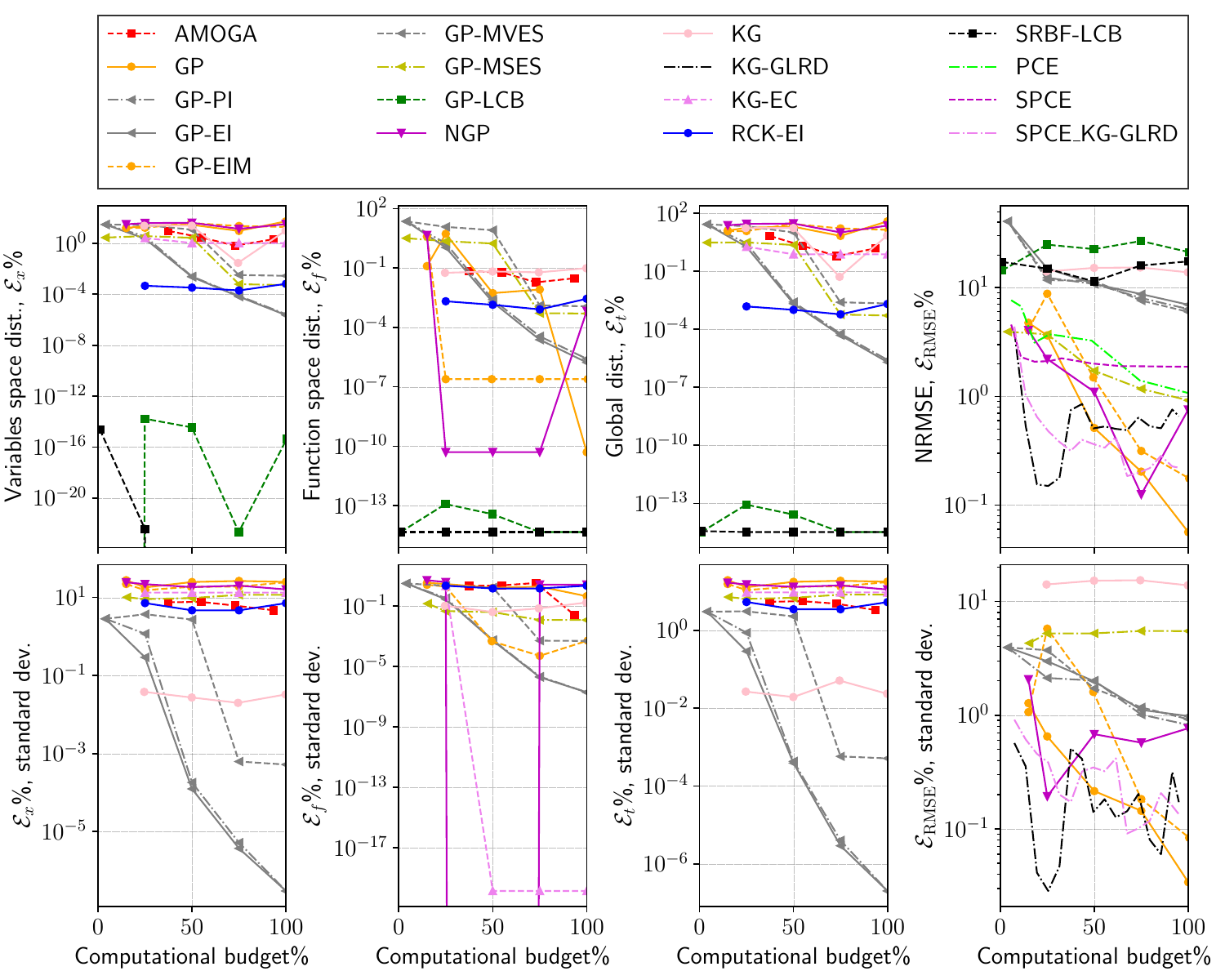}
\caption{MF4.2 metrics convergence}
\label{fig_ch06:l1-mf4.2}
\end{figure}

\begin{figure}[H]
\centering
\includegraphics[width=0.9\textwidth]{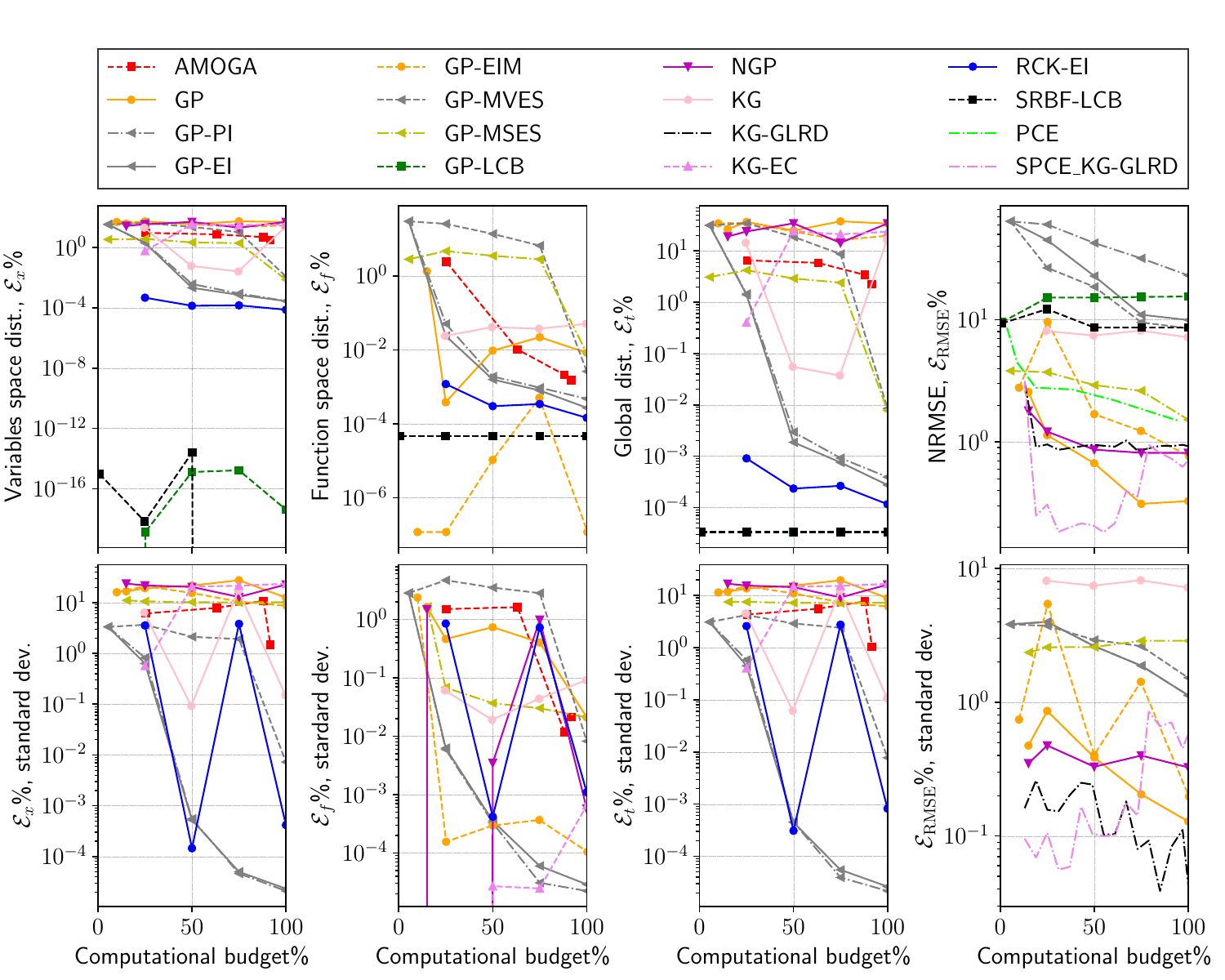}
\caption{MF4.3 metrics convergence}
\label{fig_ch06:l1-mf4.3}
\centering
\includegraphics[width=0.9\textwidth]{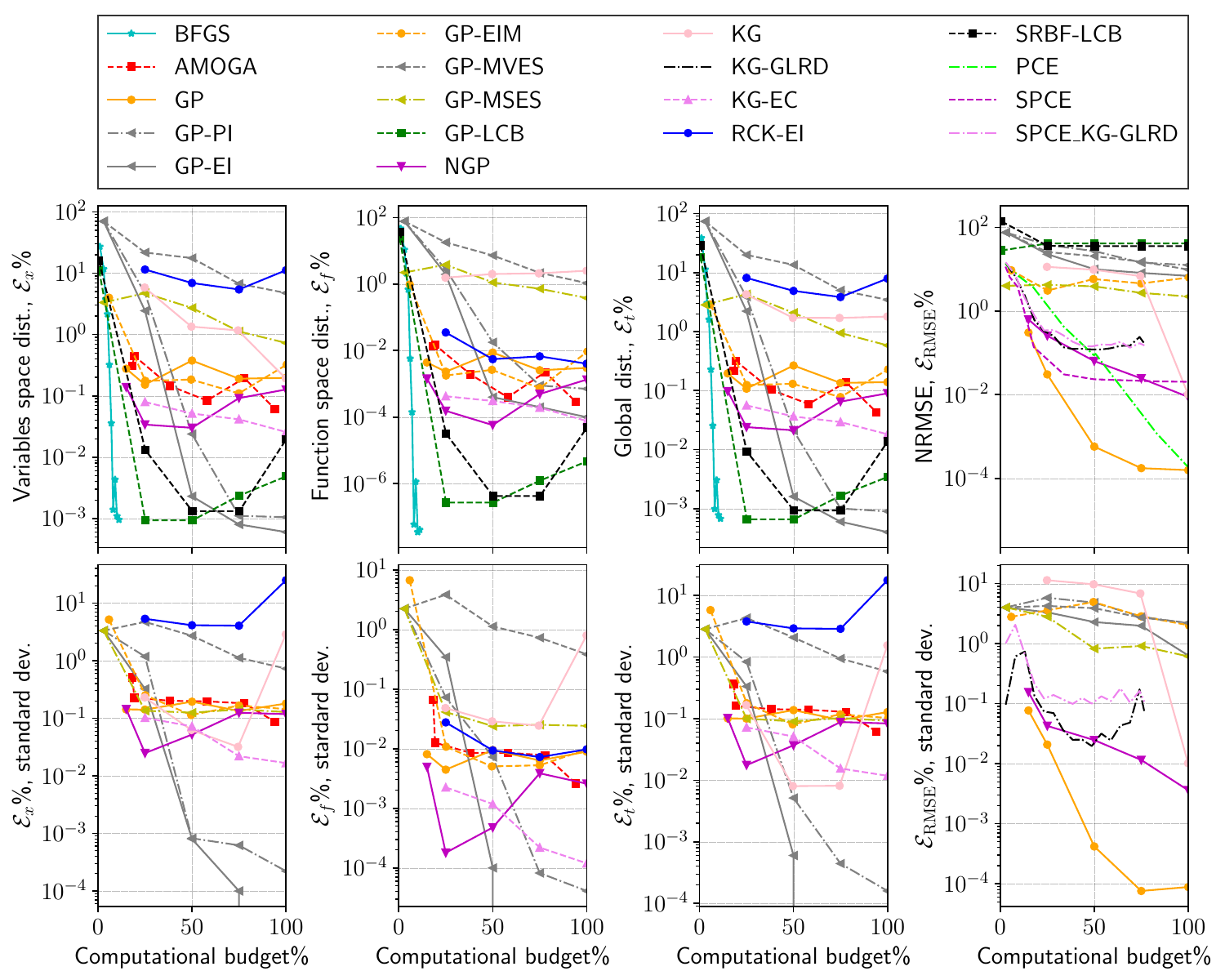}
\caption{MF5.1 metrics convergence}
\label{fig_ch06:l1-mf5.1}
\end{figure}

For a globally accurate representation, once sufficient training data is gathered to construct a fourth-order polynomial expansion, the PCE and SPCE models should match the truth function exactly with the error dropping accordingly. The full-order PCE surrogate model with HF training data alone requires 30 training points to match the exact function in two dimensions and 2002 points in ten dimensions due to an oversampling ratio of two. For MF models some of these points may be low-fidelity and thus computational savings can be achieved. This expected behavior can be seen in the figures. In general, sparse PCEs (SPCE) outperform full-order PCEs up to a point when the penalty coefficient, $\lambda$, introduces error into the fit and the RMSE cannot drop further as can be seen in the top right-hand corners of Figures~\ref{fig_ch06:l1-mf2.1}, \ref{fig_ch06:l1-mf3.1}, \ref{fig_ch06:l1-mf4.1}, \ref{fig_ch06:l1-mf4.2}, and \ref{fig_ch06:l1-mf5.1}.
Among the non-polynomial-based methods, GP and NGP are working the best to represent the function globally as indicated by relatively low NRMSE values.

\begin{figure}[H]
\centering
\includegraphics[width=0.9\textwidth]{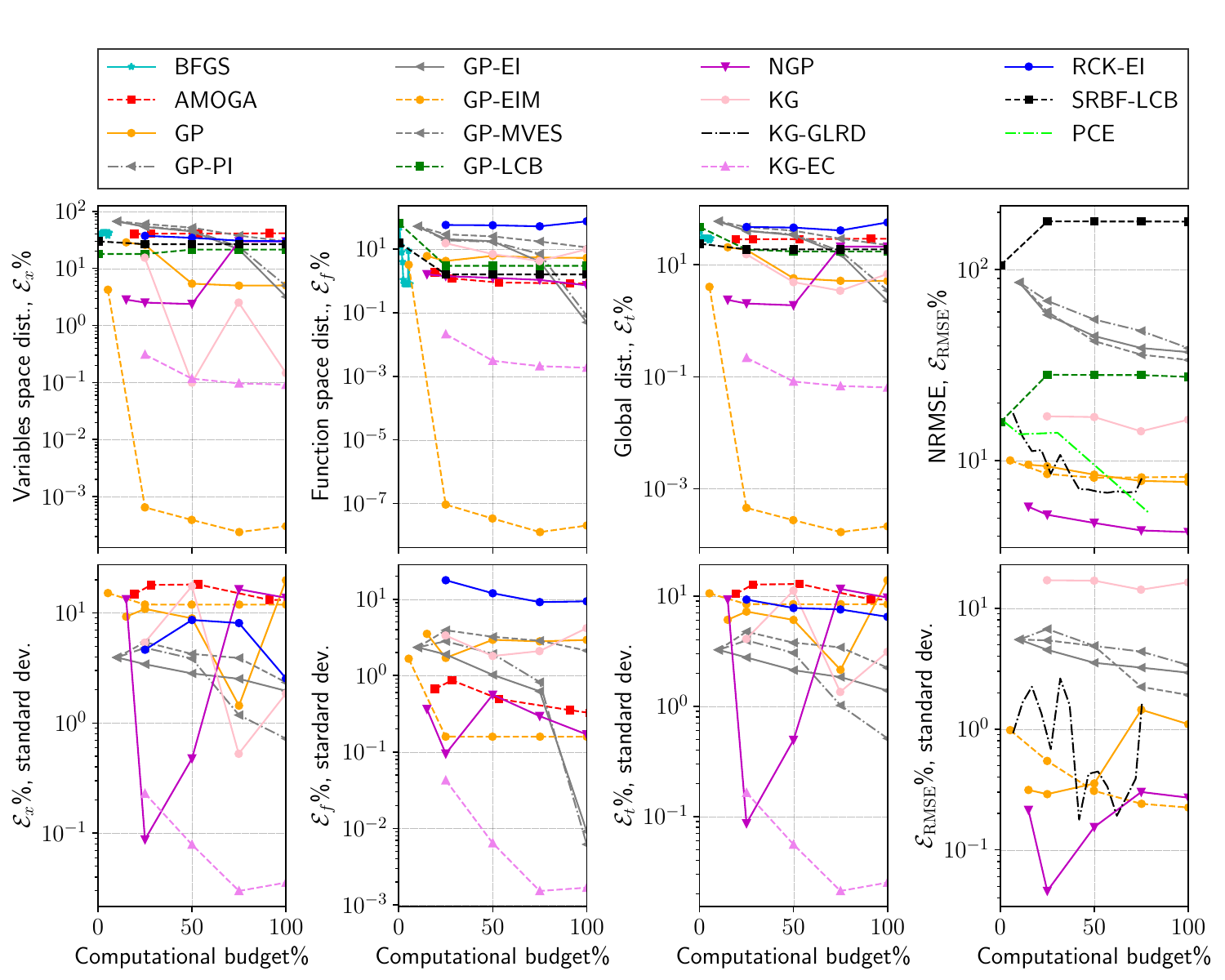}
\caption{MF3.2 metrics convergence}
\label{fig_ch06:l1-mf3.2}
\centering
\includegraphics[width=0.9\textwidth]{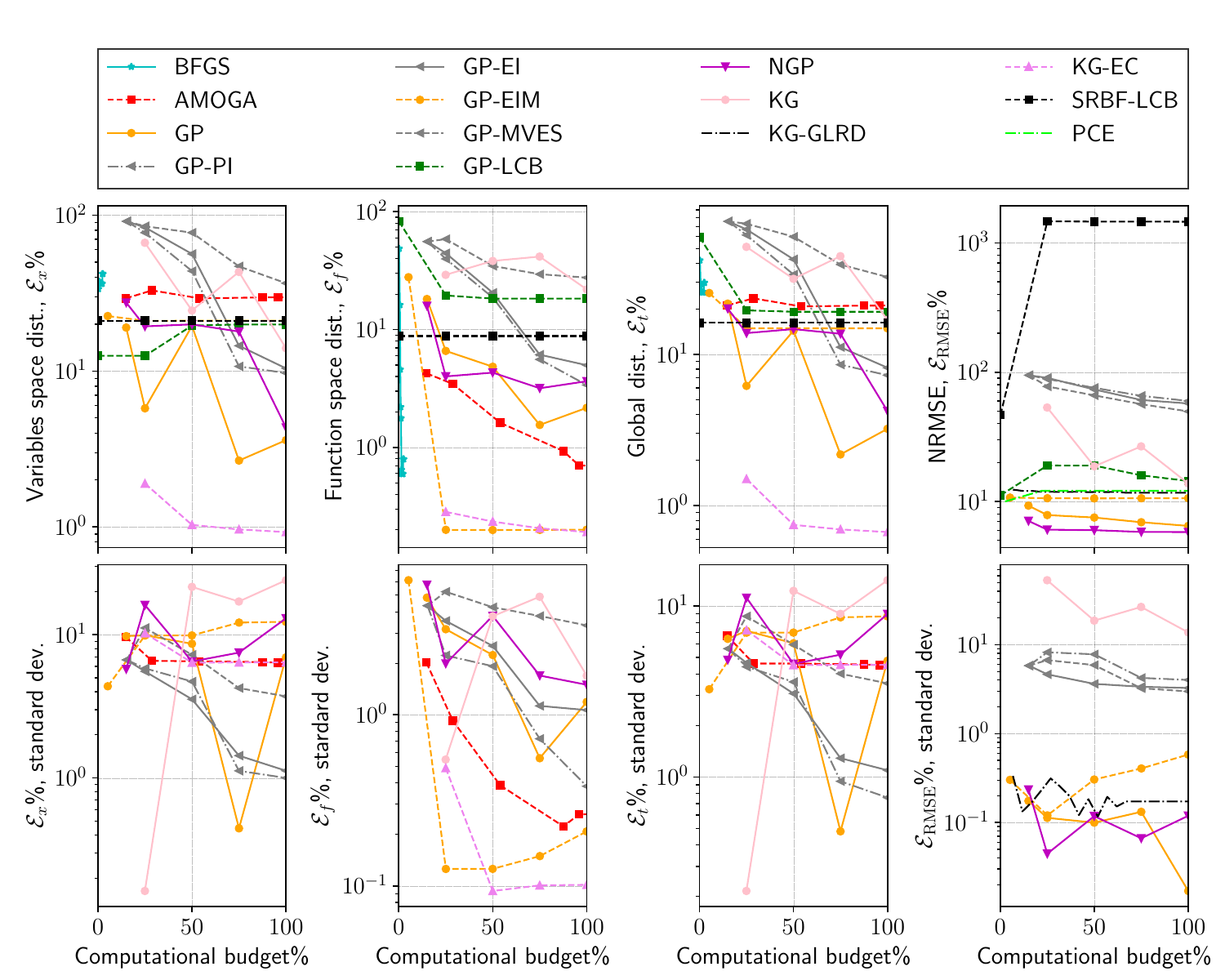}
\caption{MF3.3 metrics convergence}
\label{fig_ch06:l1-mf3.3}
\end{figure}

The Rastrigin function results are shown in Figures~\ref{fig_ch06:l1-mf3.1} to \ref{fig_ch06:l1-mf3.3} for 2, 5, and 10 dimensions, respectively. For the optimization, BFGS is doing well as are GP-EIM (for all three dimensions) and KG-EC (in 5D and 10D) and RCK-EI (in 2D). For building globally accurate surrogate models the best methods vary with the dimensionality of the problem. In 2D GP, NGP and SPCE are performing best. In 5D NGP, PCE and KG-GLRD seem to work well whereas in 10D NGP and GP are best though in reality most methods are close to each other and all convergence rates are very flat.

The Heterogeneous function results are shown in Figures~\ref{fig_ch06:l1-mf4.1} to \ref{fig_ch06:l1-mf4.3} for 1, 2, and 3 dimensions, respectively. The suitability of the two- and three-dimensional case for the assessment of optimization performances is somewhat questionable since the optimum lies on a boundary and is the entire boundary line in 2D. Depending on the initial DOE some methods may find that optimum right away as can be seen in Figure~\ref{fig_ch06:l1-mf4.2}. For the optimization in 1D GP-EIM is best with RCK-EI and GP-EI somewhat trailing behind. In 3D SRBF-LCB and GP-LCB are performing best. 

Looking at the global accuracy via NRMSE, in one dimension GP and GP-EIM are doing well, in two dimensions KG-GLRD is performing well early on, and at the end GP and GP-EIM. In three dimensions SPCE\_KG-GLRD and again GP and GP-EIM are best, however, all method's convergences are fairly flat.

\begin{figure}[H]
\centering
\includegraphics[width=0.9\textwidth]{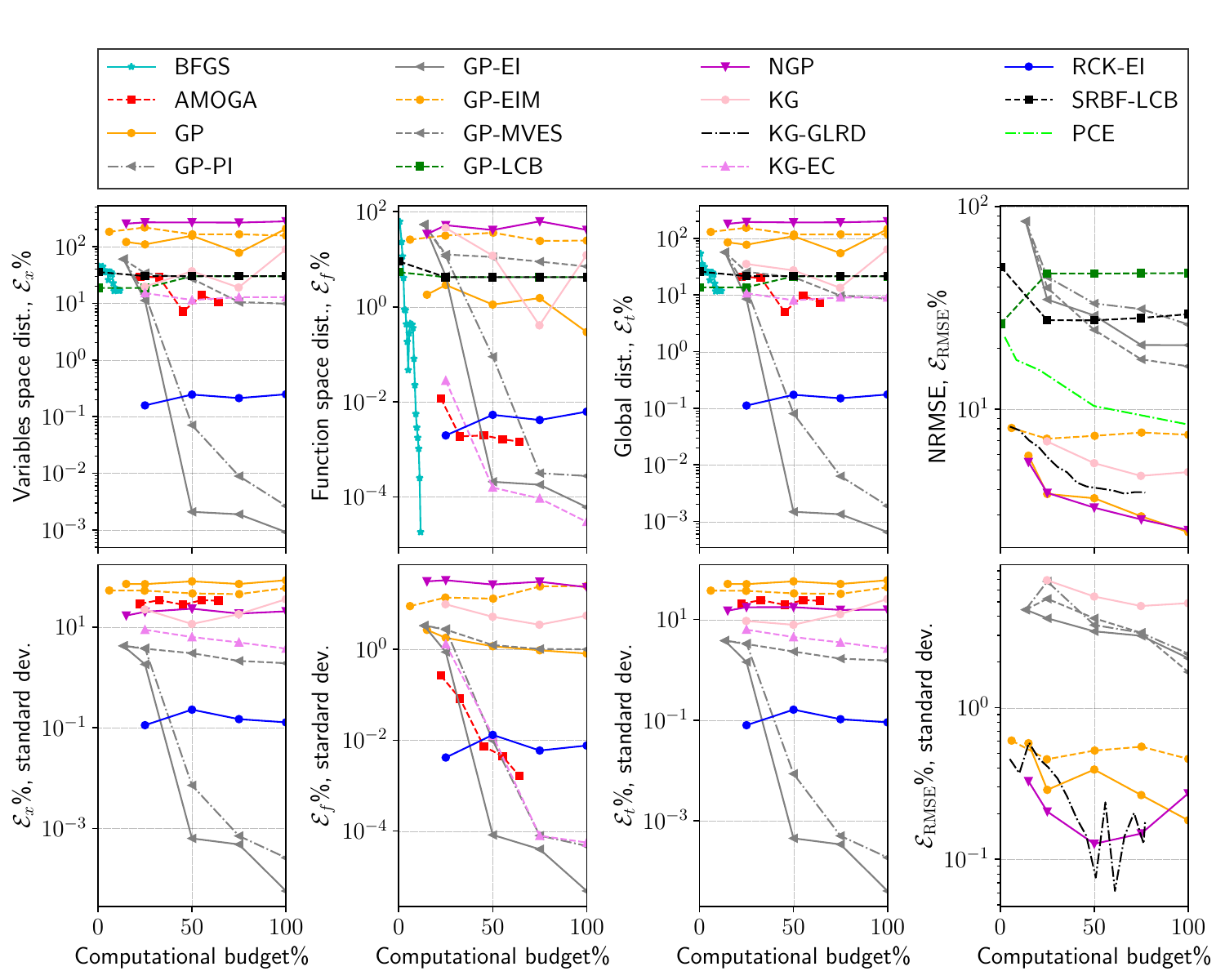}
\caption{MF5.2 metrics convergence}
\label{fig_ch06:l1-mf5.2}
\centering
\includegraphics[width=1\textwidth]{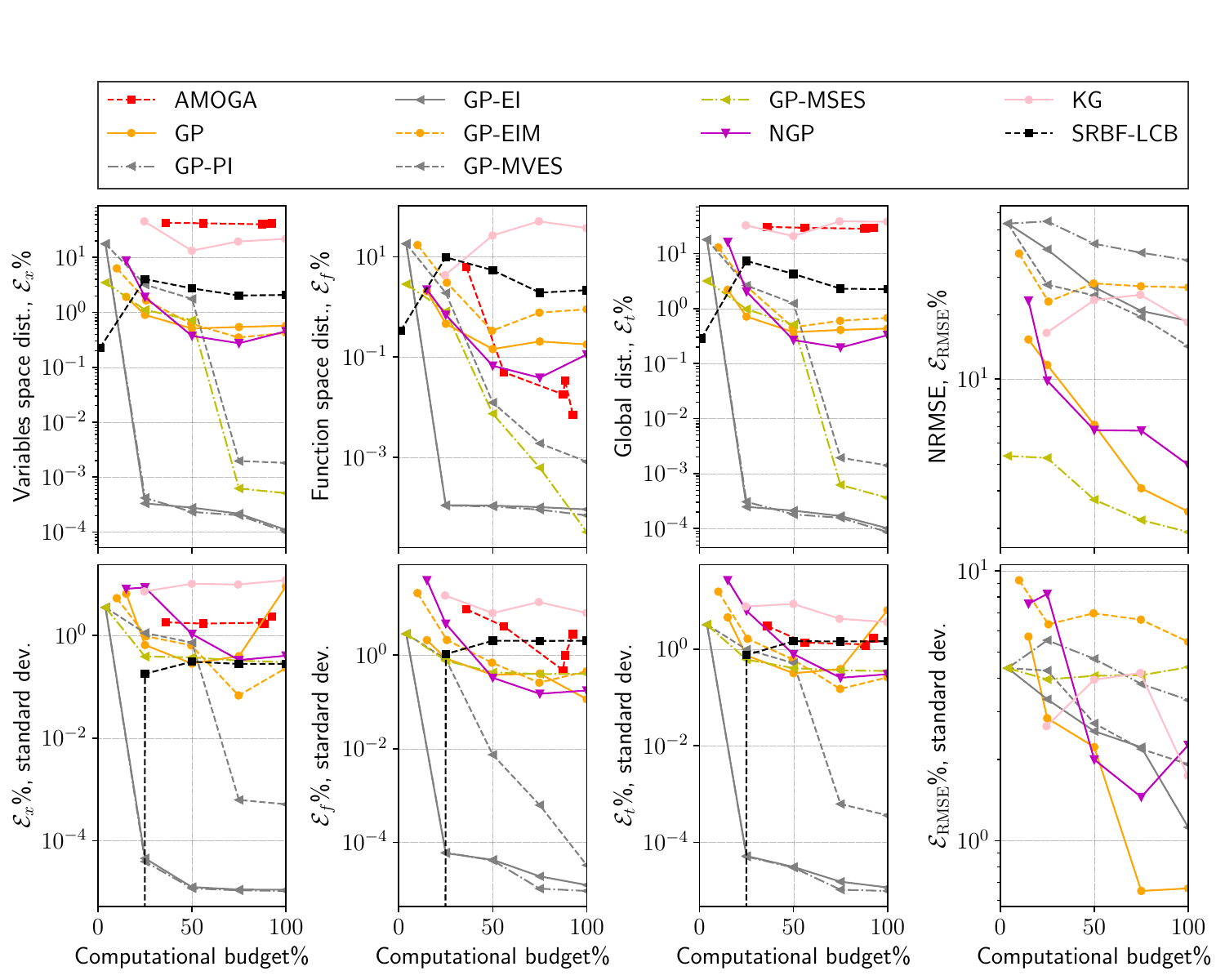}
\caption{MF6 metrics convergence}
\label{fig_ch06:l1-mf6}
\end{figure}

The Springs and Springs-masses results are shown in Figures~\ref{fig_ch06:l1-mf5.1} and \ref{fig_ch06:l1-mf5.2}, respectively. For the optimization, BFGS and GP-EI are doing well in both 2D and 4D as are GP-LCB and SRBF-LCB in 2D and RCK-EI and GP-PI in 4D. 
For global accuracy GP and PCE are looking good in 2D and NGP, GP and KG-GLRD are best in 4D, though the convergence of all methods is somewhat sluggish.

Lastly, the Pacioreck results are shown in Figure~\ref{fig_ch06:l1-mf6}. For the optimization GP-PI, GP-EI, and by the end, GP-MSES are performing the best whereas for the global accuracy GP-MSES, GP, and NGP are superior to the other methods. 

Table~\ref{tab_ch06:l1_app} gives a comprehensive overview of the results for all multi-fidelity L1 benchmark problems for all methods presented in Table~\ref{tab_ch06:l1_methods}.
Four different budget levels are highlighted (25\%, 50\%, 75\%, and 100\% of the total budget) and the percent error is color-coded (orange if greater than 10\%, yellow if between 1 and 10\%, green if between 0.1 and 1\%, light blue if between 0.01 and 0.1\% and lastly dark blue if less than 0.01\%). This should give the reader a quick visual overview of the various performances.

\begin{table}[H]
\centering
\caption{Overview of multi-fidelity methods results for L1 benchmark}
\label{tab_ch06:l1_app}
\tiny
\begin{tabular}{l
p{0.15cm}@{}p{0.15cm}@{}p{0.15cm}@{}p{0.15cm}@{}p{0.15cm}@{}p{0.15cm}@{}p{0.15cm}@{}p{0.15cm}@{}p{0.15cm}@{}p{0.15cm}@{}p{0.15cm}@{}p{0.15cm}@{}p{0.15cm}@{}p{0.15cm}
@{}p{0.3cm}@{}
p{0.15cm}@{}p{0.15cm}@{}p{0.15cm}@{}p{0.15cm}@{}p{0.15cm}@{}p{0.15cm}@{}p{0.15cm}@{}p{0.15cm}@{}p{0.15cm}@{}p{0.15cm}@{}p{0.15cm}@{}p{0.15cm}@{}p{0.15cm}@{}p{0.15cm}
@{}p{0.3cm}@{}
p{0.15cm}@{}p{0.15cm}@{}p{0.15cm}@{}p{0.15cm}@{}p{0.15cm}@{}p{0.15cm}@{}p{0.15cm}@{}p{0.15cm}@{}p{0.15cm}@{}p{0.15cm}@{}p{0.15cm}@{}p{0.15cm}@{}p{0.15cm}@{}p{0.15cm}
@{}p{0.3cm}@{}
p{0.15cm}@{}p{0.15cm}@{}p{0.15cm}@{}p{0.15cm}@{}p{0.15cm}@{}p{0.15cm}@{}p{0.15cm}@{}p{0.15cm}@{}p{0.15cm}@{}p{0.15cm}@{}p{0.15cm}@{}p{0.15cm}@{}p{0.15cm}@{}p{0.15cm}
c
}
\toprule
\bf Method & \multicolumn{14}{c}{$\mathcal{E}_x$} & &\multicolumn{14}{c}{$\mathcal{E}_f$} & &\multicolumn{14}{c}{$\mathcal{E}_t$} & &\multicolumn{14}{c}{$\mathcal{E}_{\rm RMSE}$} & \bf B\%\\
\midrule
\multirow{4}{*}{BFGS           }&\nmk&\nmk&\omk&\bmk&\bmk&\rmk&\rmk&\rmk&\nmk&\nmk&\nmk&\bmk&\rmk&\nmk&&\nmk&\nmk&\bmk&\bmk&\bmk&\ymk&\ymk&\ymk&\nmk&\nmk&\nmk&\bmk&\bmk&\nmk&&\nmk&\nmk&\omk&\bmk&\bmk&\rmk&\rmk&\rmk&\nmk&\nmk&\nmk&\bmk&\rmk&\nmk&&\nmk&\nmk&\nmk&\nmk&\nmk&\nmk&\nmk&\nmk&\nmk&\nmk&\nmk&\nmk&\nmk&\nmk&25\\
                                &\nmk&\nmk&\bmk&\nmk&\nmk&\nmk&\nmk&\nmk&\nmk&\nmk&\nmk&\nmk&\nmk&\nmk&&\nmk&\nmk&\bmk&\nmk&\nmk&\nmk&\nmk&\nmk&\nmk&\nmk&\nmk&\nmk&\nmk&\nmk&&\nmk&\nmk&\bmk&\nmk&\nmk&\nmk&\nmk&\nmk&\nmk&\nmk&\nmk&\nmk&\nmk&\nmk&&\nmk&\nmk&\nmk&\nmk&\nmk&\nmk&\nmk&\nmk&\nmk&\nmk&\nmk&\nmk&\nmk&\nmk&50\\
                                &\nmk&\nmk&\nmk&\nmk&\nmk&\nmk&\nmk&\nmk&\nmk&\nmk&\nmk&\nmk&\nmk&\nmk&&\nmk&\nmk&\nmk&\nmk&\nmk&\nmk&\nmk&\nmk&\nmk&\nmk&\nmk&\nmk&\nmk&\nmk&&\nmk&\nmk&\nmk&\nmk&\nmk&\nmk&\nmk&\nmk&\nmk&\nmk&\nmk&\nmk&\nmk&\nmk&&\nmk&\nmk&\nmk&\nmk&\nmk&\nmk&\nmk&\nmk&\nmk&\nmk&\nmk&\nmk&\nmk&\nmk&75\\
                                &\nmk&\nmk&\nmk&\nmk&\nmk&\nmk&\nmk&\nmk&\nmk&\nmk&\nmk&\nmk&\nmk&\nmk&&\nmk&\nmk&\nmk&\nmk&\nmk&\nmk&\nmk&\nmk&\nmk&\nmk&\nmk&\nmk&\nmk&\nmk&&\nmk&\nmk&\nmk&\nmk&\nmk&\nmk&\nmk&\nmk&\nmk&\nmk&\nmk&\nmk&\nmk&\nmk&&\nmk&\nmk&\nmk&\nmk&\nmk&\nmk&\nmk&\nmk&\nmk&\nmk&\nmk&\nmk&\nmk&\nmk&100\\
\hline
\multirow{4}{*}{AMOGA          }&\omk&\omk&\nmk&\rmk&\rmk&\omk&\rmk&\rmk&\omk&\omk&\omk&\ymk&\rmk&\rmk&&\omk&\omk&\nmk&\gmk&\ymk&\ymk&\omk&\omk&\omk&\gmk&\omk&\gmk&\gmk&\omk&&\omk&\omk&\nmk&\rmk&\rmk&\omk&\rmk&\rmk&\omk&\omk&\omk&\ymk&\rmk&\rmk&&\nmk&\nmk&\nmk&\nmk&\nmk&\nmk&\nmk&\nmk&\nmk&\nmk&\nmk&\nmk&\nmk&\nmk&25\\
                                &\omk&\ymk&\rmk&\rmk&\rmk&\omk&\rmk&\rmk&\omk&\omk&\omk&\ymk&\omk&\rmk&&\ymk&\gmk&\gmk&\gmk&\gmk&\ymk&\omk&\omk&\ymk&\gmk&\gmk&\bmk&\bmk&\gmk&&\ymk&\ymk&\rmk&\rmk&\rmk&\omk&\rmk&\rmk&\omk&\omk&\omk&\ymk&\omk&\rmk&&\nmk&\nmk&\nmk&\nmk&\nmk&\nmk&\nmk&\nmk&\nmk&\nmk&\nmk&\nmk&\nmk&\nmk&50\\
                                &\omk&\ymk&\rmk&\rmk&\rmk&\omk&\rmk&\rmk&\ymk&\ymk&\omk&\ymk&\rmk&\rmk&&\ymk&\bmk&\bmk&\gmk&\gmk&\ymk&\ymk&\omk&\ymk&\gmk&\bmk&\bmk&\bmk&\gmk&&\ymk&\gmk&\rmk&\rmk&\rmk&\ymk&\rmk&\rmk&\ymk&\ymk&\omk&\ymk&\omk&\rmk&&\nmk&\nmk&\nmk&\nmk&\nmk&\nmk&\nmk&\nmk&\nmk&\nmk&\nmk&\nmk&\nmk&\nmk&75\\
                                &\ymk&\gmk&\rmk&\rmk&\rmk&\omk&\rmk&\rmk&\ymk&\omk&\omk&\gmk&\nmk&\rmk&&\gmk&\bmk&\bmk&\gmk&\ymk&\ymk&\ymk&\ymk&\ymk&\gmk&\bmk&\bmk&\nmk&\bmk&&\ymk&\gmk&\omk&\rmk&\rmk&\ymk&\rmk&\rmk&\ymk&\omk&\omk&\gmk&\nmk&\rmk&&\nmk&\nmk&\nmk&\nmk&\nmk&\nmk&\nmk&\nmk&\nmk&\nmk&\nmk&\nmk&\nmk&\nmk&100\\
\hline
\multirow{4}{*}{GP             }&\gmk&\omk&\ymk&\rmk&\rmk&\gmk&\rmk&\omk&\gmk&\rmk&\rmk&\ymk&\rmk&\ymk&&\bmk&\ymk&\bmk&\gmk&\ymk&\bmk&\omk&\omk&\bmk&\omk&\bmk&\bmk&\omk&\ymk&&\gmk&\omk&\gmk&\rmk&\rmk&\gmk&\rmk&\omk&\gmk&\rmk&\rmk&\ymk&\rmk&\ymk&&\bmk&\omk&\bmk&\ymk&\omk&\ymk&\omk&\omk&\ymk&\omk&\omk&\gmk&\omk&\rmk&25\\
                                &\bmk&\omk&\gmk&\rmk&\rmk&\bmk&\omk&\rmk&\gmk&\rmk&\rmk&\ymk&\rmk&\ymk&&\bmk&\gmk&\bmk&\gmk&\ymk&\bmk&\omk&\omk&\bmk&\bmk&\bmk&\bmk&\omk&\ymk&&\bmk&\ymk&\gmk&\rmk&\rmk&\bmk&\omk&\rmk&\gmk&\rmk&\rmk&\ymk&\rmk&\ymk&&\bmk&\omk&\bmk&\ymk&\ymk&\bmk&\omk&\omk&\ymk&\ymk&\ymk&\bmk&\omk&\omk&50\\
                                &\bmk&\ymk&\gmk&\omk&\rmk&\bmk&\omk&\omk&\bmk&\omk&\rmk&\ymk&\rmk&\ymk&&\bmk&\gmk&\bmk&\bmk&\ymk&\bmk&\omk&\omk&\bmk&\bmk&\gmk&\bmk&\omk&\ymk&&\bmk&\ymk&\gmk&\omk&\omk&\bmk&\omk&\omk&\bmk&\omk&\rmk&\ymk&\rmk&\ymk&&\bmk&\omk&\bmk&\ymk&\ymk&\bmk&\omk&\omk&\bmk&\ymk&\ymk&\bmk&\omk&\omk&75\\
                                &\bmk&\ymk&\gmk&\rmk&\rmk&\bmk&\omk&\omk&\bmk&\rmk&\rmk&\ymk&\rmk&\ymk&&\bmk&\bmk&\bmk&\bmk&\gmk&\bmk&\omk&\omk&\bmk&\bmk&\bmk&\bmk&\ymk&\ymk&&\bmk&\ymk&\gmk&\omk&\rmk&\bmk&\omk&\omk&\bmk&\rmk&\rmk&\ymk&\rmk&\ymk&&\bmk&\omk&\bmk&\bmk&\ymk&\bmk&\omk&\omk&\bmk&\gmk&\ymk&\bmk&\omk&\omk&100\\
\hline
\multirow{4}{*}{GP-PI          }&\omk&\omk&\ymk&\rmk&\rmk&\omk&\rmk&\rmk&\omk&\omk&\omk&\omk&\rmk&\bmk&&\ymk&\omk&\ymk&\rmk&\rmk&\ymk&\rmk&\rmk&\omk&\omk&\gmk&\omk&\rmk&\bmk&&\ymk&\omk&\ymk&\rmk&\rmk&\omk&\rmk&\rmk&\omk&\omk&\omk&\omk&\rmk&\bmk&&\omk&\rmk&\rmk&\rmk&\rmk&\rmk&\rmk&\rmk&\rmk&\rmk&\rmk&\rmk&\rmk&\rmk&25\\
                                &\gmk&\omk&\ymk&\rmk&\rmk&\omk&\rmk&\rmk&\bmk&\bmk&\bmk&\gmk&\gmk&\bmk&&\bmk&\omk&\ymk&\rmk&\rmk&\ymk&\rmk&\rmk&\bmk&\bmk&\bmk&\gmk&\gmk&\bmk&&\gmk&\omk&\ymk&\rmk&\rmk&\ymk&\rmk&\rmk&\bmk&\bmk&\bmk&\gmk&\gmk&\bmk&&\omk&\rmk&\rmk&\rmk&\rmk&\rmk&\rmk&\rmk&\omk&\rmk&\rmk&\rmk&\rmk&\rmk&50\\
                                &\gmk&\gmk&\gmk&\omk&\rmk&\bmk&\rmk&\rmk&\bmk&\bmk&\bmk&\bmk&\bmk&\bmk&&\bmk&\gmk&\gmk&\omk&\rmk&\bmk&\omk&\omk&\bmk&\bmk&\bmk&\bmk&\bmk&\bmk&&\gmk&\gmk&\gmk&\omk&\rmk&\bmk&\rmk&\omk&\bmk&\bmk&\bmk&\bmk&\bmk&\bmk&&\omk&\omk&\rmk&\rmk&\rmk&\rmk&\rmk&\rmk&\omk&\omk&\rmk&\rmk&\rmk&\rmk&75\\
                                &\gmk&\bmk&\bmk&\omk&\rmk&\bmk&\omk&\omk&\bmk&\bmk&\bmk&\bmk&\bmk&\bmk&&\bmk&\bmk&\bmk&\omk&\rmk&\bmk&\gmk&\omk&\bmk&\bmk&\bmk&\bmk&\bmk&\bmk&&\gmk&\bmk&\bmk&\omk&\rmk&\bmk&\omk&\omk&\bmk&\bmk&\bmk&\bmk&\bmk&\bmk&&\omk&\omk&\rmk&\rmk&\rmk&\rmk&\rmk&\rmk&\omk&\omk&\rmk&\rmk&\rmk&\rmk&100\\
\hline
\multirow{4}{*}{GP-EI          }&\gmk&\omk&\ymk&\rmk&\rmk&\rmk&\rmk&\rmk&\ymk&\omk&\omk&\omk&\rmk&\bmk&&\gmk&\omk&\gmk&\omk&\rmk&\omk&\rmk&\rmk&\gmk&\omk&\gmk&\omk&\omk&\bmk&&\gmk&\omk&\ymk&\rmk&\rmk&\rmk&\rmk&\rmk&\ymk&\omk&\omk&\omk&\omk&\bmk&&\omk&\rmk&\rmk&\rmk&\rmk&\rmk&\rmk&\rmk&\omk&\rmk&\rmk&\rmk&\rmk&\rmk&25\\
                                &\bmk&\omk&\bmk&\rmk&\rmk&\rmk&\rmk&\rmk&\bmk&\bmk&\bmk&\bmk&\bmk&\bmk&&\bmk&\gmk&\bmk&\omk&\rmk&\omk&\rmk&\rmk&\bmk&\bmk&\bmk&\bmk&\bmk&\bmk&&\bmk&\omk&\bmk&\omk&\rmk&\rmk&\rmk&\rmk&\bmk&\bmk&\bmk&\bmk&\bmk&\bmk&&\omk&\rmk&\rmk&\rmk&\rmk&\rmk&\rmk&\rmk&\omk&\rmk&\rmk&\rmk&\rmk&\rmk&50\\
                                &\bmk&\ymk&\bmk&\omk&\rmk&\omk&\rmk&\rmk&\bmk&\bmk&\bmk&\bmk&\bmk&\bmk&&\bmk&\gmk&\bmk&\omk&\rmk&\omk&\omk&\omk&\bmk&\bmk&\bmk&\bmk&\bmk&\bmk&&\bmk&\ymk&\bmk&\omk&\rmk&\omk&\rmk&\rmk&\bmk&\bmk&\bmk&\bmk&\bmk&\bmk&&\ymk&\omk&\rmk&\rmk&\rmk&\rmk&\rmk&\rmk&\omk&\omk&\rmk&\omk&\rmk&\rmk&75\\
                                &\bmk&\gmk&\bmk&\gmk&\rmk&\bmk&\omk&\rmk&\bmk&\bmk&\bmk&\bmk&\bmk&\bmk&&\bmk&\gmk&\bmk&\bmk&\rmk&\bmk&\gmk&\omk&\bmk&\bmk&\bmk&\bmk&\bmk&\bmk&&\bmk&\gmk&\bmk&\gmk&\rmk&\bmk&\omk&\omk&\bmk&\bmk&\bmk&\bmk&\bmk&\bmk&&\gmk&\omk&\rmk&\rmk&\rmk&\omk&\rmk&\rmk&\ymk&\omk&\omk&\omk&\rmk&\rmk&100\\
\hline
\multirow{4}{*}{GP-EIM         }&\gmk&\ymk&\ymk&\rmk&\rmk&\bmk&\bmk&\rmk&\bmk&\rmk&\rmk&\ymk&\rmk&\omk&&\bmk&\bmk&\bmk&\ymk&\gmk&\bmk&\bmk&\ymk&\bmk&\bmk&\bmk&\bmk&\rmk&\omk&&\bmk&\ymk&\ymk&\rmk&\rmk&\bmk&\bmk&\rmk&\bmk&\rmk&\rmk&\ymk&\rmk&\omk&&\bmk&\omk&\ymk&\omk&\omk&\omk&\omk&\rmk&\omk&\omk&\omk&\omk&\omk&\rmk&25\\
                                &\bmk&\bmk&\gmk&\omk&\rmk&\bmk&\bmk&\rmk&\bmk&\rmk&\rmk&\ymk&\rmk&\ymk&&\bmk&\bmk&\bmk&\bmk&\gmk&\bmk&\bmk&\ymk&\bmk&\bmk&\bmk&\bmk&\rmk&\ymk&&\bmk&\bmk&\gmk&\omk&\rmk&\bmk&\bmk&\rmk&\bmk&\rmk&\rmk&\ymk&\rmk&\ymk&&\bmk&\omk&\ymk&\omk&\omk&\omk&\omk&\rmk&\ymk&\omk&\omk&\omk&\omk&\rmk&50\\
                                &\bmk&\gmk&\ymk&\ymk&\omk&\bmk&\bmk&\rmk&\bmk&\rmk&\rmk&\ymk&\rmk&\ymk&&\bmk&\bmk&\bmk&\bmk&\bmk&\bmk&\bmk&\ymk&\bmk&\bmk&\bmk&\bmk&\rmk&\ymk&&\bmk&\gmk&\ymk&\ymk&\omk&\bmk&\bmk&\rmk&\bmk&\rmk&\rmk&\gmk&\rmk&\ymk&&\bmk&\omk&\omk&\omk&\omk&\omk&\omk&\rmk&\ymk&\ymk&\omk&\omk&\omk&\rmk&75\\
                                &\bmk&\gmk&\gmk&\ymk&\gmk&\bmk&\bmk&\rmk&\bmk&\rmk&\rmk&\ymk&\rmk&\ymk&&\bmk&\bmk&\bmk&\bmk&\bmk&\bmk&\bmk&\ymk&\bmk&\bmk&\bmk&\bmk&\rmk&\ymk&&\bmk&\gmk&\gmk&\ymk&\gmk&\bmk&\bmk&\rmk&\bmk&\rmk&\rmk&\ymk&\rmk&\ymk&&\bmk&\omk&\omk&\omk&\omk&\omk&\omk&\rmk&\ymk&\ymk&\ymk&\omk&\omk&\rmk&100\\
\hline
\multirow{4}{*}{GP-MVES        }&\omk&\omk&\rmk&\rmk&\rmk&\rmk&\rmk&\rmk&\omk&\rmk&\rmk&\rmk&\rmk&\omk&&\ymk&\ymk&\gmk&\rmk&\rmk&\rmk&\rmk&\rmk&\omk&\rmk&\rmk&\rmk&\rmk&\omk&&\omk&\omk&\rmk&\rmk&\rmk&\rmk&\rmk&\rmk&\omk&\rmk&\rmk&\rmk&\rmk&\omk&&\rmk&\rmk&\rmk&\rmk&\rmk&\rmk&\rmk&\rmk&\omk&\rmk&\rmk&\rmk&\rmk&\rmk&25\\
                                &\bmk&\omk&\rmk&\rmk&\rmk&\rmk&\rmk&\rmk&\omk&\rmk&\rmk&\rmk&\rmk&\omk&&\bmk&\ymk&\gmk&\omk&\rmk&\omk&\rmk&\rmk&\omk&\omk&\rmk&\omk&\rmk&\gmk&&\bmk&\omk&\rmk&\rmk&\rmk&\rmk&\rmk&\rmk&\omk&\rmk&\rmk&\rmk&\rmk&\omk&&\omk&\rmk&\rmk&\rmk&\rmk&\rmk&\rmk&\rmk&\omk&\rmk&\rmk&\rmk&\rmk&\rmk&50\\
                                &\bmk&\omk&\gmk&\rmk&\rmk&\rmk&\rmk&\rmk&\omk&\bmk&\rmk&\omk&\rmk&\bmk&&\bmk&\ymk&\bmk&\omk&\rmk&\omk&\rmk&\rmk&\gmk&\bmk&\omk&\omk&\omk&\bmk&&\bmk&\omk&\gmk&\rmk&\rmk&\rmk&\rmk&\rmk&\omk&\bmk&\omk&\omk&\omk&\bmk&&\omk&\omk&\omk&\rmk&\rmk&\rmk&\rmk&\rmk&\omk&\omk&\omk&\rmk&\rmk&\rmk&75\\
                                &\bmk&\ymk&\gmk&\omk&\rmk&\omk&\rmk&\rmk&\bmk&\bmk&\gmk&\omk&\omk&\bmk&&\bmk&\gmk&\bmk&\omk&\rmk&\omk&\rmk&\rmk&\bmk&\bmk&\bmk&\omk&\omk&\bmk&&\bmk&\ymk&\gmk&\omk&\rmk&\omk&\rmk&\rmk&\bmk&\bmk&\bmk&\omk&\omk&\bmk&&\omk&\omk&\omk&\rmk&\rmk&\rmk&\rmk&\rmk&\ymk&\omk&\omk&\omk&\rmk&\rmk&100\\
\hline
\multirow{4}{*}{GP-MSES        }&\rmk&\ymk&\ymk&\rmk&\rmk&\omk&\nmk&\nmk&\rmk&\rmk&\rmk&\omk&\nmk&\omk&&\omk&\ymk&\gmk&\gmk&\rmk&\ymk&\nmk&\nmk&\ymk&\gmk&\gmk&\ymk&\nmk&\omk&&\rmk&\ymk&\ymk&\rmk&\rmk&\omk&\nmk&\nmk&\rmk&\rmk&\rmk&\omk&\nmk&\omk&&\omk&\omk&\ymk&\omk&\rmk&\rmk&\nmk&\nmk&\rmk&\rmk&\rmk&\omk&\nmk&\rmk&25\\
                                &\rmk&\ymk&\ymk&\rmk&\rmk&\omk&\nmk&\nmk&\rmk&\rmk&\rmk&\omk&\nmk&\omk&&\omk&\gmk&\gmk&\bmk&\rmk&\ymk&\nmk&\nmk&\ymk&\gmk&\gmk&\gmk&\nmk&\omk&&\omk&\gmk&\ymk&\omk&\rmk&\omk&\nmk&\nmk&\rmk&\rmk&\rmk&\omk&\nmk&\omk&&\omk&\omk&\ymk&\omk&\rmk&\rmk&\nmk&\nmk&\rmk&\rmk&\rmk&\omk&\nmk&\rmk&50\\
                                &\omk&\ymk&\ymk&\rmk&\rmk&\omk&\nmk&\nmk&\rmk&\rmk&\rmk&\omk&\nmk&\omk&&\omk&\bmk&\gmk&\bmk&\rmk&\ymk&\nmk&\nmk&\ymk&\gmk&\gmk&\gmk&\nmk&\omk&&\omk&\gmk&\ymk&\omk&\rmk&\omk&\nmk&\nmk&\rmk&\rmk&\rmk&\omk&\nmk&\omk&&\omk&\omk&\ymk&\omk&\rmk&\rmk&\nmk&\nmk&\rmk&\rmk&\rmk&\omk&\nmk&\rmk&75\\
                                &\omk&\gmk&\ymk&\rmk&\rmk&\omk&\nmk&\nmk&\rmk&\rmk&\rmk&\omk&\nmk&\omk&&\omk&\bmk&\gmk&\bmk&\rmk&\ymk&\nmk&\nmk&\ymk&\gmk&\gmk&\gmk&\nmk&\ymk&&\omk&\gmk&\ymk&\omk&\rmk&\omk&\nmk&\nmk&\rmk&\rmk&\rmk&\omk&\nmk&\omk&&\omk&\omk&\ymk&\omk&\rmk&\rmk&\nmk&\nmk&\rmk&\rmk&\rmk&\omk&\nmk&\rmk&100\\
\hline
\multirow{4}{*}{GP-LCB         }&\ymk&\nmk&\rmk&\omk&\omk&\omk&\rmk&\rmk&\rmk&\bmk&\bmk&\bmk&\rmk&\nmk&&\ymk&\nmk&\gmk&\gmk&\gmk&\omk&\omk&\rmk&\ymk&\bmk&\bmk&\bmk&\omk&\nmk&&\ymk&\nmk&\rmk&\omk&\omk&\omk&\rmk&\rmk&\rmk&\bmk&\bmk&\bmk&\rmk&\nmk&&\rmk&\nmk&\rmk&\rmk&\omk&\rmk&\rmk&\rmk&\rmk&\rmk&\rmk&\rmk&\rmk&\nmk&25\\
                                &\ymk&\nmk&\rmk&\omk&\omk&\omk&\rmk&\rmk&\rmk&\bmk&\bmk&\bmk&\rmk&\nmk&&\ymk&\nmk&\ymk&\ymk&\gmk&\ymk&\omk&\rmk&\rmk&\bmk&\bmk&\bmk&\omk&\nmk&&\ymk&\nmk&\rmk&\omk&\omk&\omk&\rmk&\rmk&\rmk&\bmk&\bmk&\bmk&\rmk&\nmk&&\omk&\nmk&\rmk&\rmk&\omk&\rmk&\rmk&\rmk&\rmk&\rmk&\rmk&\rmk&\rmk&\nmk&50\\
                                &\ymk&\nmk&\omk&\rmk&\omk&\omk&\rmk&\rmk&\rmk&\bmk&\bmk&\bmk&\rmk&\nmk&&\gmk&\nmk&\bmk&\ymk&\gmk&\ymk&\omk&\rmk&\rmk&\bmk&\bmk&\bmk&\omk&\nmk&&\ymk&\nmk&\omk&\omk&\omk&\omk&\rmk&\rmk&\rmk&\bmk&\bmk&\bmk&\rmk&\nmk&&\omk&\nmk&\rmk&\rmk&\rmk&\rmk&\rmk&\rmk&\rmk&\rmk&\rmk&\rmk&\rmk&\nmk&75\\
                                &\ymk&\nmk&\ymk&\omk&\omk&\omk&\rmk&\rmk&\rmk&\bmk&\bmk&\bmk&\rmk&\nmk&&\gmk&\nmk&\bmk&\ymk&\gmk&\omk&\omk&\rmk&\rmk&\bmk&\bmk&\bmk&\omk&\nmk&&\ymk&\nmk&\ymk&\omk&\omk&\omk&\rmk&\rmk&\rmk&\bmk&\bmk&\bmk&\rmk&\nmk&&\omk&\nmk&\rmk&\rmk&\rmk&\rmk&\rmk&\rmk&\rmk&\rmk&\rmk&\rmk&\rmk&\nmk&100\\
\hline
\multirow{4}{*}{NGP            }&\gmk&\ymk&\omk&\rmk&\rmk&\omk&\omk&\rmk&\gmk&\rmk&\rmk&\gmk&\rmk&\omk&&\bmk&\bmk&\bmk&\bmk&\omk&\omk&\omk&\omk&\bmk&\bmk&\bmk&\bmk&\rmk&\ymk&&\gmk&\gmk&\omk&\omk&\rmk&\omk&\omk&\rmk&\gmk&\rmk&\rmk&\gmk&\rmk&\omk&&\gmk&\omk&\bmk&\ymk&\omk&\omk&\omk&\omk&\ymk&\omk&\omk&\ymk&\omk&\omk&25\\
                                &\gmk&\ymk&\ymk&\omk&\rmk&\ymk&\omk&\rmk&\gmk&\rmk&\rmk&\gmk&\rmk&\ymk&&\bmk&\bmk&\bmk&\bmk&\ymk&\ymk&\omk&\omk&\bmk&\bmk&\bmk&\bmk&\rmk&\gmk&&\bmk&\ymk&\ymk&\omk&\rmk&\ymk&\omk&\rmk&\bmk&\rmk&\rmk&\gmk&\rmk&\ymk&&\bmk&\omk&\bmk&\gmk&\ymk&\gmk&\omk&\omk&\ymk&\omk&\ymk&\gmk&\omk&\omk&50\\
                                &\gmk&\bmk&\ymk&\rmk&\rmk&\bmk&\rmk&\rmk&\bmk&\rmk&\rmk&\gmk&\rmk&\ymk&&\bmk&\bmk&\bmk&\bmk&\ymk&\bmk&\omk&\omk&\bmk&\bmk&\bmk&\bmk&\rmk&\gmk&&\gmk&\bmk&\ymk&\omk&\rmk&\bmk&\rmk&\rmk&\bmk&\rmk&\rmk&\gmk&\rmk&\ymk&&\gmk&\omk&\bmk&\gmk&\ymk&\bmk&\omk&\omk&\ymk&\ymk&\ymk&\gmk&\omk&\omk&75\\
                                &\bmk&\bmk&\gmk&\rmk&\rmk&\ymk&\rmk&\omk&\bmk&\rmk&\rmk&\ymk&\rmk&\ymk&&\bmk&\bmk&\bmk&\bmk&\ymk&\ymk&\ymk&\omk&\bmk&\bmk&\bmk&\bmk&\rmk&\ymk&&\bmk&\bmk&\gmk&\omk&\rmk&\ymk&\rmk&\omk&\bmk&\rmk&\rmk&\gmk&\rmk&\ymk&&\bmk&\omk&\bmk&\gmk&\ymk&\bmk&\omk&\omk&\ymk&\ymk&\ymk&\bmk&\omk&\omk&100\\
\hline
\multirow{4}{*}{KG             }&\rmk&\omk&\rmk&\rmk&\rmk&\omk&\rmk&\rmk&\rmk&\rmk&\rmk&\omk&\rmk&\rmk&&\rmk&\omk&\omk&\ymk&\rmk&\ymk&\rmk&\rmk&\rmk&\gmk&\gmk&\omk&\rmk&\omk&&\rmk&\omk&\rmk&\rmk&\rmk&\omk&\rmk&\rmk&\rmk&\rmk&\rmk&\omk&\rmk&\rmk&&\rmk&\rmk&\rmk&\omk&\rmk&\rmk&\rmk&\rmk&\rmk&\rmk&\omk&\rmk&\omk&\rmk&25\\
                                &\rmk&\omk&\rmk&\rmk&\rmk&\omk&\gmk&\rmk&\rmk&\rmk&\gmk&\omk&\rmk&\rmk&&\rmk&\omk&\omk&\omk&\rmk&\omk&\omk&\rmk&\rmk&\gmk&\gmk&\omk&\rmk&\rmk&&\rmk&\omk&\rmk&\rmk&\rmk&\omk&\omk&\rmk&\rmk&\rmk&\gmk&\omk&\rmk&\rmk&&\rmk&\rmk&\rmk&\omk&\rmk&\omk&\rmk&\rmk&\rmk&\rmk&\omk&\omk&\omk&\rmk&50\\
                                &\rmk&\omk&\rmk&\rmk&\rmk&\ymk&\omk&\rmk&\rmk&\gmk&\gmk&\omk&\rmk&\rmk&&\rmk&\omk&\omk&\omk&\omk&\omk&\omk&\rmk&\omk&\gmk&\gmk&\omk&\ymk&\rmk&&\rmk&\omk&\rmk&\rmk&\rmk&\omk&\omk&\rmk&\rmk&\gmk&\gmk&\omk&\rmk&\rmk&&\rmk&\rmk&\rmk&\omk&\rmk&\omk&\rmk&\rmk&\rmk&\rmk&\omk&\omk&\omk&\rmk&75\\
                                &\rmk&\omk&\rmk&\rmk&\rmk&\omk&\ymk&\rmk&\rmk&\omk&\rmk&\ymk&\rmk&\rmk&&\rmk&\omk&\omk&\omk&\omk&\omk&\omk&\rmk&\omk&\gmk&\gmk&\omk&\rmk&\rmk&&\rmk&\omk&\rmk&\rmk&\rmk&\omk&\omk&\rmk&\rmk&\omk&\rmk&\omk&\rmk&\rmk&&\rmk&\rmk&\omk&\omk&\omk&\omk&\rmk&\rmk&\rmk&\rmk&\omk&\gmk&\omk&\rmk&100\\
\hline
\multirow{4}{*}{KG-GLRD        }&\nmk&\nmk&\nmk&\nmk&\nmk&\nmk&\nmk&\nmk&\nmk&\nmk&\nmk&\nmk&\nmk&\nmk&&\nmk&\nmk&\nmk&\nmk&\nmk&\nmk&\nmk&\nmk&\nmk&\nmk&\nmk&\nmk&\nmk&\nmk&&\nmk&\nmk&\nmk&\nmk&\nmk&\nmk&\nmk&\nmk&\nmk&\nmk&\nmk&\nmk&\nmk&\nmk&&\gmk&\nmk&\gmk&\omk&\rmk&\omk&\rmk&\rmk&\ymk&\ymk&\ymk&\ymk&\omk&\nmk&25\\
                                &\nmk&\nmk&\nmk&\nmk&\nmk&\nmk&\nmk&\nmk&\nmk&\nmk&\nmk&\nmk&\nmk&\nmk&&\nmk&\nmk&\nmk&\nmk&\nmk&\nmk&\nmk&\nmk&\nmk&\nmk&\nmk&\nmk&\nmk&\nmk&&\nmk&\nmk&\nmk&\nmk&\nmk&\nmk&\nmk&\nmk&\nmk&\nmk&\nmk&\nmk&\nmk&\nmk&&\bmk&\nmk&\ymk&\omk&\rmk&\ymk&\omk&\rmk&\ymk&\ymk&\ymk&\ymk&\omk&\nmk&50\\
                                &\nmk&\nmk&\nmk&\nmk&\nmk&\nmk&\nmk&\nmk&\nmk&\nmk&\nmk&\nmk&\nmk&\nmk&&\nmk&\nmk&\nmk&\nmk&\nmk&\nmk&\nmk&\nmk&\nmk&\nmk&\nmk&\nmk&\nmk&\nmk&&\nmk&\nmk&\nmk&\nmk&\nmk&\nmk&\nmk&\nmk&\nmk&\nmk&\nmk&\nmk&\nmk&\nmk&&\bmk&\nmk&\ymk&\omk&\rmk&\ymk&\omk&\rmk&\ymk&\ymk&\ymk&\ymk&\omk&\nmk&75\\
                                &\nmk&\nmk&\nmk&\nmk&\nmk&\nmk&\nmk&\nmk&\nmk&\nmk&\nmk&\nmk&\nmk&\nmk&&\nmk&\nmk&\nmk&\nmk&\nmk&\nmk&\nmk&\nmk&\nmk&\nmk&\nmk&\nmk&\nmk&\nmk&&\nmk&\nmk&\nmk&\nmk&\nmk&\nmk&\nmk&\nmk&\nmk&\nmk&\nmk&\nmk&\nmk&\nmk&&\bmk&\nmk&\ymk&\omk&\rmk&\ymk&\omk&\rmk&\ymk&\ymk&\ymk&\ymk&\omk&\nmk&100\\
\hline
\multirow{4}{*}{KG-EC          }&\nmk&\nmk&\rmk&\rmk&\rmk&\gmk&\ymk&\omk&\gmk&\omk&\ymk&\gmk&\rmk&\nmk&&\nmk&\nmk&\gmk&\gmk&\omk&\bmk&\gmk&\ymk&\bmk&\bmk&\bmk&\bmk&\gmk&\nmk&&\nmk&\nmk&\rmk&\rmk&\rmk&\gmk&\ymk&\omk&\gmk&\omk&\ymk&\gmk&\rmk&\nmk&&\nmk&\nmk&\nmk&\nmk&\nmk&\nmk&\nmk&\nmk&\nmk&\nmk&\nmk&\nmk&\nmk&\nmk&25\\
                                &\nmk&\nmk&\rmk&\rmk&\rmk&\gmk&\ymk&\omk&\gmk&\omk&\rmk&\gmk&\rmk&\nmk&&\nmk&\nmk&\gmk&\gmk&\gmk&\bmk&\bmk&\ymk&\bmk&\bmk&\bmk&\bmk&\bmk&\nmk&&\nmk&\nmk&\rmk&\rmk&\rmk&\gmk&\gmk&\ymk&\gmk&\ymk&\rmk&\gmk&\omk&\nmk&&\nmk&\nmk&\nmk&\nmk&\nmk&\nmk&\nmk&\nmk&\nmk&\nmk&\nmk&\nmk&\nmk&\nmk&50\\
                                &\nmk&\nmk&\rmk&\rmk&\rmk&\gmk&\gmk&\ymk&\gmk&\omk&\rmk&\gmk&\rmk&\nmk&&\nmk&\nmk&\gmk&\gmk&\gmk&\bmk&\bmk&\ymk&\bmk&\bmk&\bmk&\bmk&\bmk&\nmk&&\nmk&\nmk&\rmk&\rmk&\rmk&\gmk&\gmk&\ymk&\gmk&\ymk&\rmk&\gmk&\omk&\nmk&&\nmk&\nmk&\nmk&\nmk&\nmk&\nmk&\nmk&\nmk&\nmk&\nmk&\nmk&\nmk&\nmk&\nmk&75\\
                                &\nmk&\nmk&\rmk&\rmk&\rmk&\gmk&\gmk&\ymk&\bmk&\omk&\rmk&\gmk&\rmk&\nmk&&\nmk&\nmk&\bmk&\gmk&\gmk&\bmk&\bmk&\ymk&\bmk&\bmk&\bmk&\bmk&\bmk&\nmk&&\nmk&\nmk&\omk&\rmk&\rmk&\bmk&\gmk&\ymk&\bmk&\ymk&\rmk&\gmk&\omk&\nmk&&\nmk&\nmk&\nmk&\nmk&\nmk&\nmk&\nmk&\nmk&\nmk&\nmk&\nmk&\nmk&\nmk&\nmk&100\\
\hline
\multirow{4}{*}{KG-AEI         }&\bmk&\nmk&\omk&\nmk&\nmk&\nmk&\nmk&\nmk&\nmk&\nmk&\nmk&\nmk&\nmk&\nmk&&\bmk&\nmk&\bmk&\nmk&\nmk&\nmk&\nmk&\nmk&\nmk&\nmk&\nmk&\nmk&\nmk&\nmk&&\bmk&\nmk&\gmk&\nmk&\nmk&\nmk&\nmk&\nmk&\nmk&\nmk&\nmk&\nmk&\nmk&\nmk&&\nmk&\nmk&\nmk&\nmk&\nmk&\nmk&\nmk&\nmk&\nmk&\nmk&\nmk&\nmk&\nmk&\nmk&25\\
                                &\bmk&\nmk&\omk&\nmk&\nmk&\nmk&\nmk&\nmk&\nmk&\nmk&\nmk&\nmk&\nmk&\nmk&&\bmk&\nmk&\bmk&\nmk&\nmk&\nmk&\nmk&\nmk&\nmk&\nmk&\nmk&\nmk&\nmk&\nmk&&\bmk&\nmk&\bmk&\nmk&\nmk&\nmk&\nmk&\nmk&\nmk&\nmk&\nmk&\nmk&\nmk&\nmk&&\nmk&\nmk&\nmk&\nmk&\nmk&\nmk&\nmk&\nmk&\nmk&\nmk&\nmk&\nmk&\nmk&\nmk&50\\
                                &\bmk&\nmk&\omk&\nmk&\nmk&\nmk&\nmk&\nmk&\nmk&\nmk&\nmk&\nmk&\nmk&\nmk&&\bmk&\nmk&\bmk&\nmk&\nmk&\nmk&\nmk&\nmk&\nmk&\nmk&\nmk&\nmk&\nmk&\nmk&&\bmk&\nmk&\bmk&\nmk&\nmk&\nmk&\nmk&\nmk&\nmk&\nmk&\nmk&\nmk&\nmk&\nmk&&\nmk&\nmk&\nmk&\nmk&\nmk&\nmk&\nmk&\nmk&\nmk&\nmk&\nmk&\nmk&\nmk&\nmk&75\\
                                &\bmk&\nmk&\omk&\nmk&\nmk&\nmk&\nmk&\nmk&\nmk&\nmk&\nmk&\nmk&\nmk&\nmk&&\bmk&\nmk&\bmk&\nmk&\nmk&\nmk&\nmk&\nmk&\nmk&\nmk&\nmk&\nmk&\nmk&\nmk&&\bmk&\nmk&\gmk&\nmk&\nmk&\nmk&\nmk&\nmk&\nmk&\nmk&\nmk&\nmk&\nmk&\nmk&&\nmk&\nmk&\nmk&\nmk&\nmk&\nmk&\nmk&\nmk&\nmk&\nmk&\nmk&\nmk&\nmk&\nmk&100\\
\hline
\multirow{4}{*}{RCK-EI         }&\gmk&\nmk&\omk&\rmk&\nmk&\rmk&\rmk&\nmk&\gmk&\bmk&\bmk&\rmk&\ymk&\nmk&&\bmk&\nmk&\bmk&\gmk&\nmk&\rmk&\rmk&\nmk&\bmk&\bmk&\bmk&\gmk&\bmk&\nmk&&\gmk&\nmk&\omk&\rmk&\nmk&\rmk&\rmk&\nmk&\gmk&\bmk&\bmk&\omk&\ymk&\nmk&&\nmk&\nmk&\nmk&\nmk&\nmk&\nmk&\nmk&\nmk&\nmk&\nmk&\nmk&\nmk&\nmk&\nmk&25\\
                                &\gmk&\nmk&\ymk&\rmk&\nmk&\rmk&\rmk&\nmk&\gmk&\bmk&\bmk&\omk&\ymk&\nmk&&\bmk&\nmk&\bmk&\gmk&\nmk&\rmk&\rmk&\nmk&\bmk&\bmk&\bmk&\bmk&\bmk&\nmk&&\gmk&\nmk&\ymk&\rmk&\nmk&\rmk&\rmk&\nmk&\gmk&\bmk&\bmk&\omk&\ymk&\nmk&&\nmk&\nmk&\nmk&\nmk&\nmk&\nmk&\nmk&\nmk&\nmk&\nmk&\nmk&\nmk&\nmk&\nmk&50\\
                                &\gmk&\nmk&\ymk&\rmk&\nmk&\rmk&\rmk&\nmk&\gmk&\bmk&\bmk&\omk&\ymk&\nmk&&\bmk&\nmk&\bmk&\bmk&\nmk&\rmk&\rmk&\nmk&\bmk&\bmk&\bmk&\bmk&\bmk&\nmk&&\gmk&\nmk&\ymk&\rmk&\nmk&\rmk&\rmk&\nmk&\bmk&\bmk&\bmk&\omk&\ymk&\nmk&&\nmk&\nmk&\nmk&\nmk&\nmk&\nmk&\nmk&\nmk&\nmk&\nmk&\nmk&\nmk&\nmk&\nmk&75\\
                                &\gmk&\nmk&\ymk&\rmk&\nmk&\rmk&\rmk&\nmk&\bmk&\bmk&\bmk&\rmk&\ymk&\nmk&&\bmk&\nmk&\bmk&\bmk&\nmk&\rmk&\rmk&\nmk&\bmk&\bmk&\bmk&\bmk&\bmk&\nmk&&\gmk&\nmk&\ymk&\omk&\nmk&\rmk&\rmk&\nmk&\bmk&\bmk&\bmk&\omk&\ymk&\nmk&&\nmk&\nmk&\nmk&\nmk&\nmk&\nmk&\nmk&\nmk&\nmk&\nmk&\nmk&\nmk&\nmk&\nmk&100\\
\hline
\multirow{4}{*}{SRBF-LCB       }&\omk&\nmk&\rmk&\omk&\omk&\rmk&\rmk&\rmk&\rmk&\bmk&\bmk&\gmk&\rmk&\omk&&\ymk&\nmk&\bmk&\gmk&\gmk&\omk&\omk&\omk&\ymk&\bmk&\bmk&\bmk&\omk&\omk&&\ymk&\nmk&\omk&\omk&\omk&\rmk&\rmk&\rmk&\rmk&\bmk&\bmk&\bmk&\rmk&\omk&&\omk&\nmk&\rmk&\rmk&\rmk&\rmk&\rmk&\rmk&\rmk&\rmk&\rmk&\rmk&\rmk&\nmk&25\\
                                &\ymk&\nmk&\rmk&\omk&\omk&\rmk&\rmk&\rmk&\rmk&\bmk&\bmk&\bmk&\rmk&\omk&&\ymk&\nmk&\bmk&\gmk&\gmk&\omk&\omk&\omk&\ymk&\bmk&\bmk&\bmk&\omk&\omk&&\ymk&\nmk&\omk&\omk&\omk&\rmk&\rmk&\rmk&\rmk&\bmk&\bmk&\bmk&\rmk&\omk&&\omk&\nmk&\rmk&\rmk&\rmk&\rmk&\rmk&\rmk&\omk&\rmk&\omk&\rmk&\rmk&\nmk&50\\
                                &\ymk&\nmk&\omk&\omk&\omk&\rmk&\rmk&\rmk&\rmk&\bmk&\bmk&\bmk&\rmk&\omk&&\gmk&\nmk&\bmk&\bmk&\gmk&\omk&\omk&\omk&\ymk&\bmk&\bmk&\bmk&\omk&\omk&&\ymk&\nmk&\omk&\omk&\omk&\rmk&\rmk&\rmk&\rmk&\bmk&\bmk&\bmk&\rmk&\omk&&\omk&\nmk&\rmk&\rmk&\rmk&\rmk&\rmk&\rmk&\omk&\rmk&\omk&\rmk&\rmk&\nmk&75\\
                                &\ymk&\nmk&\omk&\omk&\omk&\rmk&\rmk&\rmk&\rmk&\bmk&\bmk&\gmk&\rmk&\omk&&\gmk&\nmk&\bmk&\bmk&\gmk&\omk&\omk&\omk&\ymk&\bmk&\bmk&\bmk&\omk&\omk&&\ymk&\nmk&\omk&\omk&\omk&\rmk&\rmk&\rmk&\rmk&\bmk&\bmk&\gmk&\rmk&\omk&&\omk&\nmk&\rmk&\rmk&\rmk&\rmk&\rmk&\rmk&\ymk&\rmk&\omk&\rmk&\rmk&\nmk&100\\
\hline
\multirow{4}{*}{PCE            }&\nmk&\nmk&\nmk&\nmk&\nmk&\nmk&\nmk&\nmk&\nmk&\nmk&\nmk&\nmk&\nmk&\nmk&&\nmk&\nmk&\nmk&\nmk&\nmk&\nmk&\nmk&\nmk&\nmk&\nmk&\nmk&\nmk&\nmk&\nmk&&\nmk&\nmk&\nmk&\nmk&\nmk&\nmk&\nmk&\nmk&\nmk&\nmk&\nmk&\nmk&\nmk&\nmk&&\bmk&\nmk&\bmk&\omk&\omk&\omk&\rmk&\rmk&\omk&\omk&\omk&\omk&\rmk&\nmk&25\\
                                &\nmk&\nmk&\nmk&\nmk&\nmk&\nmk&\nmk&\nmk&\nmk&\nmk&\nmk&\nmk&\nmk&\nmk&&\nmk&\nmk&\nmk&\nmk&\nmk&\nmk&\nmk&\nmk&\nmk&\nmk&\nmk&\nmk&\nmk&\nmk&&\nmk&\nmk&\nmk&\nmk&\nmk&\nmk&\nmk&\nmk&\nmk&\nmk&\nmk&\nmk&\nmk&\nmk&&\bmk&\nmk&\nmk&\bmk&\omk&\omk&\rmk&\rmk&\omk&\omk&\omk&\gmk&\rmk&\nmk&50\\
                                &\nmk&\nmk&\nmk&\nmk&\nmk&\nmk&\nmk&\nmk&\nmk&\nmk&\nmk&\nmk&\nmk&\nmk&&\nmk&\nmk&\nmk&\nmk&\nmk&\nmk&\nmk&\nmk&\nmk&\nmk&\nmk&\nmk&\nmk&\nmk&&\nmk&\nmk&\nmk&\nmk&\nmk&\nmk&\nmk&\nmk&\nmk&\nmk&\nmk&\nmk&\nmk&\nmk&&\bmk&\nmk&\nmk&\nmk&\bmk&\ymk&\omk&\rmk&\omk&\omk&\omk&\gmk&\rmk&\nmk&75\\
                                &\nmk&\nmk&\nmk&\nmk&\nmk&\nmk&\nmk&\nmk&\nmk&\nmk&\nmk&\nmk&\nmk&\nmk&&\nmk&\nmk&\nmk&\nmk&\nmk&\nmk&\nmk&\nmk&\nmk&\nmk&\nmk&\nmk&\nmk&\nmk&&\nmk&\nmk&\nmk&\nmk&\nmk&\nmk&\nmk&\nmk&\nmk&\nmk&\nmk&\nmk&\nmk&\nmk&&\nmk&\nmk&\nmk&\nmk&\nmk&\ymk&\nmk&\rmk&\ymk&\ymk&\omk&\bmk&\omk&\nmk&100\\
\hline
\multirow{4}{*}{SPCE           }&\nmk&\nmk&\nmk&\nmk&\nmk&\nmk&\nmk&\nmk&\nmk&\nmk&\nmk&\nmk&\nmk&\nmk&&\nmk&\nmk&\nmk&\nmk&\nmk&\nmk&\nmk&\nmk&\nmk&\nmk&\nmk&\nmk&\nmk&\nmk&&\nmk&\nmk&\nmk&\nmk&\nmk&\nmk&\nmk&\nmk&\nmk&\nmk&\nmk&\nmk&\nmk&\nmk&&\bmk&\nmk&\bmk&\nmk&\nmk&\omk&\nmk&\nmk&\omk&\omk&\nmk&\ymk&\nmk&\nmk&25\\
                                &\nmk&\nmk&\nmk&\nmk&\nmk&\nmk&\nmk&\nmk&\nmk&\nmk&\nmk&\nmk&\nmk&\nmk&&\nmk&\nmk&\nmk&\nmk&\nmk&\nmk&\nmk&\nmk&\nmk&\nmk&\nmk&\nmk&\nmk&\nmk&&\nmk&\nmk&\nmk&\nmk&\nmk&\nmk&\nmk&\nmk&\nmk&\nmk&\nmk&\nmk&\nmk&\nmk&&\bmk&\nmk&\nmk&\nmk&\nmk&\gmk&\nmk&\nmk&\omk&\omk&\nmk&\gmk&\nmk&\nmk&50\\
                                &\nmk&\nmk&\nmk&\nmk&\nmk&\nmk&\nmk&\nmk&\nmk&\nmk&\nmk&\nmk&\nmk&\nmk&&\nmk&\nmk&\nmk&\nmk&\nmk&\nmk&\nmk&\nmk&\nmk&\nmk&\nmk&\nmk&\nmk&\nmk&&\nmk&\nmk&\nmk&\nmk&\nmk&\nmk&\nmk&\nmk&\nmk&\nmk&\nmk&\nmk&\nmk&\nmk&&\bmk&\nmk&\nmk&\nmk&\nmk&\gmk&\nmk&\nmk&\omk&\omk&\nmk&\gmk&\nmk&\nmk&75\\
                                &\nmk&\nmk&\nmk&\nmk&\nmk&\nmk&\nmk&\nmk&\nmk&\nmk&\nmk&\nmk&\nmk&\nmk&&\nmk&\nmk&\nmk&\nmk&\nmk&\nmk&\nmk&\nmk&\nmk&\nmk&\nmk&\nmk&\nmk&\nmk&&\nmk&\nmk&\nmk&\nmk&\nmk&\nmk&\nmk&\nmk&\nmk&\nmk&\nmk&\nmk&\nmk&\nmk&&\nmk&\nmk&\nmk&\nmk&\nmk&\bmk&\nmk&\nmk&\omk&\omk&\nmk&\gmk&\nmk&\nmk&100\\
\hline
\multirow{4}{*}{SPCE\_KG-GLRD  }&\nmk&\nmk&\nmk&\nmk&\nmk&\nmk&\nmk&\nmk&\nmk&\nmk&\nmk&\nmk&\nmk&\nmk&&\nmk&\nmk&\nmk&\nmk&\nmk&\nmk&\nmk&\nmk&\nmk&\nmk&\nmk&\nmk&\nmk&\nmk&&\nmk&\nmk&\nmk&\nmk&\nmk&\nmk&\nmk&\nmk&\nmk&\nmk&\nmk&\nmk&\nmk&\nmk&&\gmk&\nmk&\bmk&\nmk&\nmk&\omk&\nmk&\nmk&\ymk&\ymk&\ymk&\ymk&\nmk&\nmk&25\\
                                &\nmk&\nmk&\nmk&\nmk&\nmk&\nmk&\nmk&\nmk&\nmk&\nmk&\nmk&\nmk&\nmk&\nmk&&\nmk&\nmk&\nmk&\nmk&\nmk&\nmk&\nmk&\nmk&\nmk&\nmk&\nmk&\nmk&\nmk&\nmk&&\nmk&\nmk&\nmk&\nmk&\nmk&\nmk&\nmk&\nmk&\nmk&\nmk&\nmk&\nmk&\nmk&\nmk&&\bmk&\nmk&\gmk&\nmk&\nmk&\ymk&\nmk&\nmk&\ymk&\ymk&\ymk&\ymk&\nmk&\nmk&50\\
                                &\nmk&\nmk&\nmk&\nmk&\nmk&\nmk&\nmk&\nmk&\nmk&\nmk&\nmk&\nmk&\nmk&\nmk&&\nmk&\nmk&\nmk&\nmk&\nmk&\nmk&\nmk&\nmk&\nmk&\nmk&\nmk&\nmk&\nmk&\nmk&&\nmk&\nmk&\nmk&\nmk&\nmk&\nmk&\nmk&\nmk&\nmk&\nmk&\nmk&\nmk&\nmk&\nmk&&\bmk&\nmk&\nmk&\nmk&\nmk&\ymk&\nmk&\nmk&\ymk&\ymk&\ymk&\ymk&\nmk&\nmk&75\\
                                &\nmk&\nmk&\nmk&\nmk&\nmk&\nmk&\nmk&\nmk&\nmk&\nmk&\nmk&\nmk&\nmk&\nmk&&\nmk&\nmk&\nmk&\nmk&\nmk&\nmk&\nmk&\nmk&\nmk&\nmk&\nmk&\nmk&\nmk&\nmk&&\nmk&\nmk&\nmk&\nmk&\nmk&\nmk&\nmk&\nmk&\nmk&\nmk&\nmk&\nmk&\nmk&\nmk&&\bmk&\nmk&\nmk&\nmk&\nmk&\ymk&\nmk&\nmk&\ymk&\ymk&\ymk&\ymk&\nmk&\nmk&100\\
\midrule
MF Benchmark & \rotatebox[origin=c]{90}{1.1} & \rotatebox[origin=c]{90}{1.2} & \rotatebox[origin=c]{90}{2.1} & \rotatebox[origin=c]{90}{2.2} & \rotatebox[origin=c]{90}{2.3}  & \rotatebox[origin=c]{90}{3.1} & \rotatebox[origin=c]{90}{3.2} & \rotatebox[origin=c]{90}{3.3}  & \rotatebox[origin=c]{90}{4.1} & \rotatebox[origin=c]{90}{4.2} & \rotatebox[origin=c]{90}{4.3}  & \rotatebox[origin=c]{90}{5.1} & \rotatebox[origin=c]{90}{5.2} & \rotatebox[origin=c]{90}{6} & & \rotatebox[origin=c]{90}{1.1} & \rotatebox[origin=c]{90}{1.2} & \rotatebox[origin=c]{90}{2.1} & \rotatebox[origin=c]{90}{2.2} & \rotatebox[origin=c]{90}{2.3}  & \rotatebox[origin=c]{90}{3.1} & \rotatebox[origin=c]{90}{3.2} & \rotatebox[origin=c]{90}{3.3}  & \rotatebox[origin=c]{90}{4.1} & \rotatebox[origin=c]{90}{4.2} & \rotatebox[origin=c]{90}{4.3}  & \rotatebox[origin=c]{90}{5.1} & \rotatebox[origin=c]{90}{5.2} & \rotatebox[origin=c]{90}{6} & & \rotatebox[origin=c]{90}{1.1} & \rotatebox[origin=c]{90}{1.2} & \rotatebox[origin=c]{90}{2.1} & \rotatebox[origin=c]{90}{2.2} & \rotatebox[origin=c]{90}{2.3}  & \rotatebox[origin=c]{90}{3.1} & \rotatebox[origin=c]{90}{3.2} & \rotatebox[origin=c]{90}{3.3}  & \rotatebox[origin=c]{90}{4.1} & \rotatebox[origin=c]{90}{4.2} & \rotatebox[origin=c]{90}{4.3}  & \rotatebox[origin=c]{90}{5.1} & \rotatebox[origin=c]{90}{5.2} & \rotatebox[origin=c]{90}{6} & & \rotatebox[origin=c]{90}{1.1} & \rotatebox[origin=c]{90}{1.2} & \rotatebox[origin=c]{90}{2.1} & \rotatebox[origin=c]{90}{2.2} & \rotatebox[origin=c]{90}{2.3}  & \rotatebox[origin=c]{90}{3.1} & \rotatebox[origin=c]{90}{3.2} & \rotatebox[origin=c]{90}{3.3}  & \rotatebox[origin=c]{90}{4.1} & \rotatebox[origin=c]{90}{4.2} & \rotatebox[origin=c]{90}{4.3}  & \rotatebox[origin=c]{90}{5.1} & \rotatebox[origin=c]{90}{5.2} & \rotatebox[origin=c]{90}{6}  & \\
 \bottomrule
 \\
\multicolumn{58}{p{10.5cm}}{\rmk \, = $\mathcal{E}_y>$10\%, \omk \, = 1\%$\leq \mathcal{E}_y\leq$10\%, \ymk \, = 0.1\%$\leq \mathcal{E}_y\leq$1\%, \gmk \, = 0.01\%$\leq \mathcal{E}_y\leq$0.1\%, \bmk \, = $\mathcal{E}_y<$0.01\%, \nmk \, = not evaluated}\\
\end{tabular}
\end{table}

\section{Conclusions}\label{sec:conclusion}
This paper proposes and discusses a benchmark suite of fourteen problems of class L1, specifically designed to assess the efficiency and effectiveness of multifidelity optimization methods. The proposed benchmarks are meant to stress multifidelity optimization methods, dealing with various challenges typical of complex real-world optimization problems such as handling localized, multimodal, and discontinuous behaviors in the objective functions, as well as handling the possible presence of noise in the objective functions. A distinctive characteristic of the proposed problems is their ease of implementation. Specifically, these problems are based on analytic functions (such as Forrester, Rosenbrock, shifted-rotated Rastrigin, and Paciorek) and the solution of a spring-mass system. The set has been implemented and validated for Matlab, Fortran and Python environments and the codes are publicly available to the community.

The paper also provides recommended experimental setups to support the rigorous test and comparison of different multifidelity methods. A quite large variety of methods was used to solve the benchmark problems and the results have been reported and discussed. These include surrogate lean and surrogate rich approaches. The former include BFGS and genetic algorithms approaches. The latter includes both offline and online collection of training data. Online approaches use adaptive sampling techniques, where several acquisition functions are used for comparison. Surrogate models include numerous Gaussian processes, radial-basis functions, and polynomial chaos expansions. The performances of the various methods vary significantly for each problem and with dimensionality. In particular it is observed that BFGS, GP-EIM, NGP, and KG-EC seem to overall perform best for search and  optimization problems, whereas NGP, GP, and PCE are best suited for the construction of surrogate models that target high global accuracy. 
\section*{Acknowledgements} 
The work was conducted in collaboration with the NATO Science and Technology Organization task group AVT-331. The authors thank the AIAA Multidisciplinary Design Optimization Technical Committee (MDO TC) and the community of experts and practitioners that joined the AIAA Workshop on Multifidelity modelling for design and uncertainty quantification for the valuable interactions. 
Distribution A: Approved for public release, distribution unlimited. Case number AFRL-2024-6631. 

\bibliographystyle{abbrv}  
\bibliography{references}  

@inproceedings{hwang2019large,
  title={Large-scale multidisciplinary design optimization—review and recommendations},
  author={Hwang, John T and Jain, Arnav V and Ha, Tae H},
  booktitle={AIAA Aviation 2019 Forum},
  pages={3106},
  year={2019},
  doi = {10.2514/6.2019-3106}
}

@article{2018-JMeST-Li_etal,
  title={Kriging-based multi-fidelity optimization via information fusion with uncertainty},
  author={Li, Chengshan and Wang, Peng and Dong, Huachao},
  journal={Journal of Mechanical Science and Technology},
  volume={32},
  number={1},
  pages={245--259},
  year={2018},
  publisher={Springer},
  doi = {10.1007/s12206-017-1225-7}
}

@article{2016-JCS-Liu_etal,
  title={A multi-fidelity surrogate-model-assisted evolutionary algorithm for computationally expensive optimization problems},
  author={Liu, Bo and Koziel, Slawomir and Zhang, Qingfu},
  journal={Journal of Computational Science},
  volume={12},
  pages={28--37},
  year={2016},
  publisher={Elsevier},
  doi = {10.1016/j.jocs.2015.11.004}
}

@article{sarkar2019multifidelity,
  title={Multifidelity and multiscale Bayesian framework for high-dimensional engineering design and calibration},
  author={Sarkar, Soumalya and Mondal, Sudeepta and Joly, Michael and Lynch, Matthew E and Bopardikar, Shaunak D and Acharya, Ranadip and Perdikaris, Paris},
  journal={Journal of Mechanical Design},
  volume={141},
  number={12},
  pages={121001},
  year={2019},
  publisher={American Society of Mechanical Engineers},
  doi = {10.1115/1.4044598}
}

@article{perdikaris2017nonlinear,
  title={Nonlinear information fusion algorithms for data-efficient multi-fidelity modelling},
  author={Perdikaris, Paris and Raissi, Maziar and Damianou, Andreas and Lawrence, Neil D and Karniadakis, George Em},
  journal={Proceedings of the Royal Society A: Mathematical, Physical and Engineering Sciences},
  volume={473},
  number={2198},
  pages={20160751},
  year={2017},
  publisher={The Royal Society Publishing},
  doi = {10.1098/rspa.2016.0751}
}

@article{ng2014multifidelity,
  title={Multifidelity approaches for optimization under uncertainty},
  author={Ng, Leo WT and Willcox, Karen E},
  journal={International Journal for Numerical Methods in Engineering},
  volume={100},
  number={10},
  pages={746--772},
  year={2014},
  publisher={Wiley Online Library},
  doi = {10.1002/nme.4761}
}

@article{li2022machine,
  title={Machine learning in aerodynamic shape optimization},
  author={Li, Jichao and Du, Xiaosong and Martins, Joaquim RRA},
  journal={Progress in Aerospace Sciences},
  volume={134},
  pages={100849},
  year={2022},
  publisher={Elsevier},
  doi = {10.1016/j.paerosci.2022.100849}
}

@article{khatouri2022metamodeling,
  title={Metamodeling techniques for CPU-intensive simulation-based design optimization: a survey},
  author={Khatouri, Hanane and Benamara, Tariq and Breitkopf, Piotr and Demange, Jean},
  journal={Advanced Modeling and Simulation in Engineering Sciences},
  volume={9},
  number={1},
  pages={1},
  year={2022},
  publisher={Springer},
  doi = {10.1186/s40323-022-00214-y}
}

@article{2022-CMAME-Guo_etal,
  title={Multi-fidelity regression using artificial neural networks: efficient approximation of parameter-dependent output quantities},
  author={Guo, Mengwu and Manzoni, Andrea and Amendt, Maurice and Conti, Paolo and Hesthaven, Jan S},
  journal={Computer Methods in Applied Mechanics and Engineering},
  volume={389},
  pages={114378},
  year={2022},
  publisher={Elsevier},
  doi = {10.1016/j.cma.2021.114378}
}

@article{2019-IISE-Song_etal,
  title={Multi-fidelity sampling for efficient simulation-based decision making in manufacturing management},
  author={Song, Jie and Qiu, Yunzhe and Xu, Jie and Yang, Feng},
  journal={IISE Transactions},
  volume={51},
  number={7},
  pages={792--805},
  year={2019},
  publisher={Taylor \& Francis},
  doi = {10.1080/24725854.2019.1576951}
}

@inproceedings{de2007state,
  title={State-of-the-art and future trends in multidisciplinary design optimization},
  author={de Weck, Olivier and Agte, Jeremy and Sobieszczanski-Sobieski, Jaroslaw and Arendsen, Paul and Morris, Alan and Spieck, Martin},
  booktitle={48th Aiaa/Asme/Asce/Ahs/Asc Structures, Structural Dynamics, and Materials Conference},
  pages={1905},
  year={2007},
  doi = {10.2514/6.2007-1905}
}

@inproceedings{vanaret2017consequences,
  title={On the consequences of the "no free lunch" theorem for optimization on the choice of an appropriate mdo architecture},
  author={Vanaret, Charlie and Gallard, Fran{\c{c}}ois and Martins, Joaquim RRA},
  booktitle={18th AIAA/ISSMO Multidisciplinary Analysis and Optimization Conference},
  pages={3148},
  year={2017},
  doi = {10.2514/6.2017-3148}
}

@article{foumani2023multi,
  title={Multi-fidelity cost-aware Bayesian optimization},
  author={Foumani, Zahra Zanjani and Shishehbor, Mehdi and Yousefpour, Amin and Bostanabad, Ramin},
  journal={Computer Methods in Applied Mechanics and Engineering},
  volume={407},
  pages={115937},
  year={2023},
  publisher={Elsevier},
  doi = {10.1016/j.cma.2023.115937}
}

@article{2021-CMAME-Zhang_etal,
  title={Multi-fidelity deep neural network surrogate model for aerodynamic shape optimization},
  author={Zhang, Xinshuai and Xie, Fangfang and Ji, Tingwei and Zhu, Zaoxu and Zheng, Yao},
  journal={Computer Methods in Applied Mechanics and Engineering},
  volume={373},
  pages={113485},
  year={2021},
  publisher={Elsevier},
  doi = {10.1016/j.cma.2020.113485}
}

@article{2019-IJCFD-Serani_etal,
  title={Adaptive multi-fidelity sampling for {CFD}-based optimisation via radial basis function metamodels},
  author={Serani, A and Pellegrini, R and Wackers, J and Jeanson, C-E and Queutey, P and Visonneau, Michel and Diez, M},
  journal={International Journal of Computational Fluid Dynamics},
  volume={33},
  number={6-7},
  pages={237--255},
  year={2019},
  publisher={Taylor \& Francis},
  doi = {10.1080/10618562.2019.1683164}
}

@article{he2021multi,
  title={Multi-objective efficient global optimization of expensive simulation-based problem in presence of simulation failures},
  author={He, Youwei and Sun, Jinju and Song, Peng and Wang, Xuesong},
  journal={Engineering with Computers},
  pages={1--26},
  year={2021},
  publisher={Springer},
  doi = {10.1007/s00366-021-01351-5}
}

@article{peherstorfer2018survey,
  title={Survey of multifidelity methods in uncertainty propagation, inference, and optimization},
  author={Peherstorfer, Benjamin and Willcox, Karen and Gunzburger, Max},
  journal={Siam Review},
  volume={60},
  number={3},
  pages={550--591},
  year={2018},
  publisher={SIAM},
  doi = {10.1137/16M1082469}
}

@article{lucia2004reduced,
  title={Reduced-order modeling: new approaches for computational physics},
  author={Lucia, David J and Beran, Philip S and Silva, Walter A},
  journal={Progress in Aerospace Sciences},
  volume={40},
  number={1-2},
  pages={51--117},
  year={2004},
  publisher={Elsevier},
  doi = {10.1016/j.paerosci.2003.12.001}
}

@article{kennedy2001bayesian,
  title={Bayesian calibration of computer models},
  author={Kennedy, Marc C and O'Hagan, Anthony},
  journal={Journal of the Royal Statistical Society: Series B (Statistical Methodology)},
  volume={63},
  number={3},
  pages={425--464},
  year={2001},
  publisher={Wiley Online Library},
  doi = {10.1111/1467-9868.00294}
}

@article{park2017remarks,
  title={Remarks on multi-fidelity surrogates},
  author={Park, Chanyoung and Haftka, Raphael T and Kim, Nam H},
  journal={Structural and Multidisciplinary Optimization},
  volume={55},
  number={3},
  pages={1029--1050},
  year={2017},
  publisher={Springer},
  doi = {10.1007/s00158-016-1550-y}
}

@article{2021-KBS-Liu_etal,
  title={Multi-fidelity global optimization using a data-mining strategy for computationally intensive black-box problems},
  author={Liu, Jie and Dong, Huachao and Wang, Peng},
  journal={Knowledge-Based Systems},
  volume={227},
  pages={107212},
  year={2021},
  publisher={Elsevier},
  doi = {10.1016/j.knosys.2021.107212}
}

@article{2021-JMLR-Moss_etal,
  title={{GIBBON}: General-purpose information-based bayesian optimisation},
  author={Moss, Henry B and Leslie, David S and Gonzalez, Javier and Rayson, Paul},
  journal={Journal of Machine Learning Research},
  volume={22},
  number={235},
  pages={1--49},
  year={2021},
  doi = {}
}

@article{2020-ASOC-Li_etal,
  title={Multi-fidelity meta-optimization for nature inspired optimization algorithms},
  author={Li, Hui and Huang, Zhiguo and Liu, Xiao and Zeng, Chenbo and Zou, Peng},
  journal={Applied Soft Computing},
  volume={96},
  pages={106619},
  year={2020},
  publisher={Elsevier},
  doi = {10.1016/j.asoc.2020.106619}
}

@article{2020-SMO-Yi_etal,
  title={A multi-fidelity {RBF} surrogate-based optimization framework for computationally expensive multi-modal problems with application to capacity planning of manufacturing systems},
  author={Yi, Jin and Shen, Yichi and Shoemaker, Christine A},
  journal={Structural and Multidisciplinary Optimization},
  volume={62},
  number={4},
  pages={1787--1807},
  year={2020},
  publisher={Springer},
  doi = {10.1007/s00158-020-02575-7}
}

@article{2018-IEEE-Wang_etal,
  title={A generic test suite for evolutionary multifidelity optimization},
  author={Wang, Handing and Jin, Yaochu and Doherty, John},
  journal={IEEE Transactions on Evolutionary Computation},
  volume={22},
  number={6},
  pages={836--850},
  year={2017},
  publisher={IEEE},
  doi = {10.1109/TEVC.2017.2758360}
}

@article{2017-SMO-Durantin_etal,
  title={Multifidelity surrogate modeling based on radial basis functions},
  author={Durantin, C{\'e}dric and Rouxel, Justin and D{\'e}sid{\'e}ri, Jean-Antoine and Gli{\`e}re, Alain},
  journal={Structural and Multidisciplinary Optimization},
  volume={56},
  number={5},
  pages={1061--1075},
  year={2017},
  publisher={Springer},
  doi = {10.1007/s00158-017-1703-7}
}

@article{2017-AIAA-Cai_etal,
  title={Adaptive radial-basis-function-based multifidelity metamodeling for expensive black-box problems},
  author={Cai, Xiwen and Qiu, Haobo and Gao, Liang and Wei, Li and Shao, Xinyu},
  journal={AIAA journal},
  volume={55},
  number={7},
  pages={2424--2436},
  year={2017},
  publisher={American Institute of Aeronautics and Astronautics},
  doi = {10.2514/1.J055649}
}

@inproceedings{rumpfkeil2020-AIAA,
  title={Multi-Fidelity, Gradient-enhanced, and Locally Optimized Sparse Polynomial Chaos and Kriging Surrogate Models Applied to Benchmark Problems},
  author={Rumpfkeil, Markus P and Beran, Philip S},
  booktitle={AIAA SciTech 2020 Forum},
  year={2020-0677},
  doi = {10.2514/6.2020-0677}
}

@inproceedings{grassi2021-AIAA,
  title={Resource Aware Multifidelity Active Learning for Efficient Optimization},
  author={Grassi, Francesco and Manganini, Giorgio and Garraffa, Michele and Mainini, Laura},
  booktitle={AIAA SciTech 2021 Forum},
  pages={0894},
  year={2021},
  doi = {10.48550/arXiv.2007.04674}
}

@article{Bryson2018a,
  title={Multifidelity quasi-newton method for design optimization},
  author={Bryson, Dean E and Rumpfkeil, Markus P},
  journal={AIAA Journal},
  volume={56},
  number={10},
  pages={4074--4086},
  year={2018},
  publisher={American Institute of Aeronautics and Astronautics},
  doi = {10.2514/6.2018-0102}
}

@inproceedings{bryson2016-AIAA,
  title={Variable-fidelity surrogate modeling of lambda wing transonic aerodynamic performance},
  author={Bryson, Dean E and Rumpfkeil, Markus P},
  booktitle={54th AIAA Aerospace Sciences Meeting},
  pages={0294},
  year={2016},
  doi = {10.2514/6.2016-0294}
}

@article{bellman1957dynamic,
  title={Dynamic programming},
  author={Bellman, Richard},
  journal={Press Princeton, New Jersey},
  year={1957},
  doi = {10.2307/j.ctv1nxcw0f}
}

@article{koziel2012scaling,
  title={Scaling properties of multi-fidelity shape optimization algorithms},
  author={Koziel, Slawomir and Leifsson, Leifur},
  journal={Procedia Computer Science},
  volume={9},
  pages={832--841},
  year={2012},
  publisher={Elsevier},
  doi = {10.1016/j.procs.2012.04.089}
}

@article{lee2024comprehensive,
  title={A comprehensive multi-fidelity surrogate framework based on Gaussian process for datasets with heterogeneous responses},
  author={Lee, Juyoung and Lee, Mingyu and Lee, Bong Jae and Lee, Ikjin},
  journal={Knowledge-Based Systems},
  volume={295},
  pages={111827},
  year={2024},
  publisher={Elsevier},
  doi = {10.1016/j.knosys.2024.111827}
}

@article{pellegrini2022multi,
  title={A multi-fidelity active learning method for global design optimization problems with noisy evaluations},
  author={Pellegrini, Riccardo and Wackers, Jeroen and Broglia, Riccardo and Serani, Andrea and Visonneau, Michel and Diez, Matteo},
  journal={Engineering with Computers},
  pages={1--24},
  year={2022},
  publisher={Springer},
  doi = {10.1007/s00366-022-01728-0}
}

@article{fernandez2019-AIAA,
  title={Issues in deciding whether to use multifidelity surrogates},
  author={Giselle Fern{\'a}ndez-Godino, M and Park, Chanyoung and Kim, Nam H and Haftka, Raphael T},
  journal={AIAA Journal},
  volume={57},
  number={5},
  pages={2039--2054},
  year={2019},
  publisher={American Institute of Aeronautics and Astronautics},
  doi = {10.2514/1.J057750}
}

@inproceedings{beran2020comparison,
  title={Comparison of multi-fidelity approaches for military vehicle design},
  author={Beran, Philip S and Bryson, Dean and Thelen, Andrew S and Diez, Matteo and Serani, Andrea},
  booktitle={AIAA AVIATION 2020 FORUM},
  pages={3158},
  year={2020},
  doi = {10.2514/6.2020-3158}
}

@article{forrester2007-PRSA,
  title={Multi-fidelity optimization via surrogate modelling},
  author={Forrester, Alexander IJ and S{\'o}bester, Andr{\'a}s and Keane, Andy J},
  journal={Proceedings of the Royal Society {A}: mathematical, Physical and Engineering sciences},
  volume={463},
  number={2088},
  pages={3251--3269},
  year={2007},
  publisher={The Royal Society London},
  doi = {10.1098/rspa.2007.1900}
}

@article{DiscForr2018,
  title={Linking {G}aussian process regression with data-driven manifold embeddings for nonlinear data fusion},
  author={Lee, Seungjoon and Dietrich, Felix and Karniadakis, George E and Kevrekidis, Ioannis G},
  journal={Interface Focus},
  volume={9},
  number={3},
  pages={20180083},
  year={2019},
  publisher={The Royal Society},
  doi = {10.1098/rsfs.2018.0083}
}

@article{JumpForr,
  author    = {Maziar Raissi and
               George E. Karniadakis},
  title     = {Deep Multi-fidelity Gaussian Processes},
  journal   = {CoRR},
  year      = {2016},
  doi = {10.48550/arXiv.1604.07484}
}

@inproceedings{ficini2021-AIAA,
  title={An Adaptive N-Fidelity Metamodel for Design Optimization via Gaussian Process Regression},
  author={Ficini, S and Iemma, U and Pellegrini, R and Serani, A and Diez, M},
  booktitle={AIAA AVIATION 2021 Forum},
  pages={},
  year={2021},
  doi = {}
}

@article{aguilera2004general,
  title={General n-dimensional rotations},
  author={Aguilera, Antonio and P{\'e}rez-Aguila, Ricardo},
  year={2004},
  publisher={V{\'a}clav Skala-UNION Agency},
}

@article{clark2016-AIAA,
  title={Engineering design exploration using locally optimized covariance Kriging},
  author={Clark Jr, Daniel L and Bae, Ha-Rok and Gobal, Koorosh and Penmetsa, Ravi},
  journal={AIAA Journal},
  volume={54},
  number={10},
  pages={3160--3175},
  year={2016},
  publisher={American Institute of Aeronautics and Astronautics},
  doi = {10.2514/1.J054860}
}

@article{Toal2014,
  year = {2014},
  month = {dec},
  publisher = {Springer Science and Business Media {LLC}},
  volume = {51},
  number = {6},
  pages = {1223--1245},
  author = {David J. J. Toal},
  title = {Some considerations regarding the use of multi-fidelity Kriging in the construction of surrogate models},
  journal = {Structural and Multidisciplinary Optimization},
  doi = {10.1007/s00158-014-1209-5}
}

@inproceedings{quagliarella2004aerodynamic,
  title={Aerodynamic shape design using hybrid evolutionary computation and fitness approximation},
  author={Quagliarella, Domenico and Iannelli, Pierluigi and Vitagliano, Pier Luigi and Chinnici, Giorgio},
  booktitle={AIAA 1st Intelligent Systems Technical Conference},
  pages={6514},
  year={2004},
  doi = {10.2514/6.2004-6514}
}

@InProceedings{yildiz2022-ECCOMAS,
  Title = {Advanced Experiments on Gaussian Process-based Multi-fidelity Methods over Diverse Mathematical Characteristics},
  Author                   = {Yildiz, S. and Pehlivan-Solak, H. and Diez, M. and Goren, O. and Nikbay, M.},
  Booktitle                = {8th European Congress on Computational Methods in Applied Sciences and Engineering, ECCOMAS},
  Year                     = {2022},
  Address                  = {Oslo, Norway},
  Month                    = {jun},
  doi = {10.23967/eccomas.2022.039}
}

@article{ruan2020variable,
  title={Variable-fidelity probability of improvement method for efficient global optimization of expensive black-box problems},
  author={Ruan, Xiongfeng and Jiang, Ping and Zhou, Qi and Hu, Jiexiang and Shu, Leshi},
  journal={Structural and Multidisciplinary Optimization},
  volume={62},
  number={6},
  pages={3021--3052},
  year={2020},
  publisher={Springer},
  doi = {10.1007/s00158-020-02646-9}
}

@article{huang2006sequential,
  title={Sequential kriging optimization using multiple-fidelity evaluations},
  author={Huang, Deng and Allen, Theodore T and Notz, William I and Miller, R Allen},
  journal={Structural and Multidisciplinary Optimization},
  volume={32},
  number={5},
  pages={369--382},
  year={2006},
  publisher={Springer},
  doi = {10.1007/s00158-005-0587-0},
}

@inproceedings{takeno2020multi,
  title={Multi-fidelity Bayesian optimization with max-value entropy search and its parallelization},
  author={Takeno, Shion and Fukuoka, Hitoshi and Tsukada, Yuhki and Koyama, Toshiyuki and Shiga, Motoki and Takeuchi, Ichiro and Karasuyama, Masayuki},
  booktitle={International Conference on Machine Learning},
  pages={9334--9345},
  year={2020},
  doi = {},
}

@inproceedings{emukit2019,
  author = {Paleyes, Andrei and Pullin, Mark and Mahsereci, Maren and Lawrence, Neil and González, Javier},
  title = {Emulation of physical processes with Emukit},
  booktitle = {Second Workshop on Machine Learning and the Physical Sciences, NeurIPS},
  year = {2019},
  doi = {10.48550/arXiv.2110.13293}
}

@techreport{dakota2021,
  title={Dakota, A Multilevel Parallel Object-Oriented Framework for Design Optimization, Parameter Estimation, Uncertainty Quantification, and Sensitivity Analysis: Version 6.15 User’s Manual},
  author={Adams, B. M. and Bohnhoff, W. J. and Dalbey, K. R. and Ebeida, M. S. and Eddy, J. P. and Eldred, M. S. and Hooper, R. W. and Hough, P. D. and Hu, K. T. and Jakeman, J. D. and Khalil, M. and Maupin, K. A. and Monschke, J. A. and Ridgway, E. M. and Rushdi, A. A. and Seidl, D. T. and Stephens, J. A. and Swiler, L. P. and Winokur, J. G.},
  group = {csg},
  year = {2021},
  institution = {Sandia},
  month = nov,
  number = {SAND2020-12495},
  doi = {10.2172/1829573}
}

@misc{Rumpfkeil2016,
  AUTHOR = "M. P. Rumpfkeil and P. Beran",
  TITLE = "{Construction of Multi-Fidelity Surrogate Models for Aerodynamic Databases}",
  howpublished = "{Proceedings of the Ninth International Conference on Computational Fluid Dynamics, 
  ICCFD9, Istanbul, Turkey, July 11-15}",
  year = 2016,
  doi = {}
}

@article{Minisci:2013,
  title = {Robust Design of a Re-entry Unmanned Space Vehicle by Multi-fidelity Evolution Control},
  author = {Minisci, E. and Vasile, M.},
  journal = {AIAA Journal},
  volume = {51},
  number = {6},
  year = {2013},
  pages = {1284--1295},
  doi = {10.1145/2001576.2001671}
}

@article{RCKthelen2020,
  title = {Multifidelity Flutter Prediction using Regression Cokriging with Adaptive Sampling},
  author = {Thelen, A. and Leifsson, L. and Beran, P.},
  journal = {Journal of {F}luids and {S}tructures},
  volume = {97},
  pages = {103081},
  year = {2020},
  doi = {10.1016/j.jfluidstructs.2020.103081}
}

@article{wackers2022efficient,
  title={Efficient initialization for multi-fidelity surrogate-based optimization},
  author={Wackers, Jeroen and Pellegrini, Riccardo and Serani, Andrea and Visonneau, Michel and Diez, Matteo},
  journal={Journal of {O}cean {E}ngineering and {M}arine {E}nergy},
  pages={1--17},
  year={2022},
  publisher={Springer},
  doi = {10.1007/s40722-022-00268-5}
}

@article{Bryson2017b,
  Author = {Bryson, Dean E. and Rumpfkeil, Markus P.},
  Journal = {Aerospace Science and Technology},
  Pages = {121--136},
  Title = {All-at-once approach to multifidelity polynomial chaos expansion surrogate modeling},
  Volume = {70},
  Year = {2017},
  doi = {10.1016/j.ast.2017.07.043}
}

@article{Rumpfkeil2020b,
  AUTHOR = "M. P. Rumpfkeil and P. Beran",
  TITLE = "{Multi-Fidelity Sparse Polynomial Chaos Surrogate Models Applied to Flutter Databases}",
  journal = "{AIAA} Journal",
  volume = "58, No.3",
  pages = "1292-1303",
  year = 2020,
  doi = {10.2514/1.J058452}
}

@article{Rumpfkeil2022a,
  AUTHOR = "M.P. Rumpfkeil and D. Bryson and P. Beran",
  TITLE = "{Multi-Fidelity Sparse Polynomial Chaos and Kriging Surrogate Models Applied to Analytical Benchmark Problems}",
  journal = "Algorithms",
  volume = "15(3)",
  pages = "101",
  year = 2022,
  doi = {10.3390/a15030101}
}

@article{serani2016-ASOC,
  title={Parameter selection in synchronous and asynchronous deterministic particle swarm optimization for ship hydrodynamics problems},
  author={Serani, Andrea and Leotardi, Cecilia and Iemma, Umberto and Campana, Emilio F and Fasano, Giovanni and Diez, Matteo},
  journal={Applied {S}oft {C}omputing},
  volume={49},
  pages={313--334},
  year={2016},
  publisher={Elsevier},
  doi = {10.1016/j.asoc.2016.08.028}
}

\end{document}